\documentclass{article}
\usepackage{amsmath}
\usepackage{amssymb}
\usepackage{ntheorem}
\usepackage{mathdots}
\usepackage[super]{nth}
\usepackage{graphicx}
\usepackage{blkarray}
\usepackage{multirow}
\usepackage{mathtools}
\usepackage{subcaption}

\title{Projections, embeddings and stability}
\author{Pelle Olsson}
\date{\today}

\newtheorem{thm}{Theorem}
\newtheorem{prop}[thm]{Proposition}
\newtheorem{cor}[thm]{Corollary}
\newtheorem{lemma}[thm]{Lemma}
\newtheorem{define}[thm]{Definition}

{\theorembodyfont{\rmfamily} \newtheorem{exa}[thm]{Example}}
{\theorembodyfont{\rmfamily} \newtheorem{rem}[thm]{Remark}}

\newenvironment{acknowledgements}{%
  \begin{abstract}
}{%
  \end{abstract}
}

\begin{document}

\maketitle

\begin{abstract}
In the present work, we demonstrate how the pseudoinverse concept from linear algebra can be used to represent
and analyze the boundary conditions of linear systems of partial differential equations.  This approach 
has theoretical and practical implications; the theory applies even if the boundary operator is rank deficient,
or near rank deficient.  If desired, the pseudoinverse can be implemented directly using standard tools
like Matlab.  We also introduce a new and simplified version of the semidiscrete approximation of the linear 
PDE system, which completely avoids taking the time derivative of the boundary data, cf.~\cite{po:spps2}.
The stability results of \cite{po:spps1} are generalized to nondiagonal summation-by-parts norms.
Another key result is the extension of summation-by-parts operators to multi-domains by means of carefully
crafted embedding operators.  No extra numerical boundary conditions are required at the grid interfaces.
The aforementioned pseudoinverse allows for a compact representation of these multi-block operators,
which preserves all relevant properties of the single-block operators.  The embedding operators can be constructed
for multiple space dimensions.  Numerical results for the two-dimensional Maxwell's equations are presented, 
and they show very good agreement with theory.
\end{abstract}

\begin{acknowledgements}
The author is indebted to Prof.~Ken Kreutz-Delgado, Department of Electrical and Computer Engineering, 
UC San Diego, CA, for sharing his lecture notes for ECE~174, which offer the theoretical framework for 
implementing boundary conditions as pseudoinverses.  The numerical results presented in Section~\ref{sec:numerics}
are provided courtesy of MSc.~Gustav Eriksson, Department of Scientific Computing, Uppsala University,
Sweden.
\end{acknowledgements}

\section{Introduction}
\label{sec:intro}
The focus of the present study is summation-by-parts (SBP) methods for the model problem
\begin{align}
u_t(x,t) + Qu(x,t) &= 0, \quad t > 0, \ x\in\Omega \label{eq:model} \\
Lu(x,t) &= g(t), \quad t \geq 0, \ x\in\Gamma \nonumber \\
u(x,0) &= f(x). \nonumber
\end{align}
We will restrict ourselves to the case where $\Omega$ is a subset in $\mathbb{R}$ or $\mathbb{R}^2$;
$\Gamma$ refers to the boundary of $\Omega$.  The differential operator $Q=Q(\partial)$ is assumed 
to be semibounded, cf.~\cite{kl:ibvpnse}, leading to a well-posed problem in the sense that any 
solution of (\ref{eq:model}) must satisfy an energy estimate
\[
\|u(\cdot,t)\|^2 \leq Ke^{ct}\left(\|f\|^2 + \int_0^t\|g(\cdot,\tau)\|^2_\Gamma d\tau\right),
\]
where $\|\cdot\|$ and $\|\cdot\|_\Gamma$ are the $L^2$-norms on $\Omega$ and $\Gamma$.  This kind
of estimate typically follows from an integration-by-parts procedure (divergence theorem in multiple
space dimensions) and a properly designed boundary condition operator $L$.

We have adopted an operator centric approach for analyzing summation-by-parts 
difference methods and their boundary conditions, which are implemented by means of
projections.  This technique is also used to define multi-block difference operators.  The analysis 
is based on two key concepts: {\em pseudoinverses} and {\em embedding} operators.  The latter are used to define
multi-block difference operators that satisfy summation by parts given the existence of
SBP operators for the individual blocks.  There is no need to construct "extra" boundary conditions
at the grid interfaces;  the embedding operators will handle this automatically.

Rather than considering the state spaces as some space $\mathbb{R}^m$ equipped with a special scalar product 
$(\cdot,\cdot)_H$, we will regard the pair $[\mathbb{R}^m, (\cdot,\cdot)_H]$ as an inner product space
$V$ in its own right; the scalar product $(\cdot,\cdot) \equiv (\cdot,\cdot)_H$ will {\em define}
$V$.  In this context, difference operators, boundary operators, projections
and embeddings will be treated as mappings between inner products spaces with well-defined adjoints 
and pseudoinverses.  This will lead to a systematic and concise notation for multi-block operators.
Since the theory relies on well-defined concepts like adjoints and pseudoinverses, the resulting
discretizations can be implemented directly using matrix algebra, e.~g. Matlab or similar packages.

In what follows, $x, y \in \mathbb{R}^m$ and $z, w \in \mathbb{R}^n$ will be viewed as vectors 
in some linear spaces $V$ and $W$ with scalar products
\begin{equation}
\label{eq:scal}
(x,y) \equiv x^TH_1y, \quad \langle z,w\rangle \equiv z^TH_2w,
\end{equation}
where we have dropped the usual subscripts ${}_{H_i}$ of the scalar products.  
A linear operator $T$ is a mapping
\begin{equation}
\label{eq:opdef}
\nonumber
T:V \rightarrow W
\end{equation}
between two inner product spaces (possibly identical).

Sections~\ref{sec:proj},~\ref{sec:adjoint} briefly review some basic properties of projections
and adjoint operators.  They can be ignored by readers who are familiar with these concepts.
Section~\ref{sec:pseudo} gives a detailed account of the theory of pseudoinverses,
which typically is presented in the context of the Euclidean scalar product
\[
(x,y) \equiv x^Ty \in \mathbb{R}, \qquad x, y \in \mathbb{R}^{m},
\]
e.~g., \cite{rp:gifm, aa:rmpp, bic:ctgi}.  To obtain difference operators 
$D\in\mathbb{R}^{m\times m}$ that fulfil summation by parts, it is necessary to work with weighted scalar products:
\[
(x,y) \equiv x^THy \in \mathbb{R}, \qquad H \in \mathbb{R}^{m \times m},
\]
where $H$ is symmetric positive definite (SPD) \cite{ks:fefdhpde}.  Combining the
theory in \cite{aa:rmpp} with the operator centric approach of \cite{kkd:ece174}
leads to a formulation of the Moore-Penrose conditions \cite{ehm:rgam,rp:gifm} that is 
particularly useful in the subsequent stability theory, where they will play a fundamental role.

In Section~\ref{sec:sbp}, we have collected some well-known results about summation-by-parts
operators.  We also introduce permutations matrices $J_r$ that will prove useful in subsequent sections.  We then
set the scene by defining the basic inner product space that will be used again and again
throughout the presentation.  The next section is concerned with the semidiscrete version of
the model problem \eqref{eq:model}.  The key result is a simplified form, which significantly reduces
the implementation complexity of the corresponding method found in \cite{po:spps1}.

The main objective of Section~\ref{sec:bc} is to establish the pseudoinverse as a tool
for constructing the boundary projection.  Numerous examples illustrate the general theory.  
In Section~\ref{sec:multiblock} we turn our interest to multi-block theory.  The concept of 
a {\em multiset} turns out to be very helpful in this context.  Preparing for multi-block scalar
products and difference operators, we introduce so-called {\em augmented} state spaces and 
embedding operators.  This is followed by a detailed analysis of the resulting multi-block
operators.

Sections~\ref{sec:2dim} and~\ref{sec:2dimmb} are devoted to the extension of the previous results to two 
space dimensions.  The results presented in the one-dimensional case carry over to the two-dimensional 
case verbatim.  The technique used in two space dimensions can be extended to higher-dimensional
spaces.  We also obtain a generalization of the main stability result in \cite{po:spps1}.

Finally, in Section~\ref{sec:numerics} we tie all loose ends together by applying the theory to
the two-dimensional Maxwell's equations on a curvilinear four-block domain.  Numerical results
confirming the expected convergence rate are presented.

\section{Projections}
\label{sec:proj}
This brief section lists two fundamental results of inner product spaces.  Without proof we state
the well-known projection theorem:

\begin{thm}
\label{thm:proj}
Let $x \in V$ be an inner product space (\ref{eq:scal}) and let ${\cal L}$ be a linear manifold in $V$.
There is a unique vector $\hat{x}\in{\cal L}$ such that
\begin{equation}
\label{eq:proj}
\nonumber
(x-\hat{x}, y) = 0, \quad y \in {\cal L}.
\end{equation}
\end{thm}

\begin{rem}
\label{rem:proj}
The vector $\hat{x}$ is known as the projection of $x$ onto ${\cal L}$.  This equation defines
the projection of $x$ onto ${\cal L}$.  In words: The projection of $x$ onto ${\cal L}$ is a vector
$\hat{x}\in {\cal L}$ such that $x - \hat{x}$ is orthogonal to all $y \in {\cal L}$.  It is important to note that 
orthogonality is always expressed with respect to the inner product defined for the particular the 
vector space.  In our case this will invariably be a scalar product different from the standard Euclidean 
inner product. \hfill$\Box$
\end{rem}

\begin{thm}
\label{thm:mindist}
Let $x \in V$ be a vector and let ${\cal L}$ be a linear manifold in $V$.
If $x = \hat{x} + \tilde{x}$ where $\hat{x}$ is the projection of $x$ onto ${\cal L}$, then
\begin{equation}
\label{eq:mindist}
\nonumber
\| x - y \| > \| x -\hat{x} \|, \quad y \in {\cal L}
\end{equation}
iff $y \neq \hat{x}$.
\end{thm}

\noindent
{\bf Proof:} By construction, $(x - \hat{x}, y) = 0$ for all $y \in {\cal L}$.  Furthermore, 
$\hat{x} - y \in {\cal L}$.  Thus,
\[
\| x - y \|^2 = \|(\hat{x} - y) + x - \hat{x} \|^2 = \| \hat{x} - y \|^2 + \| x - \hat{x} \|^2
\ge \| x - \hat{x} \|^2_H
\]
for all $y \in {\cal L}$ with strict inequality iff $y \neq \hat{x}$. \hfill$\Box$

\section{Adjoint operators}
\label{sec:adjoint}
We will recast the relevant results of \cite{aa:rmpp} so as to apply to weighted scalar products.

\begin{define}
\label{def:adjoint}
Let $T:V \rightarrow W$. The adjoint operator $T^*:W \rightarrow V$ is defined as
\[
(x, T^*z) \equiv \langle Tx, z \rangle.
\]
\hfill$\Box$
\end{define}

\begin{prop}
\label{prop:adjoint}
Let $T:V \rightarrow W$. The adjoint operator $T^*:W \rightarrow V$ is unique and satisfies $T^{**}=T$.
\end{prop}

\noindent
{\bf Proof:} The proof will be broken down into four simple steps:  existence, linearity, uniqueness 
and reflexivity.

{\em Existence}:  
From the definition of the adjoint operator and the inner product~(\ref{eq:scal}):
\[
x^TH_1T^*z \quad \Longleftrightarrow \quad H_1T^* = T^TH_2
\]
Thus, 
\begin{equation}
\label{eq:adjtransp}
T^* = H_1^{-1}T^TH_2. 
\end{equation}

{\em Linearity}: Follows immediately from the definition of $T^*$ and the linearity of $T$.  This 
proves existence of the adjoint operator $T^*:W \rightarrow V$.

{\em Uniqueness}: Suppose that there are two operators 
$T^*_1, T^*_2: W \rightarrow V$ that fulfill Definition~\ref{def:adjoint}.
Hence,
\[
(x,T_1^*z) = (x,T_2^*z),
\]
which proves uniqueness (repeat the arguments from the existence proof).

{\em Reflexivity}: Let $S = T^*$.  Thus,
\[
\langle z, Tx \rangle = \langle Tx, z \rangle = (x, T^*z) = (x, Sz) = (Sz, x) = 
\langle z, S^*x \rangle = \langle z, T^{**}x \rangle,
\]
which is true for all $x \in V, z \in W$. This concludes the proof. \hfill$\Box$

\begin{prop}
\label{prop:inv}
Let $T:V \rightarrow V$ be invertible.  The inverse of $T^*$ exists and satisfies 
\[
\left(T^*\right)^{-1} = \left(T^{-1}\right)^* 
\]
\end{prop}

\noindent
{\bf Proof:} Since $T$ is invertible, the adjoint of $T^{-1}$ is well defined:
\[
(x, \left(T^{-1}\right)^*y) = (T^{-1}x, y).
\]
It follows that
\[
(x, T^*\left(T^{-1}\right)^*y) = (Tx, \left(T^{-1}\right)^*y) = (T^{-1}Tx, y) = (x, y).
\]
Hence, $T^*\left(T^{-1}\right)^* = I$.  Similarly,
\[
(x, \left(T^{-1}\right)^*T^*y) = (x, y).
\]
In summary,
\[
\left(T^{-1}\right)^*T^* = T^*\left(T^{-1}\right)^* = I,
\]
which proves that $\left(T^*\right)^{-1} = \left(T^{-1}\right)^*$ \hfill$\Box$
\vskip 0.2 cm
The following propositions are presented without proof:

\begin{prop}
\label{prop:adjadd}
Let $S,T:V \rightarrow W$.  Then $(S + T)^* = S^* + T^*$.
\end{prop}

\begin{prop}
\label{prop:adjmul}
Let $S:U \rightarrow V$ and $T:V \rightarrow W$.  
Then $(TS)^* = S^*T^*$.
\end{prop}

Before proceeding a few more definitions will be required.

\begin{define}
\label{def:nr}
Let $T:V \rightarrow W$.  Define four linear manifolds:
\begin{alignat*}{1}
{\cal N}(T)       & = \{ x \in V : Tx = 0 \} \\
{\cal R}(T)       & = \{ z \in W : z = Tx \ \mbox{for some} \ x \in V \} \\
{\cal N}(T)^\perp & = \{ x \in V : (x,y) = 0, \quad y \in {\cal N}(T) \} \\
{\cal R}(T)^\perp & = \{ z \in W : \langle z,w \rangle = 0, \quad w \in {\cal R}(T) \}.
\end{alignat*}
${\cal N}(T)$ is known as the null space of $T$ and ${\cal R}(T)$ is the range of $T$.
\hfill$\Box$
\end{define}

The three subsequent propositions are very important and will be used frequently when deriving 
the pseudoinverse of $T$.

\begin{prop}
\label{prop:nr1}
Let $T:V \rightarrow W$.  Then
\begin{alignat}{1}
{\cal N}(T^*) & = {\cal R}(T)^\perp     \label{eq:nr11} \\
{\cal R}(T^*) & = {\cal N}(T)^\perp     \label{eq:nr12} \\
{\cal N}(T)   & = {\cal R}(T^*)^\perp   \label{eq:nr13} \\
{\cal R}(T)   & = {\cal N}(T^*)^\perp.  \label{eq:nr14}
\end{alignat}
\end{prop}

\noindent
{\bf Proof:} If we can prove (\ref{eq:nr11}) and (\ref{eq:nr12}), then (\ref{eq:nr13}) and (\ref{eq:nr14}) 
follow by replacing  $T \rightarrow T^*$ and using $T^{**} = T$.  

Suppose that $z \in {\cal N}(T^*)$,
i.e.,  $T^*z = 0$.  Thus,
\[
0=(x, T^*z) = \langle Tx,z \rangle,
\]
for all $x\in V$, which shows that $z \in {\cal R}(T)^\perp$, that is, ${\cal N}(T^*) \subset {\cal R}(T)^\perp$.

Conversely, if $z \in {\cal R}(T)^\perp$, then
\[
0 = \langle Tx,z \rangle = (x, T^*z)
\]
for all $x\in V$.  Choose $x = T^*z$, which implies
\[
 (T^*z, T^*z) = 0 \quad \Longleftrightarrow \quad T^*z = 0 \quad \Longleftrightarrow \quad z \in {\cal N}(T^*),
\]
since $(\cdot, \cdot)$ is positive definite.  Hence, ${\cal R}(T)^\perp \subset {\cal N}(T^*)$, which 
establishes (\ref{eq:nr11})

Now, (\ref{eq:nr12}) follows by forming the orthogonal complement of (\ref{eq:nr11}) and 
replacing $T \rightarrow T^*$:
\[
{\cal R}(T^*) = {\cal N}(T^{**}) = [\mbox{Prop. \ref{prop:adjoint}}] = {\cal N}(T),
\]
which concludes the proof. \hfill$\Box$

\begin{prop}
\label{prop:decomp}
Let $T:V \rightarrow W$ and $z$ a vector in $W$.  Then $z$ can 
be uniquely decomposed as
\begin{equation}
\label{eq:decomp}
\nonumber
z = \hat{z} + \tilde{z},
\end{equation}
where $\hat{z}$ is the orthogonal projection of $z$ 
onto ${\cal R}(T)$ and $\tilde{z} \in {\cal N}(T^*)$.
\end{prop}

\noindent
{\bf Proof:} Apply Theorem \ref{thm:proj} with ${\cal L} = {\cal R}(T)$.  Hence, there is a unique 
vector $\hat{z}\in {\cal R}(T)$ that satisfies
\[
\langle z - \hat{z}, w \rangle = 0, \quad w \in {\cal R}(T).
\]
This means that
\[
\tilde{z} \equiv z - \hat{z} \in {\cal R}(T)^\perp = [\mbox{Prop.~\ref{prop:nr1}}] = {\cal N}(T^*),
\]
which proves the claim. \hfill$\Box$

\begin{prop}
\label{prop:nr2}
Let $T:V \rightarrow W$.  Then
\begin{alignat}{1}
{\cal R}(T)   & = {\cal R}(TT^*)   \label{eq:nr21} \\
{\cal R}(T^*) & = {\cal R}(T^*T)   \label{eq:nr22} \\
{\cal N}(T)   & = {\cal N}(T^*T)   \label{eq:nr23} \\
{\cal N}(T^*) & = {\cal N}(TT^*).  \label{eq:nr24}
\end{alignat}
\end{prop}

\noindent
{\bf Proof:} Let $x \in {\cal N}(T)$.  Then
\[
 Tx = 0 \quad \Longrightarrow \quad T^*Tx = 0 \quad \Longleftrightarrow \quad x \in {\cal N}(T^*T).
\]
Hence ${\cal N}(T) \subset {\cal N}(T^*T)$.  Conversely, suppose $x \in {\cal N}(T^*T)$.  Then
\[
T^*Tx = 0 \quad \Longrightarrow \quad (x,T^*Tx) = 0.
\]
From the definition of $T^*$ it follows that
\[
\langle Tx,Tx\rangle = 0 \quad \Longleftrightarrow \quad Tx = 0 \quad \Longleftrightarrow \quad x \in {\cal N}(T).
\] 
This shows that ${\cal N}(T^*T) \subset {\cal N}(T)$, which proves (\ref{eq:nr23}).  To prove 
(\ref{eq:nr24}) it suffices to substitute $T \rightarrow T^*$ in (\ref{eq:nr23}) and then use 
$T^{**} = T$.  Applying Proposition~\ref{prop:nr1} twice:
\[
{\cal R}(T) = {\cal N}(T^*)^\perp = [(\ref{eq:nr24})] = {\cal N}(TT^*)^\perp = {\cal R}([TT^*]^*) = {\cal R}(TT^*),
\]
which proves (\ref{eq:nr21}). Eq.~(\ref{eq:nr22}), finally, follows by substituting $T\rightarrow T^*$ 
in the above expression. \hfill$\Box$

\begin{prop}
\label{prop:onto}
Let $T:V \rightarrow W$.  Then $T$ is onto iff $T^*$ is one-to-one.
\end{prop}

\noindent
{\bf Proof:}
Suppose that $T$ is onto.  Let $z \in {\cal N}(T^*)$.  Then 
\[
T^*z = 0 \quad \Longleftrightarrow \quad (x, T^*z) = 0  
         \quad \Longleftrightarrow \quad \langle Tx, z \rangle = 0
\]
for all $x\in v$.  But $T$ onto means that we can choose $x = x_0$ such that $z = Tx_0$.  Thus,
\[
\| z \|^2 = \langle z, z \rangle = 0 \quad \Longleftrightarrow \quad z = 0,
\]
i.e., ${\cal N}(T^*) = \{ 0 \}$, which shows that $T^*$ is one-to-one.

Conversely, suppose that $T^*$ is one-to-one.  Let $z$ 
be an arbitrary vector in $\mathbb{R}^n$.  By Proposition~\ref{prop:decomp}:
\[
z = \hat{z} + \tilde{z}, \quad \hat{z} \mbox{ is the projection of $z$ onto } {\cal R}(T),
    \quad \tilde{z} \in {\cal N}(T^*).
\]
Now $T^*$ is one-to-one, then by definition:
\[
{\cal N}(T^*) = \{ 0 \}  \quad \Longrightarrow \quad \tilde{z} = 0.
\]
Hence,
\[ 
z = \hat{z}.
\]
Since $\hat{z}$ is the projection of $z$ onto ${\cal R}(T)$, there is a vector 
$x \in \mathbb{R}^m$ such that
\[
z = \hat{z} = Tx,
\]
where $z \in W$ is arbitrary.  This shows that $T$ is onto.
\hfill$\Box$

\subsection{Self-adjoint operators}
The core of the existence proof of the pseudoinverse depends on the spectral theorem for
self-adjoint operators.  Let us begin by recalling the definition of a self-adjoint operator.

\begin{define}
\label{def:selfadj}
An operator $T:V \rightarrow V$ is self-adjoint if $T^* = T$. \hfill$\Box$
\end{define}

\begin{rem}
\label{rem:selfadj}
Self-adjoint operators are only defined for operators that map an inner product space onto itself.
This implies that there is only one scalar product involved in the definition and the criterion
becomes
\begin{equation}
\label{eq:selfadj}
\nonumber
(Tx, y) = (x, Ty), \quad x, y \in V.
\end{equation} \hfill$\Box$
\end{rem}

The spectral theorem of self-adjoint operators can be formulated as (proof omitted):

\begin{thm}
\label{thm:spectral}
The operator $T:V \rightarrow V$ is self-adjoint iff its eigenvalues 
$\lambda_j$ are real and the corresponding eigenvectors $e_j$ are orthonormal:
\begin{equation}
\label{eq:spectral}
\nonumber
Te_j = \lambda_je_j, \quad (e_i, e_j) = \delta_{ij}, \quad j = 1, \ldots, m.
\end{equation}
\end{thm}

When analyzing the pseudoinverse, we will frequently encounter operators of the 
form $T^*T$ and $T^*T + \delta^2I$.  We conclude this section by collecting two simple
results.

\begin{prop}
\label{prop:sa1}
Let $T:V \rightarrow W$.  Then $T^*T:V \rightarrow V$ and $TT^*:W \rightarrow W$ are self-adjoint.
\end{prop}

\noindent
{\bf Proof:} Follows immediately from Definition~\ref{def:adjoint}. \hfill$\Box$

\begin{prop}
\label{prop:sa2}
Let $T:V \rightarrow W$.  Then $T^*T + \delta^2I$
and $TT^* + \delta^2I$ are self-adjoint and nonsingular for $\delta \neq 0$.
\end{prop}

\noindent
{\bf Proof:}
Self-adjointness follows immediately from Propositions~\ref{prop:sa1} and~\ref{prop:adjadd}.

To prove nonsingularity, suppose that $(T^*T + \delta^2I)x = 0$ for some $x \in V$.  Thus,
\[
0 = (x, (T^*T + \delta^2I)x) = (x, T^*Tx) + \delta^2(x,x) = \langle Tx, Tx \rangle + \delta^2(x,x)
\]
Since scalar products are positive definite, the above expression is true only if
\[
\langle Tx, Tx \rangle = \delta^2(x,x) = 0.
\]
From the last equality it follows immediately that $(x,x) = 0 \Longleftrightarrow x = 0$ whenever 
$\delta \neq 0$.  This shows that $T^*T + \delta^2I$ is non-singular.  The case $TT^* + \delta^2I$
is handled analogously. \hfill$\Box$

\begin{prop}
\label{prop:inv1}
Let $T:V \rightarrow W$. If $T$ is onto, then  $TT^*$ is invertible.  
Similarly, if $T$ is one-to-one, then $T^*T$ is invertible.
\end{prop}

\noindent
{\bf Proof:} Suppose that $T$ is onto.  From Proposition~\ref{prop:nr2} it follows that $TT^*$ is onto.
Hence, $(TT^*)^*$ is one-to-one by Proposition~\ref{prop:onto}.  But
\[
TT^* = [\mbox{Prop.~\ref{prop:sa1}}] = (TT^*)^*.
\]
Thus, $TT^*$ is both onto and one-to-one, which shows that $TT^*$ is invertible.  

If $T$ is one-to-one, then $T^*T$ is one-to-one as well (Proposition~\ref{prop:nr2}).
But $T^*T$ is the adjoint of $(T^*T)^*$.  Applying Proposition~\ref{prop:onto} to $(T^*T)^*$, it 
follows that $T^*T = (T^*T)^*$ is onto, i.~e., $T^*T$ is invertible. \hfill$\Box$

\section{Least squares and the pseudoinverse}
\label{sec:pseudo}
The theory of pseudoinverses harkens back to the 1920s with the pioneering work of Moore
\cite{ehm:rgam}.  It was later picked up by Bjerhammar \cite{ab:rrmsrgc} and
Penrose \cite{rp:gifm}.  During the 1960s, the general theory underwent rapid development,
e.~g., \cite{tg:sapm,tg:gimp,rc:rgipm,rc:ngipm}.  The remainder of this section is devoted
to the derivation of the pseudoinverse in a form that is suitable for stability analysis
when implementing boundary conditions of partial differential equations \eqref{eq:model}.

As mentioned in the introduction, we will derive the pseudoinverse by applying the pattern
set forth in \cite{kkd:ece174} to the approach used in \cite{aa:rmpp}.  This will extend the 
results to operators in inner product spaces as opposed to matrices and vectors in $\mathbb{R}^m$ using the standard Euclidean product $(x, y) = x^Ty$.
We will use the same symbol for operators and their corresponding matrix representation, unless
clarity demands that different notations be used.  Normally it should be clear from the context
what a particular symbol designates:
\[
\begin{array}{c}
T^*: \mbox{operator} \\
T^T: \mbox{matrix}.
\end{array}
\]
The relation between an adjoint operator and its matrix representation is given by (\ref{eq:adjtransp}).
To reduce clutter, we will use the same notation for norms in $V$ and $W$:
\begin{alignat*}{2}
\| x \|^2 &\equiv (x, x),               & \quad x \in V, \\
\| z \|^2 &\equiv \langle z, z \rangle, & \quad z \in W.
\end{alignat*}
This should cause no confusion, since we have adopted the convention that $x, y \in V$ and 
$z, w \in W$.

The main result of this section is a generalization of Penrose's characterization of the pseudoinverse
\cite{rp:gifm}, which will be needed in the subsequent sections when dealing with stable boundary 
conditions for semidiscrete approximations of PDEs.

We begin by establishing three equivalent characterizations of the (unique) solution to a least
squares problem (Theorems \ref{thm:ls1}, \ref{thm:ls2}).  The lemma that follows is a technical result
that proves that the limit of certain operators always exists and that this limit defines a projection
in the sense of Theorem \ref{thm:proj}.  At this point we are ready to state Theorem \ref{thm:pseudo},
which provides an explicit expression for the pseudoinverse.  This inverse will in fact return the
least square solution of Theorems \ref{thm:ls1}, \ref{thm:ls2}.  As a corollary we obtain closed formulas
for the fundamental projections onto the manifolds ${\cal R}(T)$, ${\cal N}(T)$, ${\cal R}(T^*)$ and
${\cal N}(T^*)$. After this, there follows a digression on the spectral decomposition of the pseudoinverse.
Proposition \ref{prop:altpseudo} offers an alternative expression for pseudoinverses.  The alternate
form will be used when extending Penrose's characterization of the pseudoinverse to the case of
non-standard scalar products (Theorem~\ref{thm:penrose}).

\begin{thm}
\label{thm:ls1}
Let $T:V \rightarrow W$.  For any $z \in W$ the 
following statements are equivalent:

\begin{enumerate}
\item[(i)] There exists a unique vector $\hat{x} \in V$ of minimum norm that minimizes
\begin{equation}
\label{eq:ls1a}
\| z - Tx \|^2.
\end{equation}
\item[(ii)] There exists a unique vector $\hat{x} \in {\cal R}(T^*)$ that satisfies
\begin{equation}
\label{eq:ls1b}
Tx = \hat{z},
\end{equation}
where $\hat{z}$ is the projection of $z$ onto ${\cal R}(T)$.
\end{enumerate}
\end{thm}

\noindent{\bf Proof:} The proof will be carried out in four distinct steps: equivalence of (\ref{eq:ls1a}) 
and (\ref{eq:ls1b}) for {\em any} minimizer $x_0$, not just minimum norm minimizers; existence of minimizers;
existence of minimum norm minimizers; uniqueness of minimum norm minimizers.

{\em Equivalence}: According to Proposition \ref{prop:decomp} the vector $z$ can be decomposed as
\[
z = \hat{z} + \tilde{z},
\]
where $\hat{z}$ is the projection of $z$ onto ${\cal R}(T)$ in the sense of Theorem~\ref{thm:proj}; 
$\tilde{z} \in {\cal N}(T^*)$.
Hence,
\begin{equation}
\label{eq:ls1min}
\nonumber
\| z - Tx \|^2 = \langle \hat{z} - Tx, \hat{z} - Tx \rangle + \langle \tilde{z}, \tilde{z} \rangle \geq
\| \tilde{z} \|^2 \quad \mbox{for all}\ x \in \mathbb{R}^m. 
\end{equation}
Thus, the minimization problem (\ref{eq:ls1a}) has a lower bound $\| \tilde{z} \|^2$.  If $x_0 \in V$ 
solves (\ref{eq:ls1b}), then
\[
\| z - Tx_0 \|^2 = \| \tilde{z} \|^2,
\]
which shows that the lower bound of (\ref{eq:ls1a}) is attained for $x = x_0$ if $x_0$ solves (\ref{eq:ls1b}).
Hence, existence of a solution $\hat{z} = Tx_0$ is a sufficient condition for (\ref{eq:ls1a}) to attain
its minimum at $x_0$.

Conversely, to prove necessity we note that
\[
\| z - Tx \|^2 = \| \hat{z} - Tx \|^2 + \| \tilde{z} \|^2.
\]
If minimum is attained at $x_0$, that is, if $\| z - Tx_0 \|^2 = \| \tilde{z} \|^2$, then $x_0$ must
satisfy
\[
\| \hat{z} - Tx_0 \|^2 = 0 \quad \Longleftrightarrow \quad \hat{z} = Tx_0.
\]
We have thus shown that $x$ minimizes (\ref{eq:ls1a}) iff $x$ solves (\ref{eq:ls1b}).

{\em Existence}: At this point we have not shown that there exists a minimizer $x_0$.  But this follows 
easily by observing that $\hat{z} \in {\cal R}(T)$ implies that there must be an element $x_0 \in V$ 
such that $\hat{z} = Tx_0$.  The first part of the proof thus implies that $x_0$ is a minimizer of
(\ref{eq:ls1a}).

{\em Minimum norm minimizers}: To prove existence of minimum norm minimizers we invoke Proposition 
\ref{prop:decomp} again with $T \rightarrow T^*$, $z \rightarrow x_0$ to decompose $x_0$ as
\begin{equation}
\label{eq:decompx0}
x_0 = \hat{x}_0 + \tilde{x}_0,
\end{equation}
where $\hat{x}_0$ is the projection of $x_0$ onto ${\cal R}(T^*)$ as in Theorem~\ref{thm:proj} and 
where $\tilde{x}_0 \in {\cal N}(T)$ (recall that $T^{**} = T$ by Proposition~\ref{prop:adjoint}).  
Hence,
\[
\hat{z} = Tx_0 = T(\hat{x}_0 + \tilde{x}_0) = [\tilde{x}_0 \in {\cal N}(T)] = T\hat{x}_0.
\]
In other words, $\hat{x}_0$ is a solution of (\ref{eq:ls1b}).  Consequently, it too is a minimizer of 
(\ref{eq:ls1a}) by virtue of the first part of the proof.  According to (\ref{eq:decompx0}) any minimizer 
$x_0$ satisfies
\begin{equation}
\label{eq:normdecompx0}
\| x_0 \|^2 = \| \hat{x}_0 \|^2 + \| \tilde{x}_0 \|^2 \geq \| \hat{x}_0 \|^2.
\end{equation}
Thus, given any minimizer $x_0$, $\| x_0 \|^2$ is bounded below by $\| \hat{x}_0 \|^2$, where
$\hat{x}_0$ also minimizes (\ref{eq:ls1a}) and belongs to ${\cal R}(T^*)$.  

Conversely, if 
$\| x_0 \|^2 = \| \hat{x}_0 \|^2$, then
(\ref{eq:normdecompx0}) implies
\[
\| \tilde{x}_0 \|^2 = 0 \quad \Longleftrightarrow \quad \tilde{x}_0 = 0 
                        \quad \Longleftrightarrow \quad x_0 = \hat{x}_0,
\]
which shows that $x_0 \in {\cal R}(T^*)$.  Thus $x_0$ is a minimum norm minimizer of (\ref{eq:ls1a}) 
iff $\hat{z} = Tx_0$ and  $x_0 \in {\cal R}(T^*)$.

{\em Uniqueness of minimum norm minimizers}:  Assume that there are two vectors $x_0, x_1$ that satisfy
\begin{align*}
\hat{z} &= Tx_0, \quad x_0 \in {\cal R}(T^*), \\
\hat{z} &= Tx_1, \quad x_1 \in {\cal R}(T^*).
\end{align*}
Thus,
\[
T(x_1 - x_0) = 0 \quad \Longleftrightarrow \quad x_1 - x_0 \in {\cal N}(T) = [\mbox{Prop.~\ref{prop:nr1}}]
             = {\cal R}(T^*)^\perp.
\]
But according to the hypothesis $x_1 - x_0 \in {\cal R}(T^*)$, which means that $x_1 - x_0$ is orthogonal
to itself:
\[
(x_1 - x_0, x_1 - x_0) = 0 \quad \Longleftrightarrow \quad x_1 = x_0,
\]
which concludes the proof. \hfill$\Box$

\begin{rem}
\label{rem:ls1}
Theorem~\ref{thm:ls1} is the fundamental existence theorem that will be used to prove the existence of 
the pseudoinverse of $T$.  We will construct an operator $T^+:W \rightarrow V$
such that if $x$ defined as
\[
x \equiv T^+z,
\]
then $x$ will solve $Tx = \hat{z}$ and $x \in {\cal R}(T^*)$.  Thus, $x$ is the minimizer of (\ref{eq:ls1a}).
The operator $T^+$ is known as the pseudoinverse of $T$. \hfill$\Box$
\end{rem}

The next theorem offers a third alternative for characterizing the minimum norm solution of the least
square problem, the so-called normal equations.

\begin{thm}
\label{thm:ls2}
Let $T:V \rightarrow W$.  The vector $\hat{x} \in V$ that minimizes
\begin{equation}
\label{eq:ls2a}
\| z - Tx \|^2
\end{equation}
and has minimum norm is uniquely defined and satisfies
\begin{alignat}{2}
T^*T\hat{x} & = T^*z \label{eq:ls2b1}\\
    \hat{x} & = T^*w \label{eq:ls2b2}
\end{alignat}
for some $w \in W$.
\end{thm}

\noindent{\bf Proof:} The structure of the proof is as follows:  Existence of solution 
(\ref{eq:ls2b2}) to the normal equations (\ref{eq:ls2b1}); uniqueness of solution; final step where we show 
that the solution to (\ref{eq:ls2b1}) also satisfies $T\hat{x} = \hat{z}$, where $\hat{z}$ is 
the projection of $z$ onto ${\cal R}(T)$.  The result then follows from Theorem~\ref{thm:ls1}.

{\em Existence}: Let $z \in W$.  Obviously, $T^*z \in {\cal R}(T^*)$.  But 
\[
{\cal R}(T^*) = [\mbox{Prop.~\ref{prop:nr2}}] = {\cal R}(T^*T).
\]
Thus,
\[
T^*z = T^*Tx
\]
for some $x \in V$.  As usual, we split $x$ as
\[
x = \hat{x} + \tilde{x},
\]
where $\hat{x}$ is the projection of $x$ onto ${\cal R}(T^*T)$ and $\tilde{x} \in {\cal N}(T^*T)$, i.~e.,
$\hat{x}$ satisfies
\[
T^*T\hat{x} = T^*z.
\]
Furthermore, since $\hat{x} \in {\cal R}(T^*T)$ by construction and since ${\cal R}(T^*T) = {\cal R}(T^*)$,
there must be a $w \in W$ such that
\[
\hat{x} = T^*w.
\]
This shows that $\hat{x}$ defined by (\ref{eq:ls2b2}) solves (\ref{eq:ls2b1}) for any $z\in W$, 
which proves existence.

{\em Uniqueness}: Assume that there are two solutions 
\begin{alignat}{2}
\hat{x}_1 & = T^*w_1 \label{eq:ls2b3}\\
\hat{x}_2 & = T^*w_2 \label{eq:ls2b4}
\end{alignat}
that both satisfy (\ref{eq:ls2b1}).  This implies
\[
T^*T(\hat{x}_1 - \hat{x}_2) = 0.
\]
Thus, $\hat{x}_1 - \hat{x}_2 \in {\cal N}(T^*T) = {\cal N}(T)$.  Hence,
\[
T(\hat{x}_1 - \hat{x}_2) = 0.
\]
Substituting (\ref{eq:ls2b3}) and (\ref{eq:ls2b4}) into the above expression yields
\[
TT^*(w_1 - w_2) = 0 \quad \Longleftrightarrow \quad w_1 - w_2 \in {\cal N}(TT^*) = {\cal N}(T^*).
\]
Thus, 
\[
T^*(w_1 - w_2) = 0 \quad \Longleftrightarrow \quad \hat{x}_1 = \hat{x}_2 \quad [(\ref{eq:ls2b3}, \ref{eq:ls2b4})],
\]
which proves that the solution of (\ref{eq:ls2b1}) and (\ref{eq:ls2b2}) is unique.

{\em Final step:}  Partition $z$ as
\[
z = \hat{z} + \tilde{z}, \quad \hat{z} \in {\cal R}(T), \quad \tilde{z} \in {\cal N}(T^*).
\]
Hence, there is a vector $y = \hat{y} + \tilde{y}, \ \hat{y} \in {\cal R}(T^*), \ \tilde{y} \in {\cal N}(T)$ 
such that
\begin{equation}
\label{eq:ls2b5}
\hat{z} = Ty = T\hat{y}.
\end{equation}
Multiplying (\ref{eq:ls2b5}) by $T^*$:
\[
T^*\hat{z} = T^*T\hat{y}.
\]
But $T^*\hat{z} = T^*(\hat{z} + \tilde{z}) = T^*z$ since $\tilde{z} \in {\cal N}(T^*)$.  This implies
that $\hat{y} \in {\cal R}(T^*)$ satisfies 
\[
 T^*T\hat{y} = T^*z, \quad \hat{y} = T^{*}w_0
\]
for some $w_0\in W$, i.~e., $\hat{y}$ solves (\ref{eq:ls2b1}) and (\ref{eq:ls2b2}).  We know from 
the second part of this proof that this solution is unique, whence $\hat{y} = \hat{x}$.  By
(\ref{eq:ls2b5}):
\[
T\hat{x} = \hat{z}.
\]
Thus, by Theorem~\ref{thm:ls1}, $\hat{x} \in {\cal R}(T^*T) = {\cal R}(T^*)$ must be the minimizer 
of (\ref{eq:ls2a}). \hfill$\Box$

\begin{rem}
\label{rem:ls2}
Theorem~\ref{thm:ls2} states that the solution of $T^*Tx = T^*z$ is unique if we restrict potential
solution candidates to $x\in{\cal R}(T^*)$.  If $T^*T$ is invertible we recover the familiar expression 
for the solution of the normal  equations:
\[
\hat{x} = \left(T^*T\right)^{-1}T^*z.
\]
In particular, the above formula holds if $T$ is one-to-one, cf.~Proposition~\ref{prop:inv1}. \hfill$\Box$
\end{rem}

The following technical lemma will play a key role when establishing the pseudoinverse of 
$T:V \rightarrow W$:

\begin{lemma}
\label{lemma:proj}
Let $S:V \rightarrow V$ be self-adjoint.  The limit
\begin{alignat}{2}
P_S & \equiv \lim_{\delta \rightarrow 0}\left(S + \delta I\right)^{-1}S \label{eq:lemma1} \\
    & =      \lim_{\delta \rightarrow 0} S \left(S + \delta I\right)^{-1} \label{eq:lemma2}
\end{alignat}
exists.  Furthermore,
\begin{equation}
\label{eq:lemma3}
\hat{x} = P_Sx
\end{equation}
is the projection of $x$ onto ${\cal R}(S)$.
\end{lemma}

\noindent{\bf Proof:}
The first step is to show that the inverse of $S + \delta I$ exists for $| \delta |$ sufficiently 
small.  This is obviously a necessary condition for (\ref{eq:lemma1}) and (\ref{eq:lemma2})
to exist.  According to the spectral theorem~\ref{thm:spectral} there exists a set of mutually
orthogonal unit vectors $e_j \in V$ such that
\[
Se_j = \lambda_j e_j, \quad e_j = 1,\dots, m,
\]
where $\lambda_j \in \mathbb{R}$ is an eigenvalue of $S$.  Thus,
\begin{equation}
\label{eq:lemma4}
\nonumber
(S + \delta I)e_j = (\lambda_j + \delta)e_j, \quad e_j = 1,\dots, m.
\end{equation}
It follows that 
\begin{equation}
\label{eq:lemma5}
\nonumber
\mu_j \equiv \lambda_j + \delta \neq 0
\end{equation}
if $0 < | \delta | < \delta_0$, where $\delta_0 = \min_j(|\lambda_j|), \lambda_j$ a 
{\em non-zero} eigenvalue of $S$.  Any vector $x \in V$ can be expressed as
\[
x = \sum^m_{j=1} x_je_j, \quad x_j = (e_j, x).
\]
Thus,
\[
(S + \delta I)x = 0 \quad \Longleftrightarrow \quad \sum^m_{j=1} x_ju_je_j = 0.
\]
Scalar multiplication by $e_k$ yields
\[
x_k\mu_k = 0, \quad k = 1, \ldots, m.
\]
But we showed earlier that $\mu_k \neq 0$ for $| \delta |$ sufficently small.  Hence, 
\[
x_k = 0, \quad k = 1, \ldots, m, \quad \Longleftrightarrow \quad x = 0.
\]
This shows that $S + \delta I$ is one-to-one.  But $(S + \delta I) = (S + \delta I)^*$ and thus
$S + \delta I$ is onto as well (Prop.~\ref{prop:onto}), whence $S + \delta I$ is 
invertible for $\delta : \ 0 < | \delta | < \delta_0$.

Next, we want to show that
\begin{equation}
\label{eq:lemma6}
\left(S + \delta I\right)^{-1}S = S\left( S + \delta I\right)^{-1}.
\end{equation}
To this end we observe that $e_j$ is an eigenvector of $(S + \delta I)^{-1}$ with eigenvalue 
$1/u_j$ iff $e_j$ is an eigenvector of $S + \delta I$ with eigenvalue $u_j$ .  It follows that
\begin{alignat*}{2}
\left(S + \delta I\right)^{-1}Sx & = \sum^m_{j=1} x_j\left(S + \delta I\right)^{-1}Se_j \\
                                 & = \sum^m_{j=1} x_j\lambda_j\left(S + \delta I\right)^{-1}e_j \\
								 & = \sum^m_{j=1} x_j\lambda_j \frac{1}{\mu_j}e_j 
								 \quad \left[ \frac{1}{\mu_j}e_j \mbox{ is an eigenvector of } S\right]\\
								 & = \sum^m_{j=1} x_jS\left(\frac{1}{\mu_j}e_j\right) \\
								 & = \sum^m_{j=1} x_jS\left(S + \delta I\right)^{-1}e_j \\
								 & = S\left(S + \delta I\right)^{-1}x, \quad x \in V.
\end{alignat*}
Thus, \eqref{eq:lemma6} has been established.

Before proving that the limit $P_S$ exists, we show that both members of (\ref{eq:lemma6})
are self-adjoint.  Since $S + \delta I$ is invertible for $| \delta |$ small enough,
Proposition~\ref{prop:inv} implies that 
\[
\left[\left(S + \delta I\right)^{-1}\right]^* = \left[\left(S + \delta I\right)^*\right]^{-1} = 
\left(S + \delta I\right)^{-1},
\]
where the last equality follows since $S + \delta I$ is self-adjoint.  Thus, its inverse is self-adjoint 
as well and so:
\begin{alignat*}{2}
(x, \left[\left(S + \delta I\right)^{-1}S\right]^*y) & \equiv (\left(S + \delta I\right)^{-1}Sx, y) \\
				   & = (S\left(S + \delta I\right)^{-1}x, y) \quad [ (\ref{eq:lemma6}) ] \\
				   & = (x, \left(S + \delta I\right)^{-1}Sy) \quad [ S, (S + \delta I)^{-1} \mbox{ self-adjoint} ].
\end{alignat*}
This proves that
\begin{equation}
\label{eq:lemma8}
\nonumber
P_S(\delta) \equiv \left(S + \delta I\right)^{-1}S = S\left(S + \delta I\right)^{-1}
\end{equation}
is self-adjoint.

At this point we are ready to prove that
\[
P_Sx \equiv \lim_{\delta \rightarrow 0}P_S(\delta)x = \hat{x}
\]
exists for each $x \in V$ and where $\hat{x}$ is the orthogonal projection of $x$ onto
${\cal R}(S)$.  According to Proposition~\ref{prop:decomp} we split $x$ as
\[
x = \hat{x} + \tilde{x}, \quad \hat{x} \in {\cal R}(S), \quad \tilde{x} \in {\cal N}(S^*) = {\cal N}(S).
\]
Thus,
\[
P_S(\delta)x = P_S(\delta)\hat{x}.
\]
Since $\hat{x}$ is the orthogonal projection onto ${\cal R}(S)$ there is an element $x_0 \in V$
such that
\[
\hat{x} = Sx_0.
\]
According to the spectral theorem~\ref{thm:spectral} we can decompose $x_0$ in terms of the eigenvectors of 
$S$ as follows:
\[
P_S(\delta)x = \left(S + \delta I\right)^{-1}S^2x_0
             = \sum^m_{j = 1} \frac{\lambda_j^2}{\lambda_j + \delta}(e_j, x_0)e_j.
\]
But
\[
\lim_{\delta \rightarrow 0}\frac{\lambda_j^2}{\lambda_j + \delta} = \lambda_j, \quad j = 1, \ldots, m.
\]
Hence, the limit 
\[
\lim_{\delta \rightarrow 0}P_S(\delta)x = \sum^m_{j = 1} \lambda_j(e_j, x_0)e_j
     = \sum^m_{j = 1} (e_j, x_0)Se_j = Sx_0 = \hat{x}
\]
exists for all $x\in V$, which proves (\ref{eq:lemma1}).  The case
\[
\lim_{\delta \rightarrow 0} S\left(S + \delta I\right)^{-1} \quad \left[(\ref{eq:lemma2})\right]
\]
is handled in a similar manner.  \hfill$\Box$

Next, we prove some direct consequences of Lemma~\ref{lemma:proj}.

\begin{cor}
\label{cor:proj}
Let $S:V \rightarrow V$ be self-adjoint.  Then $P_S$ defined by (\ref{eq:lemma1}),
(\ref{eq:lemma2}) is self-adjoint.
\end{cor}

\noindent
{\bf Proof:}  The adjoint of $P_S$ is defined as
\begin{alignat*}{1}
(x, P_S^*y) & \equiv (P_Sx,y) \\ 
            & = \lim_{\delta \rightarrow 0}(P_S(\delta)x, y) \quad [ P_S(\delta) \mbox{ self-adjoint} ]\\
			& = \lim_{\delta \rightarrow 0}(x, P_S(\delta)y) \\
			& = (x,P_Sy),
\end{alignat*}
which proves the corollary. \hfill$\Box$

\begin{prop}
\label{prop:projlimit}
Let $S:V \rightarrow V$ be self-adjoint.  Then
\begin{alignat}{1}
{\cal N}(P_S)   & = {\cal N}(S)   \label{eq:projlimit1} \\
{\cal R}(P_S)   & = {\cal R}(S)   \label{eq:projlimit2} 
\end{alignat}
where $P_S$ is defined by (\ref{eq:lemma1}), (\ref{eq:lemma2}).
\end{prop}

\noindent
{\bf Proof:}
Suppose that
\begin{alignat*}{2}
x \in {\cal N}(S) & \Longleftrightarrow Sx = 0 \\
                  & \Longleftrightarrow P_S(\delta)x = 0 \\
				  & \Longrightarrow P_Sx = \lim_{\delta \rightarrow 0}P_S(\delta)x = 0 \\
				  & \Longleftrightarrow x \in {\cal N}(P_S) .
\end{alignat*}
Thus
\[
{\cal N}(S) \subset {\cal N}(P_S).
\]

Conversely, suppose that $ x \in {\cal N}(P_S)$:
\begin{alignat*}{2}
0 = P_Sx & = \lim_{\delta \rightarrow 0}\left(S + \delta I\right)^{-1}Sx  \\
         & = \lim_{\delta \rightarrow 0} \sum^m_{j = 1} x_j\frac{\lambda_j}{\lambda_j + \delta}e_j 
		   = \sum^m_{j = 1} x_j\delta_j e_j,
\end{alignat*}
where
\[
\delta_j = 
\left\{
\begin{array}{cl}
0 \quad &\mbox{if } \lambda_j = 0 \\
1 \quad &\mbox{if } \lambda_j \neq 0 
\end{array}
\right. .
\]
Since $e_j$ are linearly independent (even orthogonal) it follows that
\[
x_j =0, \quad \lambda_j \neq 0.
\]
Thus,
\[
x = \sum_j x_je_j,
\]
where we agree to sum only over those indices $j$ for which $\lambda_j = 0$.  Hence,
\[
Sx = \sum_j x_jSe_j = \sum_j \lambda_jx_j = 0,
\]
which demonstrates that $x \in {\cal N}(S)$.  Thus, ${\cal N}(P_S) \subset {\cal N}(S)$ and 
we have established (\ref{eq:projlimit1}).

To prove (\ref{eq:projlimit2}) we use (\ref{eq:projlimit1}) together with Proposition~\ref{prop:nr1}:
\begin{alignat*}{2}
{\cal R}(P_S) & = {\cal N}(P^*_S)^\perp = [\mbox{Cor.~\ref{cor:proj}}] = {\cal N}(P_S)^\perp \\
              & = {\cal N}(S)^\perp = {\cal R}(S^*) = {\cal R}(S),
\end{alignat*}
which concludes the proof. \hfill$\Box$

\begin{prop}
\label{prop:proj1}
Let $S:V \rightarrow V$ be self-adjoint and let $P_S$ be given by (\ref{eq:lemma1}), 
(\ref{eq:lemma2}).  Then $P_S$ is a projection:
\[
P_S^2 = P_S.
\]
\end{prop}

\noindent
{\bf Proof:}
Let $x \in V$.  Split $x$:
\[
x = \hat{x} + \tilde{x}, \quad [\mbox{Prop.~\ref{prop:decomp}}]
\]
where $\hat{x}$ is the orthogonal projection of $x$ onto ${\cal R}(S)$ and 
$\tilde{x} \in {\cal N}(S)$.  From Lemma~\ref{lemma:proj} we know that
\[
P_Sx = \hat{x}, 
\]
where $\hat{x}$ is defined as above.  Hence,
\[
P_S^2x = P_S\hat{x} = [\mbox{(\ref{eq:projlimit1})}] = P_S(\hat{x} + \tilde{x}) = P_Sx,
\]
which shows that $P_S^2 = P_S$. \hfill$\Box$

At this point we have all the preliminary results required to construct the pseudoinverse
of $T:V \rightarrow W$.  

\begin{thm}
\label{thm:pseudo}
Let $T:V \rightarrow W$.  The pseudoinverse 
$T^+:W\rightarrow V$ defined by
\begin{alignat}{2}
T^+ & \equiv \lim_{\delta \rightarrow 0}\left(T^*T + \delta^2 I\right)^{-1}T^*   \label{eq:pseudo1} \\
    & =      \lim_{\delta \rightarrow 0} T^* \left(TT^* + \delta^2 I\right)^{-1} \label{eq:pseudo2}
\end{alignat}
always exists.  Furthermore, for any $z\in W$
\begin{equation}
\label{eq:pseudo3}
\hat{x} = T^+z
\end{equation}
minimizes
\begin{equation}
\label{eq:pseudo4}
\| z - Tx \|^2
\end{equation}
and has minimum norm.
\end{thm}

\noindent
{\bf Proof:} First we must ensure that the definitions of the pseudoinverse (\ref{eq:pseudo1}) and 
(\ref{eq:pseudo2}) make sense.  Thanks to the previous propositions and lemmas this will be very 
straightforward.  We then follow the usual pattern of splitting $z$ onto 
${\cal R}(T) \oplus {\cal N}(T^*)$ and then apply Lemma~\ref{lemma:proj} to prove that the limits in 
(\ref{eq:pseudo1}) and (\ref{eq:pseudo2}) exist.  As a by-product we obtain the splitting of $x$ onto 
${\cal R}(T^*) \oplus {\cal N}(T)$.  The results (\ref{eq:pseudo3}) and (\ref{eq:pseudo4}) will then 
be direct consequences of Theorem~\ref{thm:ls1}.

\begin{enumerate}
\item[(i)] $T^*T + \delta^2 I$ and $TT^* + \delta^2 I$ are self-adjoint and invertible.  According to 
Proposition~\ref{prop:sa2} both operators are self-adjoint and non-singular (one-to-one).  Since they
are self-adjoint, it follows from Proposition~\ref{prop:onto} that they are onto as well.  Hence, 
both operators are invertible.
\item[(ii)] $\left(T^*T + \delta^2 I\right)^{-1}$ and $\left(TT^* + \delta^2 I\right)^{-1}$ are self-adjoint.  
This follows immediately from the previous step and Proposition~\ref{prop:inv}.
\item[(iii)] $\left(T^*T + \delta^2 I\right)^{-1}T^* = T^*\left(TT^* + \delta^2 I\right)^{-1}$.  This can be 
shown as follows:
\begin{alignat*}{2}
T^* = I\cdot T^* & = \left(T^*T + \delta^2 I\right)^{-1}\left(T^*T + \delta^2 I\right)T^* \\
                 & = \left[\left(T^*T + \delta^2 I\right)^{-1}T^*\right]\left(TT^* + \delta^2 I\right).
\end{alignat*}
Hence,
\[
T^*\left(TT ^*+ \delta^2 I\right)^{-1} = \left(T^*T + \delta^2 I\right)^{-1}T^*.
\]
\item[(iv)] Let $z \in W$.  Split $z$ as
\[
z = \hat{z} + \tilde{z},
\]
where $\hat{z}$ is the projection onto ${\cal R}(T)$ and $\tilde{z} \in {\cal N}(T^*)$ 
(Proposition~\ref{prop:decomp}).  Hence, there is a vector $x \in V$ such that
\begin{equation}
\label{eq:pseudo5}
\hat{z} = Tx.
\end{equation} 
Thus,
\[
\left(T^*T + \delta^2 I\right)^{-1}T^*z = \left(T^*T + \delta^2 I\right)^{-1}T^*Tx.
\]
Let $S \equiv T^*T:V \rightarrow V$.  All requirements of Lemma~\ref{lemma:proj} are met and 
the following limit exists:
\begin{alignat*}{2}
T^+z & \equiv \lim_{\delta \rightarrow 0}\left(T^*T + \delta^2 I\right)^{-1}T^*z    \\
     & =      \lim_{\delta \rightarrow 0}\left(T^*T + \delta^2 I\right)^{-1}T^*Tx = \hat{x},
	 \quad [\mbox{ Lemma~\ref{lemma:proj} }]
\end{alignat*}
where $\hat{x}$ is the projection of $x$ onto ${\cal R}(T^*T), [\mbox{(\ref{eq:lemma3})}]$.
But ${\cal R}(T^*T) = {\cal R}(T^*), [\mbox{(\ref{eq:nr12})}]$.  We have thus arrived at the usual
decomposition of $x$:
\[
x = \hat{x} + \tilde{x},  \quad \hat{x} \in {\cal R}(T^*),  \quad \tilde{x} \in {\cal N}(T).
\]
Substituting this decomposition in (\ref{eq:pseudo5}):
\[
\hat{z} = Tx = T(\hat{x} + \tilde{x}) = T\hat{x}.
\]
Hence, $T^+z = \hat{x} \in {\cal R}(T^*)$ satisfies 
\[
Tx = \hat{z}.
\]
From Theorem~\ref{thm:ls1} it follows that (\ref{eq:pseudo3}) has minimum norm and minimizes 
(\ref{eq:pseudo4}).
\end{enumerate}

\noindent
This completes the proof. \hfill$\Box$

\begin{prop}
\label{prop:nr3}
Let $T:V \rightarrow W$ and $x \in V, z \in W$.  Then
\begin{align}
      T^+Tx & = \hat{x},   \quad \hat{x}   \mbox{ projection of $x$ on } {\cal R}(T^*) \label{eq:canproj1} \\
      TT^+z & = \hat{z},   \quad \hat{z}   \mbox{ projection of $z$ on } {\cal R}(T)   \label{eq:canproj2} \\
{\cal R}(T^+) & = {\cal R}(T^*)   \label{eq:nr31} \\
{\cal N}(T^+) & = {\cal N}(T^*).  \label{eq:nr32}
\end{align}
\end{prop}

\noindent
{\bf Proof:} By (\ref{eq:pseudo1}):
\[
T^+Tx = \lim_{\delta \rightarrow 0}\left(T^*T + \delta^2 I\right)^{-1}T^*Tx = \hat{x} ,
\]
where $\hat{x}$ is the orthogonal projection of $x\in V$ onto ${\cal R}(T^*T) = {\cal R}(T^*)$ 
according to Lemma~\ref{lemma:proj} (use (\ref{eq:lemma1}) and $S=T^*T$), which proves (\ref{eq:canproj1}).

Similarly, (Lemma~\ref{lemma:proj}, $V\rightarrow W$, (\ref{eq:lemma2}) and $S=TT^*)$:
\[
TT^+z = \lim_{\delta \rightarrow 0}TT^*\left(TT^* + \delta^2 I\right)^{-1}z = \hat{z},
\]
where $\hat{z}$ is the orthogonal projection of $z\in W$ onto ${\cal R}(TT^*) = {\cal R}(T)$.
Eq.~(\ref{eq:canproj2})~is thus proved.

Next, we want to to prove (\ref{eq:nr31}).  To this end, let $x \in {\cal R}(T^*)$.  Then
\[
\hat{x} = x,
\]
since $\hat{x}$ is the projection of $x$ onto ${\cal R}(T^*)$ by definition.  Hence,
by (\ref{eq:canproj1}):
\[
T^+Tx = x,
\]
which shows that $x \in {\cal R}(T^+)$, that is,
\[
{\cal R}(T^*) \subset  {\cal R}(T^+).
\]
To prove the reverse inclusion, assume that $x \in {\cal R}(T^+)$.  Then there is a vector
$z \in W$ such that
\begin{alignat*}{2}
x = T^+z & \equiv \lim_{\delta \rightarrow 0} \left(T^*T + \delta^2 I \right)^{-1}T^*z \\
         & = \lim_{\delta \rightarrow 0} \left(T^*T + \delta^2 I \right)^{-1}T^*\hat{z} \\
         & = \lim_{\delta \rightarrow 0} \left(T^*T + \delta^2 I \right)^{-1}T^*Ty = \hat{y},
\end{alignat*}
where $\hat{y} \in {\cal R}(T^*T) = {\cal R}(T^*)$ by Lemma~\ref{lemma:proj}.  But $x = \hat{y}$
according to the previous expression.  Thus, $x \in {\cal R}(T^*)$, that is
\[
{\cal R}(T^+) \subset  {\cal R}(T^*),
\]
which demonstrates (\ref{eq:nr31}).

It remains to prove (\ref{eq:nr32}).  Decompose $z$:
\begin{equation}
\label{ex:nr33}
z = \hat{z} + \tilde{z},
\end{equation}
where $\hat{z}$ and $\tilde{z}$ are the projections of $z$ onto ${\cal R}(T)$ and ${\cal N}(T^*)$.
If $z$ in ${\cal N}(T^+)$, then by (\ref{eq:canproj2}):
\begin{equation}
\label{ex:nr34}
\hat{z} = 0.
\end{equation}
Combining (\ref{ex:nr33}) and (\ref{ex:nr34}):
\[
z = \tilde{z} \in {\cal N}(T^*).
\]
Hence,
\[
{\cal N}(T^+) \subset {\cal N}(T^*).
\]
Conversely, suppose that $z \in {\cal N}(T^*)$.  Then
\[
T^+z \equiv \lim_{\delta \rightarrow 0}\left(T^*T + \delta^2 I \right)^{-1}T^*z = 0.
\]
Thus,
\[
{\cal N}(T^*) \subset {\cal N}(T^+),
\]
which concludes the proof. \hfill$\Box$

\begin{cor}
\label{cor:canonicalproj}
Let $T:V \rightarrow W$.  The following statements hold:
\begin{align}
\left(T^+T\right)^* & = T^+T \label{eq:canproj5} \\
\left(TT^+\right)^* & = TT^+ \label{eq:canproj6} \\
\left(T^+T\right)^2 & = T^+T \label{eq:canproj7} \\
\left(TT^+\right)^2 & = TT^+ \label{eq:canproj8}.
\end{align}
\end{cor}

\noindent{\bf Proof:}  
To prove (\ref{eq:canproj5}) - (\ref{eq:canproj8}) we note that (\ref{eq:pseudo1}) and (\ref{eq:pseudo2})
imply
\begin{alignat*}{2}
T^+T & = \lim_{\delta \rightarrow 0} \left(T^*T + \delta^2 I\right)^{-1}T^*T \\
TT^+ & = \lim_{\delta \rightarrow 0} TT^* \left(TT^* + \delta^2 I\right)^{-1}.
\end{alignat*}
Corollary~\ref{cor:proj} implies (\ref{eq:canproj5}), (\ref{eq:canproj6}).  Eqs.~(\ref{eq:canproj7}),
(\ref{eq:canproj8}) follow from Proposition~\ref{prop:proj1}. \hfill$\Box$

\begin{rem}
\label{rem:canonicalproj}
By (\ref{eq:canproj5}) - (\ref{eq:canproj8}), it would seem natural
for the pseudoinverse $T^+$ to satisfy the following conditions:
\begin{alignat*}{1}
\left(T^+T\right)^* & = T^+T \\
\left(TT^+\right)^* & = TT^+ \\
            T^+TT^+ & = T^+ \\
            TT^+T   & = T.
\end{alignat*}
This is indeed the case and will be shown in Theorem~\ref{thm:penrose}.  It turns out that the above conditions
are necessary and sufficient conditions for $T^+$ to be the pseudoinverse of $T$, which was first proved in
\cite{rp:gifm}. \hfill$\Box$
\end{rem}

\begin{exa}
\label{exa:pseudo}
Generalized division is the gist of pseudoinversion.  This is clearly illustrated by applying the spectral 
theorem to the pseudoinverse of a self-adjoint operator $T:V\rightarrow V$.  Let $e_j, j = 1, \ldots , m$ 
denote the ortho-normal eigenvectors of $T$.  Decompose $x \in V$
\[
x = \sum^m_{j=1} (e_j,x)e_j.
\]
Hence,
\begin{alignat*}{2}
T^+x & \equiv \lim_{\delta \rightarrow 0}\left(T^*T + \delta^2 I\right)^{-1}T^*x \\
     & = \lim_{\delta \rightarrow 0}\left(T^2 + \delta^2 I\right)^{-1}Tx    \\
	 & = \sum^m_{j=1} \lim_{\delta \rightarrow 0}\left(T^2 + \delta^2 I\right)^{-1}T(e_j,x)e_j \\
	 & = \sum^m_{j=1} \lim_{\delta \rightarrow 0}\frac{\lambda_j}{\lambda_j^2 + \delta^2}(e_j,x)e_j \\
	 & = \sum^m_{j=1} \lambda_j^+(e_j,x)e_j,
\end{alignat*}
where $\lambda_j^+ \in \mathbb{R}$ is defined as
\begin{equation}
\label{eq:pseudoscal}
\nonumber
\lambda_j^+ = 
\left\{
\begin{array}{cl}
0,           & \quad \lambda_j = 0, \\
1/\lambda_j, & \quad \lambda_j \neq 0.
\end{array}
\right.
\end{equation}
Thus, $\lambda_j^+$ extends scalar division to also include $0$.  The impact of this definition is clear from
the spectral decomposition of $T^+$:
\[
T^+x = \sum^m_{j=1} \lambda_j^+(e_j,x)e_j.
\]
Consequently, $T^+$ projects a vector in $V$ onto the linear manifold spanned by the eigenvectors that 
correspond to non-zero eigenvalues of the self-adjoint operator $T$. \hfill$\Box$
\end{exa}

\begin{exa}
\label{exa:pseudoproj}
What happens if $T$ itself is a projection, i.e., $T = T^2$?  Reusing the notation from the previous example:
\begin{alignat*}{2}
T^+x & \equiv \lim_{\delta \rightarrow 0}\left(T^2 + \delta^2 I\right)^{-1}Tx    \\
     & = \lim_{\delta \rightarrow 0}\left(T^2 + \delta^2 I\right)^{-1}T^3x \quad [T = T^3]   \\
	 & = \sum^m_{j=1} \lim_{\delta \rightarrow 0}\frac{\lambda_j^3}{\lambda_j^2 + \delta^2}(e_j,x)e_j \\
	 & = \sum^m_{j=1} \lambda_j(e_j,x)e_j \\
	 & = \sum^m_{j=1} T(e_j,x)e_j = Tx,
\end{alignat*}
where we used
\[
\lim_{\delta \rightarrow 0}\frac{\lambda_j^3}{\lambda_j^2 + \delta^2} = \lambda_j
\]
for any real $\lambda_j$ (in this particular case $\lambda_j=0,1$).  Thus, 
\[
T^+ = T
\]
if $T$ is a self-adjoint projection. \hfill$\Box$
\end{exa}

The next proposition shows that it is enough to be able to compute pseudoinverses of self-adjoint operators.
\begin{prop}
\label{prop:altpseudo}
Let $T:V \rightarrow W$.  Then
\begin{alignat}{1}
T^+ & = \left(T^*T\right)^+T^*             \label{eq:altpseudo1} \\
\left(T^*\right)^+ & = \left(T^+\right)^*  \label{eq:altpseudo2} \\
T^+ & = T^*\left(TT^*\right)^+             \label{eq:altpseudo3}.
\end{alignat}
\end{prop}

\noindent{\bf Proof:}
Let $z \in W$.  Split 
\[
z = \hat{z} + \tilde{z},
\]
where $\hat{z}$ and $\tilde{z}$ are the familiar projections onto ${\cal R}(T)$ and ${\cal N}(T^*)$.  
Thus, there is an $x \in V$ such that
\[
\hat{z} = Tx.
\]
The definition of the pseudoinverse (\ref{eq:pseudo1}), (\ref{eq:pseudo2}) implies
\begin{alignat*}{2}
\left(T^*T\right)^+T^*z & \equiv \lim_{\delta \rightarrow 0}\left[(T^*T)^*(T^*T) + \delta^2 I\right]^{-1}(T^*T)^*T^*z  \\
     & = \lim_{\delta \rightarrow 0}\left[(T^*T)^2 + \delta^2 I\right]^{-1}(T^*T)T^*\hat{z} \quad [T^*\tilde{z} = 0] \\
	 & = \lim_{\delta \rightarrow 0}\left[(T^*T)^2 + \delta^2 I\right]^{-1}(T^*T)^2x  \\
	 & = \sum^m_{j=1} \lim_{\delta \rightarrow 0}\frac{\lambda_j^2}{\lambda_j^2 + \delta^2}(e_j,x)e_j,
\end{alignat*}
where $\lambda_j \geq 0$ and $e_j \in \mathbb{R}^m$ are the eigenvalues and eigenvectors of the self-adjoint 
operator $T^*T$.  Since $\lambda_j \geq 0$ it follows that
\[
\lim_{\delta \rightarrow 0}\frac{\lambda_j^2}{\lambda_j^2 + \delta^2} = 
\lim_{\delta \rightarrow 0}\frac{\lambda_j}{\lambda_j + \delta^2} =
\left\{
\begin{array}{cl}
0, & \quad \lambda_j = 0, \\
1, & \quad \lambda_j \neq 0.
\end{array}
\right.
\]
Hence,
\begin{alignat*}{2}
\left(T^*T\right)^+T^*z & = \sum^m_{j=1} \lim_{\delta \rightarrow 0}\frac{\lambda_j}{\lambda_j + \delta^2}(e_j,x)e_j  \\
	 & = \sum^m_{j=1} \lim_{\delta \rightarrow 0}\left(T^*T + \delta^2 I\right)^{-1}T^*T(e_j,x)e_j  \\
	 & = \lim_{\delta \rightarrow 0}\left(T^*T + \delta^2 I\right)^{-1}T^*Tx \\
	 & = \lim_{\delta \rightarrow 0}\left(T^*T + \delta^2 I\right)^{-1}T^*\hat{z} \\
	 & = \lim_{\delta \rightarrow 0}\left(T^*T + \delta^2 I\right)^{-1}T^*z  \quad [T^*\tilde{z} = 0] \\
	 & \equiv T^+z,
\end{alignat*}
which proves (\ref{eq:altpseudo1}).

To prove (\ref{eq:altpseudo2}) we note that 
\begin{alignat*}{2}
\left(T^*\right)^+ & \equiv \lim_{\delta \rightarrow 0}\left(T^{**}T^* + \delta^2 I\right)^{-1}T^{**}  \\
                   & = \lim_{\delta \rightarrow 0}\left(TT^* + \delta^2 I\right)^{-1}T^{**}  
				       \quad [\mbox{Prop.\ref{prop:adjoint}}] \\
				   & = \lim_{\delta \rightarrow 0}\left[\left(TT^* + \delta^2 I\right)^*\right]^{-1}T^{**} \\
                   & = \lim_{\delta \rightarrow 0}\left[\left(TT^* + \delta^2 I\right)^{-1}\right]^*T^{**}				   
				       \quad [\mbox{Prop.~\ref{prop:inv}}] 	\\  
				   & = \left[\lim_{\delta \rightarrow 0}T^*\left(TT^* + \delta^2 I\right)^{-1}\right]^*	\\ 
				   & \equiv \left(T^+\right)^*,
\end{alignat*}
which shows that (\ref{eq:altpseudo2}) holds.

Finally, (\ref{eq:altpseudo1}) and (\ref{eq:altpseudo2}) can be combined to prove (\ref{eq:altpseudo3}).  
Taking the adjoint of both sides in (\ref{eq:altpseudo2}) yields
\begin{equation}
\label{eq:altpseudo4}
T^+ = \left[\left(T^*\right)^+\right]^*.
\end{equation}
Next, apply (\ref{eq:altpseudo1}) to $T^*$:
\begin{equation}
\label{eq:altpseudo5}
\nonumber
\left(T^*\right)^+ = (T^{**}T^*)^+T^{**} = (TT^*)^+T = 
\left[\left(TT^*\right)^*\right]^+T.
\end{equation}
Taking the adjoint a second time and applying (\ref{eq:altpseudo4}) proves (\ref{eq:altpseudo3}). \hfill$\Box$

We are now in a position to prove Penrose's characterization of the pseudoinverse of $T$ \cite{rp:gifm}
for finite-dimensional inner product spaces.

\begin{thm}
\label{thm:penrose}
Let $T:V \rightarrow W$.  Then $S:W \rightarrow V$ is the pseudoinverse of $T$ iff $S$ satisfies
\begin{alignat}{1}
\left(ST\right)^* & = ST \label{eq:penrose1} \\
\left(TS\right)^* & = TS \label{eq:penrose2} \\
              TST & = T  \label{eq:penrose3} \\
              STS & = S. \label{eq:penrose4}
\end{alignat}
\end{thm}

\noindent{\bf Proof:}  {\em Necessity}: This follows more or less directly from the canonical projections of
Corollary~\ref{cor:canonicalproj}.  Substituting $S = T^+$ in (\ref{eq:canproj5}) and (\ref{eq:canproj6}) 
implies (\ref{eq:penrose1}) and (\ref{eq:penrose2}).  

Next, multiply (\ref{eq:canproj1}) and 
(\ref{eq:canproj2}) with $T$ and $T^+$:
\begin{alignat*}{2}
TT^+Tx   & = T\hat{x} = Tx \quad\quad [\tilde{x} \in {\cal N}(T)] \\
T^+TT^+z & = T^+\hat{z} = T^+z. \quad[\tilde{z} \in {\cal N}(T^*) = {\cal N}(T^+), \mbox{Prop.~\ref{prop:nr3}}]
\end{alignat*}
Replacing $T^+$ with $S$ establishes (\ref{eq:penrose3}) and (\ref{eq:penrose4}), which proves that
(\ref{eq:penrose1}) - (\ref{eq:penrose4}) are necessary conditions.

{\em Sufficiency}: We will show that (\ref{eq:penrose1}) - (\ref{eq:penrose4}) imply that
\[
S = T^+TT^+ = T^+,
\]
where the last equality was established in the first part of the proof, which also showed that
\begin{equation}
\label{eq:penrose5}
TT^+(Tx) = Tx, \quad x \in V \quad \Longleftrightarrow \quad TT^+T = T.
\end{equation}
Hence,
\begin{equation}
\label{eq:penrose6}
T = [(\ref{eq:penrose3})] = T(ST) = [(\ref{eq:penrose1})] = T(ST)^* = TT^*S^*.
\end{equation}
We can now express the canonical projection $T^+T$ as
\begin{alignat}{2}
T^+T & = (T^+T)T^*S^*   & \quad [(\ref{eq:penrose6})] \nonumber \\
     & = (T^+T)^*T^*S^* & \quad [(\ref{eq:canproj5})] \nonumber \\	 
     & = T^*S^* = ST.   & \quad [(\ref{eq:penrose5})] \label{eq:penrose7}
\end{alignat}
But
\begin{equation}
\label{eq:penrose8}
S^* = [(\ref{eq:penrose4})] = \left[S(TS)\right]^* = (TS)^*S^* = [(\ref{eq:penrose2})] = TSS^*.
\end{equation}
Premultiplying (\ref{eq:penrose8}) with the canonical projection $TT^+$ yields
\[
(TT^+)S^* = (TT^+)(TSS^*) = [(\ref{eq:penrose5})] = TSS^* = [(\ref{eq:penrose8})] = S^*.
\]
Thus
\[
S = S\left( TT^+ \right)^* = [(\ref{eq:canproj6})] = STT^+ = [(\ref{eq:penrose7})] = T^+TT^+ = T^+.
\]
This completes the proof. \hfill$\Box$

\begin{rem}
\label{rem:penrose1}
It should be noted that (\ref{eq:penrose3}) and (\ref{eq:penrose4}) are invariants; they do 
not depend the scalar products.  In fact, the matrix representation is identical to the 
operator formulation.  Eq.~(\ref{eq:penrose3}) can be expressed as a set of "telescoping" 
constraints:
\begin{equation}
\label{eq:penrose9}
w = Ty,  \quad y = Sz, \quad z = Tx,  \quad \mbox{and} \quad w = Tx.
\end{equation}
Introducing ortho-normal base vectors $e_i, i=0,\ldots, m$ and $f_j, j=1,\ldots, n$ in $V$ and 
$W$, we can express $x, y \in \mathbb{R}^m$ and $z, w \in \mathbb{R}^n$ as
\[
x = \sum^m_{j=1}x_je_j, \quad y = \sum^m_{j=1}y_je_j, \quad z = \sum^n_{j=1}z_jf_j,
\quad w = \sum^n_{j=1}w_jf_j.
\]
Substituting these expressions into (\ref{eq:penrose9}) and using the linear independence of 
the base vectors we recover {\em exactly} the same equations (\ref{eq:penrose9}), but this time 
$T$ and $S$ should be interpreted as matrices in $\mathbb{R}^{n\times m}$ and 
$\mathbb{R}^{m\times n}$ with elements
\[
T_{ij} = \langle f_i, Te_j \rangle, \quad S_{ij} = (e_i, Sf_j).
\]
Thus,
\[
TST = T
\]
is true no matter if $T$ and $S$ are interpreted as operators or matrices.  The same conclusion holds for 
(\ref{eq:penrose4}). \hfill$\Box$
\end{rem}

\begin{rem}
\label{rem:penrose2}
While (\ref{eq:penrose3}) and (\ref{eq:penrose4}) can be interpreted as operator or matrix conditions, 
the same conclusion does not apply to the orthogonality constraints (\ref{eq:penrose1}) and 
(\ref{eq:penrose2}).  In operator form  we have $ST:V \rightarrow V$.  Hence, by definition
\[
(x, (ST)^*y) \equiv (STx, y), \quad x, y \in V.
\]
But $(ST)^* = ST, [(\ref{eq:penrose1})]$ and so
\[
(x, STy) \equiv (STx, y).
\]
In matrix form we have
\[
x^TH_1STy = (STx)^TH_1y = x^T(ST)^TH_1y,  \quad x, y \in \mathbb{R}^{m}.
\]
Hence, (\ref{eq:penrose1}) becomes
\begin{equation}
\label{eq:penrose10}
(ST)^TH_1 =  H_1(ST), \quad H_1 > 0 \in \mathbb{R}^{m \times m}.
\end{equation}

Similarly, for $TS:W \rightarrow W$ condition (\ref{eq:penrose1}) translates to
\begin{equation}
\label{eq:penrose11}
(TS)^TH_2 =  H_2(TS), \quad H_2 > 0 \in \mathbb{R}^{n \times n}.
\end{equation}
In summary:  The pseudoinverse definition (\ref{eq:pseudo1}), (\ref{eq:pseudo2}) and Penrose's 
characterization (\ref{eq:penrose1}) - (\ref{eq:penrose4}) carry over verbatim to finite-dimensional
inner product spaces as long as orthogonality is with respect to the inner product spaces.  
In case of matrix representation, then (\ref{eq:penrose1}), (\ref{eq:penrose2}) must be replaced 
by (\ref{eq:penrose10}), (\ref{eq:penrose11}) for the orthogonality conditions to be valid for 
general scalar products. \hfill$\Box$
\end{rem}

We close out this section with a result that states necessary and sufficient conditions for
\[
T^* = T^T
\]
to be true.  This property will simplify the actual computation of the pseudoinverse, since it 
will allow $T^+$ to be constructed without involving the norms $H_1$ or $H_2$.

\begin{prop}
\label{prop:adjtrans}
Let $T:V \rightarrow W$. Then
\begin{equation}
\label{eq:adjtransp1}
T^* = T^T
\end{equation}
iff
\begin{equation}
\label{eq:adjtransp2}
TH_1 = H_2T.
\end{equation}
\end{prop}

\noindent
{\bf Proof:} 

{\em Necessity}:  From (\ref{eq:adjtransp})
\[
H_1T^* = T^TH_2.
\]
Hence, $T^* = T^T$ implies
\[
H_1T^T = T^TH_2,
\]
and since $H_1, H_2$ are symmetric it follows that
\[
TH_1 = H_2T,
\]
which proves necessity.

{\em Sufficiency}:  Assume that (\ref{eq:adjtransp2}) holds.  Then, by (\ref{eq:adjtransp}):
\begin{alignat*}{2}
T^* & = H_1^{-1}T^TH_2 \\
    & = H_1^{-1}(H_2T)^T &\qquad [(\ref{eq:adjtransp2})] \\
    & = H_1^{-1}(TH_1)^T \\	 
	& = T^T,
\end{alignat*}
which finishes the proof.  \hfill$\Box$

\section{Difference operators and summation by parts}
\label{sec:sbp}
Boundary projections are closely related to summation-by-parts (SBP) operators.  Before constructing the 
projection operators, the basics of the SBP operators will be presented to facilitate the discussion in 
the subsequent sections.  For more details on the theory of SBP operators, the reader is referred to 
\cite{gko:tdpdm}.

First, we introduce notation that will prove useful when dealing with matrix representations of 
boundary modified difference operators.  For any matrix $A \in \mathbb{R}^{n\times m}$ let
\begin{equation}
\label{eq:aadj}
A^\tau \equiv J_nAJ_m, 
\end{equation}
where the anti-diagonal matrices $J_m, J_n$ are defined as
\begin{equation}
\label{eq:adiag}
\quad J_p \equiv 
\begin{pmatrix}
& & 1 \\
& \iddots \\
1
\end{pmatrix}
\in \mathbb{R}^{p\times p}.
\end{equation}
Thus, $A^\tau$ is obtained by reversing the rows and columns of $A$:
\[ 
A_{ij} \leftrightarrow A_{n-i+1,m-j+1}.
\]
Note that (\ref{eq:aadj}) preserves the shape of $A$ as opposed to the 
usual matrix transposition.  It follows immediately that $J = J^{-1}$ 
and if $A$ is invertible:
\[
A^{-\tau} \equiv \left(A^\tau\right)^{-1} = \left(A^{-1}\right)^\tau.
\]

We will be concerned with PDEs defined on the unit interval $[0, 1]$.  
On this interval we define a uniform grid
\begin{equation}
\label{eq:grid}
\nonumber
x_j \equiv jh, \quad j = 0, \ldots, N, \quad h \equiv \frac{1}{N},
\end{equation}
and the corresponding grid functions $u, v$:
\begin{equation}
\label{eq:uv}
u \equiv 
\begin{pmatrix}
u_0 \\
\vdots \\
u_N
\end{pmatrix}, \quad
v \equiv 
\begin{pmatrix}
v_0 \\
\vdots \\
v_N
\end{pmatrix} \in \mathbb{R}^{N+1}.
\end{equation}
Let $D \in \mathbb{R}^{(N+1)\times(N+1)}$ be the matrix representing a  consistent approximation 
of $\partial / \partial x$.  At interior grid points $x_i$, $r \leq i \leq N-r$, the difference 
operator $D$ corresponds to the standard anti-symmetric, $2p$-order accurate difference stencil
\begin{equation}
\label{eq:stencil}
(Dv)_i = \frac{1}{h}\sum^p_{j=-p}d_jv_{i+j}, \quad d_{-j} = -d_j.
\end{equation}
In the boundary regions the structure of $D$ is given by
\begin{equation}
\label{eq:bstencil}
(Dv)_i = \frac{1}{h}\sum^{s-1}_{j=0}d_{ij}v_j, \quad 0\leq i < r,
\end{equation}
with a similar expression for the upper boundary region, where we use the "anti-reflected"
difference stencils.  In block form:
\begin{equation}
\label{eq:antireflected}
D = 
\overset{\raise3pt\hbox{$\scriptstyle N+1$}}{
\begin{pmatrix}
D_L \\
D_I \\
D_R
\end{pmatrix}}
\begin{matrix}
{\scriptstyle r} \\
{\scriptstyle N+1-2r} \\
{\scriptstyle r}
\end{matrix},
\end{equation}
where $D_R \equiv -J_rD_LJ_{N+1} = -D^\tau_L$.  

A simple example will serve as an illustration.  Apply a 2nd order stencil at the lower boundary:
\[
(Dv)_0 = \frac{1}{h}\left(-\frac{3}{2}v_0 + 2v_1 -\frac{1}{2}v_2\right) = \frac{1}{h}
\begin{pmatrix}
-\frac{3}{2} &2 &-\frac{1}{2} &0 &\ldots &0
\end{pmatrix}v.
\]
At the upper boundary:
\[
(Dv)_N = \frac{1}{h}\left(\frac{3}{2}v_N - 2v_{N-1} + \frac{1}{2}v_{N-2}\right) = \frac{1}{h}
\begin{pmatrix}
0 &\ldots &0 &\frac{1}{2} &-2 &\frac{3}{2} 
\end{pmatrix}v.
\]
Clearly, the stencil at the upper boundary is obtained by reflecting the difference stencil at the 
lower boundary and by multiplying the reflected stencil by minus one.  A similar relation holds
between $(Dv)_1$ and $(Dv)_{N-1}$ and so on.

Given (\ref{eq:stencil}), the goal is to construct $D_L$ such that $D$ satisfies a summation-by-parts 
rule
\begin{equation}
\label{eq:partsum}
(u, Dv)_H \equiv u_Nv_N - u_0v_0 - (Du, v)_H
\end{equation}
for some inner product
\begin{equation}
\label{eq:innerprod}
(u, v)_H \equiv u^THv, 
\end{equation}
where $H \in \mathbb{R}^{(N+1) \times (N+1)}$ is SPD.  Define polynomial grid functions
\begin{equation}
\label{eq:aux}
\nonumber
\quad x^k \equiv 
\begin{pmatrix}
0^k \\
h^k \\
\vdots \\
(Nh)^k
\end{pmatrix}, \quad
x^0 \equiv {\bf 1} \equiv \begin{pmatrix}
1 \\
1 \\
\vdots \\
1
\end{pmatrix}, \quad
{\bf 0}\equiv \begin{pmatrix}
0 \\
0 \\
\vdots \\
0
\end{pmatrix} \in \mathbb{R}^{N+1}.
\end{equation}
The accuracy requirements of $D$ can then be expressed as 
\begin{equation}
\label{eq:acc}
Dx^0 = 0, \quad Dx^k = kx^{k-1}, \quad 0 < k \leq q,
\end{equation}
for some $q \leq 2p$.  Existence of such operators $D$ and scalar products $(\cdot,\cdot)_H$ was 
first established by Kreiss and Scherer in \cite{ks:fefdhpde,ks:oeeeda}.  Adapting their results to 
the case of two boundaries leads to
\begin{align}
H &= h
\overset{\raise3pt\hbox{$\scriptstyle \hspace*{2pt} r \hspace*{0.5em} N+1-2r \hspace*{0.5em} r$}}{
\begin{pmatrix}
 H_{L}  \\
 & I         \\
 && \hspace*{7pt}H_{L}^\tau
\end{pmatrix}}, \quad 
H_L=\left(h_{ij}\right)_{0\leq i,j < r} > 0, \quad h_{ij} = h_{ji} \label{eq:norm} \\
D &=\frac{1}{h}
\begin{pmatrix}
H_L^{-1} \\
&I \\
&&H_L^{-\tau}
\end{pmatrix}
\overset{\raise3pt\hbox{$\scriptstyle \hspace*{2pt} r \hspace*{1.5em} N+1-2r \hspace*{1.5em} r$}}{
\begin{pmatrix}
 \hspace*{7pt} Q_{L}  & \hspace*{6pt}\tilde{Q}_{I} &  \hspace*{2pt} 0 \\
-\tilde{Q}_{I}^T      & \hspace*{6pt} Q_{I}        & \hspace*{6pt} \bar{Q}_{I} \\
 \hspace*{3pt} 0      & -\bar{Q}_{I}^T             & -Q_{L}^\tau
\end{pmatrix}}
\hspace*{3pt}
\begin{matrix}
{\scriptstyle r} \\
{\scriptstyle N+1-2r} \\
{\scriptstyle r}
\end{matrix}, \label{eq:diffop}
\end{align}
where $h=1/N$ is the mesh size.  The blocks $Q_I$, $\tilde{Q}_I$ and $\bar{Q}_I$ are determined by 
the interior stencil (\ref{eq:stencil}) and thus known.  The explicit structure is:
\begin{align*}
Q_I &= 
\begin{pmatrix}
0         &d_1    &\ldots    &d_p \\
\bar{d}_1 &0      &d_1       &\ldots &d_p \\
\vdots    &\ddots &\ddots    &\ddots &          &\ddots \\
\bar{d}_p &\ldots &\bar{d}_1 &0      &d_1       &\ldots    &d_p \\
          &\ddots &          &\ddots &\ddots    &\ddots    &          &\ddots \\
          &       &\bar{d}_p &\ldots &\bar{d}_1 &0         &d_1       &\ldots    &d_p    \\
		  &       &          &\ddots &          &\ddots    &\ddots    &\ddots    &\vdots \\
		  &		  &          &       &\bar{d}_p &\ldots    &\bar{d}_1 &0         &d_1    \\
		  &       &          &       &          &\bar{d}_p &\ldots    &\bar{d}_1 &0      
\end{pmatrix} \\
\tilde{Q}_I &=
\overset{\raise3pt\hbox{$\scriptstyle \hspace*{2.0em} p \hspace*{1.0em} N+1-2r-p$}}{
\begin{pmatrix}
0                 &&0 \\
\hspace*{3pt}Q_p  &&0
\end{pmatrix}}
\hspace*{3pt}
\begin{matrix}
{\scriptstyle r-p} \\
{\scriptstyle p}
\end{matrix} \\
\bar{Q}_I &=
\overset{\raise3pt\hbox{$\scriptstyle \hspace*{.8em} p \hspace*{1.3em} r-p$}}{
\begin{pmatrix}
0                  &0 \\
\hspace*{4pt} Q_p  &0
\end{pmatrix}}
\hspace*{3pt}
\begin{matrix}
{\scriptstyle N+1-2r-p} \\
{\scriptstyle p}
\end{matrix},
\end{align*}
where
\begin{equation}
\label{eq:qp}
\nonumber
Q_p = 
\begin{pmatrix}
d_p \\
\vdots &\ddots \\
d_1    &\ldots &d_p
\end{pmatrix}, \quad \in \mathbb{R}^{p\times p}.
\end{equation}
We have deliberately used the notation $\bar{d}_j \equiv d_{-j} = -d_j, 1 \leq j \leq p$ to better 
convey the band structure of $Q_I$.  The matrices $-\tilde{Q}^T_I$ and $\bar{Q}_I$ provide the
head and tail of $Q_I$ thus ensuring that (\ref{eq:stencil}) holds for all points in the 
interior region.  Furthermore, (\ref{eq:diffop}), (\ref{eq:partsum}) and (\ref{eq:innerprod}) 
imply that
\begin{equation}
\label{eq:qL}
\nonumber
Q_L = 
\begin{pmatrix}
-1/2       &q_{01}     &\ldots  &q_{0,r-1} \\
-q_{01}    &0          &\ldots  &q_{1,r-1} \\
\vdots     &           &\ddots  &\vdots    \\
-q_{0,r-1} &-q_{1,r-1} &\dots   &0
\end{pmatrix}.
\end{equation}

It remains to prove that $D_R$ defined in (\ref{eq:antireflected}) is indeed given by the third 
block row of (\ref{eq:diffop}):
\begin{align*}
D_R \equiv -J_rD_LJ_{N+1} 
&= -\frac{1}{h}J_rH^{-1}_L
\begin{pmatrix}
Q_L &\tilde{Q}_I &0
\end{pmatrix}
\begin{pmatrix}
      &           &J_{r} \\
      &J_{N+1-2r}        \\
J_{r}
\end{pmatrix} \\ 
&=\frac{1}{h}H^{-\tau}_L
\begin{pmatrix}
0 &-J_{r}\tilde{Q}_IJ_{N+1-2r} &-Q^\tau_L 
\end{pmatrix}.
\end{align*}
But
\begin{align*}
J_{r}\tilde{Q}_IJ_{N+1-2r} &= 
\begin{pmatrix}
      &J_{p} \\
J_{r-p}       
\end{pmatrix} 
\begin{pmatrix}
0                &0 \\
\hspace*{3pt}Q_p &0
\end{pmatrix}
\begin{pmatrix}
              &J_{p} \\
J_{N+1-2r-p}       
\end{pmatrix} \\ \\ &= 
\overset{\raise3pt\hbox{$\scriptstyle N+1-2r-p \hspace*{2.0em} p$ \hspace*{3.5em}}}{
\begin{pmatrix}
0 &&J_pQ_pJ_p \\
0 &&0
\end{pmatrix}}
\begin{matrix}
{\scriptstyle p} \\
{\scriptstyle r-p}
\end{matrix}. 
\end{align*}
Since $Q_p$ is constant along its diagonals, it follows that $J_pQ_pJ_p = Q^T_p$, i.~e.,
\[
J_{r}\tilde{Q}_IJ_{N+1-2r} = \bar{Q}^T_I,
\]
which proves the claim.  

The main result in \cite{ks:fefdhpde,ks:oeeeda} can now be expressed as

\begin{thm}
\label{thm:sbp}
Given the interior $2p$th order accurate difference stencil (\ref{eq:stencil}), there exists a norm
$H$ (\ref{eq:norm}) and a $(2p-1)$th order accurate difference operator $D$ (\ref{eq:diffop}) such 
that $D$ satisfies a summation-by-parts rule (\ref{eq:partsum}) with respect to the scalar 
product (\ref{eq:innerprod}) provided $r$ is sufficiently large.
\end{thm}

\begin{rem}
\label{rem:sbp}
There are $r^2$ unknowns in total:  $H_L$ brings $0.5r(r+1)$ elements $h_{ij}$.  Similarly,
$Q_L$ provides $0.5r(r-1)$ unknowns $q_{ij}$.  These parameters are determined by the consistency requirements
(\ref{eq:acc}).  Typically, the resulting system of equations must be solved numerically or by using 
symbolic computation.  There are numerous examples in scientific literature.  \hfill$\Box$
\end{rem}

Next, we will derive an additional property of the non-trivial part of $H$.  
The definition of $D$ (\ref{eq:partsum}) imposes stringent conditions upon the scalar 
product (\ref{eq:innerprod}), which the following theorem shows.

\begin{thm}
\label{thm:constraint}
Let $D \in \mathbb{R}^{(N+1)\times(N+1)}$ be a consistent approximation of
$\partial / \partial x$ satisfying a summation-by-parts rule (\ref{eq:partsum}).  Then
\begin{equation}
\label{eq:constraint}
({\bf 1}, {\bf 1})_H = 1.
\end{equation}
\end{thm}

\noindent
{\bf Proof:} By (\ref{eq:partsum}):
\[
({\bf 1}, Dx)_H = 1\cdot Nh - 1\cdot 0 - (D{\bf 1}, x)_H.
\]
Since $D$ is a consistent approximation of $\partial / \partial x$ it follows that
\begin{align*}
Dx       & = {\bf 1} \\
D{\bf 1} & = {\bf 0}.
\end{align*}
Hence,
\[
({\bf 1}, {\bf 1})_H = 1,
\]
where we also used $Nh = 1$. \hfill$\Box$

\begin{thm}
\label{thm:normsum}
Let $H \in \mathbb{R}^{(N+1) \times (N+1)}$ be normalized as
\begin{equation}
\label{eq:normsum1}
\nonumber
H \equiv h 
\begin{pmatrix}
H_L \\
& I  \\
&& H_L^\tau
\end{pmatrix}, 
\quad H_L \in \mathbb{R}^{r \times r},
\end{equation}
where $h=1/N$ is the mesh size.  Furthermore, assume that $H$ is subject to 
(\ref{eq:constraint}).  Then
\begin{align*}
\sum^{r-1}_{i,j = 0}h_{ij} & = r - \frac{1}{2}.
\end{align*}
\end{thm}

\noindent
{\bf Proof:} 
According to the definition of $H$:
\begin{align*}
({\bf 1}, {\bf 1})_H 
& = h \left[ \sum^{r-1}_{i,j=0}h_{ij} + \sum^{r-1}_{i,j=0}h_{r-1-i,r-1-j} + N+1-2r \right] \\
& = h \left[2\sum^{r-1}_{i,j=0}h_{ij} + N+1-2r \right].
\end{align*}
Hence,
\[
\sum^{r-1}_{i,j = 0}h_{ij} = r - \frac{1}{2},
\]
where we used (\ref{eq:constraint}) and $hN=1$.  This proves the theorem. \hfill$\Box$

\subsection{The solution state space $V$}
\label{sec:state}
In the beginning of this section we defined grid functions $u$, $v$ as members of $\mathbb{R}^{N+1}$
where
\begin{equation}
\label{eq:euclidsp}
\nonumber
(u,v) = \sum_ju_jv_j
\end{equation}
is the usual Euclidean scalar product in $\mathbb{R}^{N+1}$.  From this scalar product, a second 
scalar product $(\cdot, \cdot)_H$ was established in $\mathbb{R}^{N+1}$.  We will now change 
perspective and regard the pair $[\mathbb{R}^{N+1}, (\cdot, \cdot)_H]$ as an inner product space
in its own right:

\begin{define}
\label{def:innerprod}
Let the inner product space $V$ be a real vector space with the inner product
\[
(\cdot,\cdot):V\times V \rightarrow \mathbb{R}
\]
for all vectors $u,v\in V$ given by (\ref{eq:uv}) and where
\[
(u,v) \equiv (u,v)_H.
\]
The inner product $(\cdot,\cdot)_H$ is defined in (\ref{eq:innerprod}). \hfill$\Box$
\end{define}

\begin{rem}
\label{rem:innerprod}
The notation $(\cdot,\cdot)$ is context dependent.  If $u,v$ are members of the inner product 
space $V$, then  $u$ and $v$ are interpreted as grid functions in $\mathbb{R}^{N+1}$ with 
$(u,v) = (u,v)_H = x^THy$.  On the other hand, if $u$ and $v$ belong to the inner product space 
$\mathbb{R}^{N+1}$, then $(\cdot,\cdot)$ denotes the usual Euclidean scalar product $(u,v) = x^Ty$ .
The vector space $V$ will be known as the {\em solution state} space, or {\em state} space for 
short.\hfill$\Box$
\end{rem}

Interpreting $u, v$ as state vectors in $V$, we conclude from (\ref{eq:partsum}) and 
Definition~\ref{def:innerprod} that 
\begin{equation}
\label{eq:operatorpartsum}
(u, Dv) \equiv u_Nv_N - u_0v_0 - (Du, v) .
\end{equation}
Hence, we regard $D$ as a linear mapping $D:V\rightarrow V$ of the inner product space 
$V$ into itself.  With this change of perspective, we can now apply the machinery 
developed in Sections~\ref{sec:adjoint},\ref{sec:pseudo} to the summation-by-parts 
operators $D$.

\section{Initial-boundary value problems}
\label{sec:ibvp}
Many problems in physics can be described by initial-boundary value problems (IBVP).  In the one-dimensional case
we will use the symbolic notation
\begin{align}
u_t(x,t) + Q\left(\partial\right)u(x,t) &= 0, \quad 0 < x < 1, \quad t > 0 \label{eq:ibvp} \\
L_0\left(\partial\right)u(0,t) &= g_0(t) \nonumber \\
L_1\left(\partial\right)u(1,t) &= g_1(t) \nonumber \\
u(x,0) &= f(x), \nonumber
\end{align}
where $u \in \mathbb{R}^d$.  The boundary operators $L_0$ and $L_1$ are assumed to be such that $Q$ 
is semibounded.  Examples include strongly hyperbolic and incompletely parabolic systems, 
cf.~\cite{kl:ibvpnse, gko:tdpdm}.  At this point, we will not be concerned with the specifics of the 
boundary operators $L_0$ and $L_1$.  They will be examined in more detail in later sections.  In general,
$L_0$ and $L_1$ will involve the state $u$ as well as derivatives of $u$.

\subsection{The boundary state space $V_\Gamma$}
\label{sec:bstate1d}
The state vectors $u, v \in \mathbb{R}^{(N+1)d}$ are defined as in (\ref{eq:uv}), but the scalar
values $u_j, v_j$ are replaced by
\begin{equation}
\label{eq:uvd}
\nonumber
\begin{pmatrix}
u_j^{(1)} \\
\vdots \\
u_j^{(d)}
\end{pmatrix},\quad 
\begin{pmatrix}
v_j^{(1)} \\
\vdots \\
v_j^{(d)}
\end{pmatrix} \in \mathbb{R}^{d}.
\end{equation}
Let
\begin{equation}
\label{eq:bstate}
u_\Gamma \equiv 
\begin{pmatrix}
u_0 \\
u_N
\end{pmatrix}, \quad
v_\Gamma \equiv 
\begin{pmatrix}
v_0 \\
v_N
\end{pmatrix} \in \mathbb{R}^{2d},
\end{equation}
denote the boundary grid functions.  

\begin{define}
\label{def:binnerprod}
Let the inner product space $V_\Gamma$ be a real vector space where the inner product is 
given by
\[
\langle\cdot,\cdot\rangle:V_\Gamma\times V_\Gamma \rightarrow \mathbb{R}
\]
for all vectors $u_\Gamma,v_\Gamma\in V_\Gamma$ defined in (\ref{eq:bstate}) and where
\begin{equation}
\label{eq:bprod}
\nonumber
\langle u_\Gamma,v_\Gamma \rangle \equiv u^T_\Gamma v_\Gamma = u^T_0v_0 + u^T_Nv_N.
\end{equation}
The space $V_\Gamma$ equipped with the above inner product will be referred to as the boundary 
state space. \hfill$\Box$
\end{define}

The discrete boundary conditions can be expressed as
\begin{equation}
\label{eq:semibc}
Lv = g, \quad 
v \in V, \quad
g = 
\begin{pmatrix}
g_0 \\
g_1
\end{pmatrix} \in V_\Gamma.
\end{equation}
Hence, the discrete boundary operator $L$ is a mapping $L:V\rightarrow V_\Gamma$
that can be represented by the following matrix:
\begin{equation}
\label{eq:bop}
L \equiv 
\begin{pmatrix}
L_L(D) \\
L_R(D)
\end{pmatrix} \in \mathbb{R}^{2d\times(N+1)d},
\end{equation}
where $L_L(D), L_R(D)$ are the discretizations of $L_0(\partial),L_1(\partial)$.  Note that the analytic boundary 
conditions will be exactly represented by $L$ if the boundary conditions do not depend
on derivatives of $u$.  This is the case for hyperbolic systems and Dirichlet conditions for
parabolic systems.  If derivatives are involved, then the first $r$  $d$ by $d$ blocks of $L_0(D)$
will be non-zero, cf. (\ref{eq:diffop}).  The same remark applies to the last $r$ blocks
of $L_1(D)$.  

\subsection{The semidiscrete equations}
\label{sec:semi}
Following the approach in \cite{po:spps1,mo:ipm} we discretize (\ref{eq:ibvp}) as
\begin{align}
v_t + PQ(D)(Pv+(I-P)\tilde{g}) &= (I-P)\tilde{g}_t, \qquad t > 0 \label{eq:semiibvp} \\
v(0) &= f, \nonumber
\end{align}
where it is assumed that initial data satisfy the boundary conditions  $Lf = g(0)$; $Q(D)$ is a 
discretization of the semibounded analytic operator $Q(\partial)$ that satisfies
a summation-by-parts property.  This is a consequence of the structure of $Q(\partial)$ and the
properties of $D$ defined by (\ref{eq:diffop});  $P$ is the projection operator representing the analytic
boundary conditions:
\begin{equation}
\label{eq:semiproj}
\nonumber
P = I - H^{-1}L^T(LH^{-1}L^T)^{-1}L,
\end{equation} 
where it is temporarily assumed that $L$ has
full rank so that $P$ is well defined.  The data vector $\tilde{g}$, finally, is defined implicitly through
\begin{equation}
\label{eq:tildeg}
L\tilde{g} \equiv g, \quad \tilde{g} \in V.
\end{equation}
But
\begin{align*}
P &= I - H^{-1}L^T(LH^{-1}L^TH_\Gamma H_\Gamma^{-1})^{-1}L \\
  &= I - H^{-1}L^TH_\Gamma\left[L\left(H^{-1}L^TH_\Gamma\right)\right]^{-1}L \qquad [(\ref{eq:adjtransp})]\\
  &= I - L^*(LL^*)^{-1}L \\
  &= I - L^+L,
\end{align*}
where the last equality follows from the generalized Penrose conditions (\ref{eq:penrose1}) - 
(\ref{eq:penrose4}) with $T=L$ and $S=L^*(LL^*)^{-1}$.  

From now on we drop the requirement that $L$ have full rank.  The boundary projection $P$
is then {\em defined} as
\begin{equation}
\label{eq:pseudoproj}
P \equiv I - L^+L,
\end{equation}
where $L$ is given by (\ref{eq:bop}) and $L^+$ is the pseudoinverse of $L$.  Note that $L^+$ is 
uniquely defined even if $L$ is singular.  This will 
allow uniform treatment of all possible combinations of characteristic boundary conditions
for hyperbolic systems.  Explicit examples will be provided in subsequent sections.  

\begin{rem}
\label{rem:bcpseudoproj1}
The boundary operator $L$ can be interpreted as the restriction of the state space onto the boundary
state space.  Similarly, $L^+$ injects the boundary state space into the state space.  \hfill$\Box$
\end{rem}

In \cite{po:spps1} it was shown that $(I-P)(v-\tilde{g}) = 0$ is equivalent to $Lv = g$ when
$L$ has full rank.  The next theorem extends this result to the case where $L$ is singular.

\begin{thm}
\label{thm:bc}
Let $P$, $L$ and $\tilde{g}$ be defined by (\ref{eq:pseudoproj}), (\ref{eq:semibc}), (\ref{eq:bop})
and (\ref{eq:tildeg}).  Then
\[
(I-P)(v - \tilde{g}) = 0 \quad \Longleftrightarrow \quad Lv = g.
\]
\end{thm}

\noindent
{\bf Proof:} 
\noindent
Assume that $(I-P)(v-\tilde{g}) = 0$.  Then
\[
L^+Lv = L^+L\tilde{g}.
\]
Multiplying by $L$ from the left yields
\[
LL^+Lv = LL^+L\tilde{g},
\]
which according to the Moore-Penrose conditions is the same as
\[
Lv = L\tilde{g} = g,
\]
where we used (\ref{eq:tildeg}) in the last step.  The opposite implication is obvious.  \hfill$\Box$

\begin{rem}
\label{rem:tildeg}
There is no need to compute $\tilde{g}$ explicitly in (\ref{eq:semiibvp}).  In fact,
\[
(I-P)\tilde{g} = L^+L\tilde{g} = L^+g.
\]
Thus, once the pseudoinverse is known we only need to compute $L^+g$.    \hfill$\Box$
\end{rem}

\subsection{The simplified semidiscrete form}
\label{sec:semisimple}
Equation~(\ref{eq:semiibvp}) can be rewritten as
\[
Pv_t + (I-P)v_t+ PQ(Pv+(I-P)\tilde{g}) = (I-P)\tilde{g}_t,
\]
which in turn can be expressed as
\begin{align*}
(I-P)(\tilde{g}_t - v_t) &= z \\
P\left[v_t + Q(Pv+(I-P)\tilde{g})\right] &= z.
\end{align*}
Thus, $z$ is orthogonal to itself and so
\[
\|z\|^2 = (z,z) = 0 \quad \Longleftrightarrow \quad z = 0.
\]
This implies that (\ref{eq:semiibvp}) decouples into two equations for $t > 0$:
\begin{align*}
(I-P)(v_t - \tilde{g}_t) &= 0 \\
P\left[v_t + Q(Pv+(I-P)\tilde{g})\right] &= 0.
\end{align*}
It should be noted that the decoupled system is equivalent to the original formulation~(\ref{eq:semiibvp}).

From now on, $P$ is assumed to be independent of $t$.  This is not a major restriction 
from a practical standpoint, since all the examples that we will consider later on will lead to
projection operators $P$ that are piecewise constant in time.  The arguments below can then be 
applied to each time interval.  Boundary data $g(t)$ may vary with $t$, however.  

Integrating the first equation of the decoupled system yields the necessary condition
\[
(I-P)(v - \tilde{g}) = (I-P)(f - \tilde{g}(0)),
\]
since $P$ does not depend on $t$.  Thus, the boundary conditions are satisfied for $t > 0$ iff 
initial data $f$ fulfill the boundary conditions  
\[
(I-P)(f - \tilde{g}(0)) = 0 \quad \Longleftrightarrow \quad Lf = g(0).
\]
This implies that there is no need to solve the first equation, since it is known a-priori
that the boundary conditions are satisfied ($Lv=g(t), t>0$).  Let
\[
w \equiv Pv \quad \Longrightarrow \quad w = Pw \quad \Longleftrightarrow \quad Lw = 0.
\]
The second equation of the decoupled system may then be expressed as
\begin{equation}
\label{eq:altsemi}
w_t + PQ(w+L^+g) = 0, \quad t > 0,
\end{equation}
where we used $Pv = Pw = w$.  From the definition of $w$ it follows immediately that $w(0) = Pf$,
that is, $w(0)$ satisfies the homogeneous boundary conditions $Lw(0)=0$.

Conversely, suppose that $w$ solves (\ref{eq:altsemi}) with initial data $w(0) = Pf$, where 
$f$ satisfies $Lf=g(0)$ by assumption.  Hence, 
$w = Pw$ for $t \geq 0$.  Next, define
\begin{equation}
\label{eq:vdef}
v \equiv Pw + (I-P)\tilde{g} = w + L^+g, \quad t \geq 0.
\end{equation}
Differentiate the above expression for $t > 0$ and leverage (\ref{eq:altsemi}):
\[
v_t + PQ(Pw+(I-P)\tilde{g}) = (I-P)\tilde{g}_t, \quad t > 0.
\]
The definition of $v$ implies $Pv = Pw$, whence
\[
v_t + PQ(Pv+(I-P)\tilde{g}) = (I-P)\tilde{g}_t, \quad t > 0,
\]
which is identical to (\ref{eq:semiibvp}).
Initial conditions:
\[
v(0) = Pw(0) + (I-P)\tilde{g}(0) = Pf + (I-P)\tilde{g}(0) = f,
\]
where we used the assumption $Lf=g(0)$ in the last step.  This implies that $v$ fulfills 
$Lv = g(t), t\geq 0$.  The following theorem has thus been proved:

\begin{thm}
\label{thm:altsemi}
Let $P$ be a time independent boundary projection defined by~(\ref{eq:pseudoproj}).
If the initial data $f$ satisfies the boundary condition $Lf = g(0)$, then
the semidiscrete approximation (\ref{eq:semiibvp}) and the simplified semidiscrete
approximation (\ref{eq:altsemi}) are equivalent.  
\end{thm}

The final version of the simplified semidiscrete approximation of (\ref{eq:ibvp}) can
now be expressed as
\begin{align}
w_t + PQ(D)(v) &= 0, \qquad t > 0 \label{eq:simplesemiibvp} \\
w(0) &= Pf, \nonumber
\end{align}
where $v$ is defined by (\ref{eq:vdef}); initial data is assumed to satisfy $Lf = g(0)$.

\subsection{Consistency of the semidiscrete approximation}
\label{sec:consistency}
To prove that (\ref{eq:simplesemiibvp}) is a consistent approximation of (\ref{eq:ibvp})
we need to introduce some additional notation.  Let
\begin{equation}
\label{eq:uqu}
\nonumber
u \equiv 
\begin{pmatrix}
u(0,t) \\
\vdots \\
u(1,t)
\end{pmatrix}, \qquad
Q(u) \equiv
\begin{pmatrix}
Q(\partial)u(0,t) \\
\vdots \\
Q(\partial)u(1,t)
\end{pmatrix}
\in \mathbb{R}^{(N+1)d}, 
\end{equation}
where $u(x,t) \in \mathbb{R}^d$ is a solution to the IBVP (\ref{eq:ibvp}). 
Consistency means that the residual
\begin{equation}
\label{eq:residual}
\nonumber
R(h) \equiv (Pu)_t + PQu
\end{equation}
is small.

\begin{rem}
\label{rem:consistency}
Consistency is obtained by formally substituting $u$ instead of $v$ in (\ref{eq:simplesemiibvp}).
Intuitively, $v$ should be an approximation of $u$ and since
\[
w = Pw = Pv,
\]
it follows that $w_t$ should approximate of $(Pu)_t$.
\end{rem}\hfill$\Box$

\begin{thm}
\label{thm:consistency}
The semidiscrete system (\ref{eq:simplesemiibvp}) is a consistent approximation of 
(\ref{eq:ibvp}).
\end{thm}

\noindent
{\bf Proof:}
According to definition of the residual $R(h)$:
\begin{align*}
R(h) &= (Pu)_t + PQu  \quad [P \mbox{\ const in time}] \nonumber \\
	 &= P(u_t + Qu)   \nonumber \\	
     &= P(u_t + Q(u) + \tilde{r}_\Omega) \quad [Q(u) = -u_t] \nonumber \\	
	 &= P\tilde{r}_\Omega(h).
\end{align*}
The vector $\tilde{r}_\Omega(h) \in \mathbb{R}^{(N+1)d}$ represents the truncation error of 
the discrete semibounded operator $Q$ at every grid point $x_j$, $j = 0, \ldots, N$.

Furthermore, since $u(x,t)$ satisfies the analytic boundary conditions: 
\[
Lu = g + r_\Gamma(h), 
\]
where $r_\Gamma(h)\in \mathbb{R}^{2d}$ is the truncation error of the discrete boundary
operator $L$.  This shows that (\ref{eq:simplesemiibvp}) is a consistent approximation 
of (\ref{eq:ibvp}).  \hfill$\Box$

\section{Boundary conditions and the pseudoinverse}
\label{sec:bc}
From the discussions in the previous sections, it is clear that the discretized boundary
condition $Lv=0$ can be implemented as a projection
\[
P = I - L^+L = I - L^*(LL^*)^+L,
\]
where $L:V\rightarrow V_\Gamma$, which implies $L^* = H^{-1}L^TH_\Gamma$.  Thus, in general
there is a dependency on the inner products of $V$ and $V_\Gamma$.  Suppose that $(LL^*)^+ = (LL^*)^{-1}$, then
\begin{align}
P = I - L^*(LL^*)^{-1}L &= I - H^{-1}L^TH_\Gamma\left[(LH^{-1}L^T)H_\Gamma\right]^{-1}L \nonumber \\
                        &= I - H^{-1}L^T\left[LH^{-1}L^T\right]^{-1}L, \label{eq:psimplified1}
\end{align}
i.~e., if $L$ has full rank, the boundary projection is independent of the inner product in $V_\Gamma$.
The above expression for $P$ is identical to that in \cite{po:spps1}, where it is assumed that $L$
has full rank.  Since $(AB)^+=B^+A^+$ is true only in special situations \cite{tg:gimp,aa:rmpp}, 
we cannot immediately extend (\ref{eq:psimplified1}) to cases where $L$ is rank deficient.  

\subsection{The simplified projection $P$}
\label{sec:psimplified}
We will now show that a result similar to (\ref{eq:psimplified1}) holds even if $L$ does not have full rank.  To prove stability, one must have:
\begin{equation}
\label{eq:pconditions}
Pv = v, \quad HP = P^TH,
\end{equation}
where $v$ is a solution of (\ref{eq:semiibvp}).  The first requirement ensures that the boundary conditions are fulfilled, cf.~Theorem~\ref{thm:bc}.  The second constraint will be used when proving stability using the energy method.  But $Pv = v$ and $P^TH = HP$ follow from the two Moore-Penrose conditions:
\[
LL^+L = L, \quad HL^+L = \left[L^+L\right]^TH.
 \]

Let $V_{\Gamma_i}$, $i=1,2$, denote the boundary space $V_\Gamma$ equipped with two different inner products
corresponding to $H_{\Gamma_i}$.  Consider the following mappings:
\[
L:V \rightarrow V_{\Gamma_1}, \quad L:V \rightarrow V_{\Gamma_2},
\]
with the pseudoinverses $L^+(H_{\Gamma_1})$ and $L^+(H_{\Gamma_2})$.  Note that $L^+(\cdot)$ denotes (possible) 
functional dependency on $H_{\Gamma_i}$, not matrix multiplication.  They will by necessity obey
the above Moore-Penrose conditions.  These pseudoinverses give rise to two projections $P_1$ and $P_2$, 
both of which satisfy
\[
P_iv = v \quad \Longleftrightarrow \quad LL^+(H_{\Gamma_i})Lv = 0 \quad \Longleftrightarrow \quad Lv = 0.
\]
Hence, $P_1$ and $P_2$ enforce identical boundary conditions.  Furthermore,  by the second Moore-Penrose condition:
\[
HP_i = P^T_iH.
\]
As a result, $P_1$ as well as $P_2$ satisfy (\ref{eq:pconditions}).    From an implementation perspective it is 
thus irrelevant which $H_\Gamma$ we choose;  the resulting projection $P$ will always lead to a consistent and 
stable implementation of the analytic boundary conditions.  In particular, we can choose $H_\Gamma=I$ when 
computing $P$, which implies (\ref{eq:psimplified1}) also in the general case.  In the case of a rank-deficient
$L$ it thus makes sense to speak of {\em a} projection, not {\em the} projection.  If $L$ has full rank, however, all 
projections reduce to a single projection as in (\ref{eq:psimplified1}).

Next, suppose that $L$ and $H$ satisfy
\[
LH = \bar{H}L
\]
for some $\bar{H} > 0$.  Thus,
\[
L^T = H^{-1}L^T\bar{H}.
\]
If we choose $H_\Gamma = \bar{H}$, then
\[
L^T = H^{-1}L^TH_\Gamma = L^*,
\]
whence
\begin{equation}
\label{eq:psimplified2}
P = I - L^+L = I - L^*\left[LL^*\right]^+L = I - L^T\left[LL^T\right]^+L.
\end{equation}

\begin{rem}
\label{rem:psimplified2}
Eq.~(\ref{eq:psimplified2}) is an example of a
situation where the algebraic expression is simplified if we choose $H_\Gamma = \bar{H}$ as opposed to the 
default $H_\Gamma = I$ .  This conclusion extends the corresponding result in \cite{po:spps1} to the case where
$L$ is rank deficient.  In fact, Proposition~\ref{prop:adjtrans} shows that $LH = \bar{H}L$ is a necessary condition
for (\ref{eq:psimplified2}) to be true.  In the next section we will examine a concrete example of where this condition holds.
\hfill $\Box$
\end{rem}

\subsection{Characteristic boundary conditions}
\label{sec:charbc1}
Let
\begin{equation}
\label{eq:charbc1}
L = 
\begin{pmatrix}
L_L \\
L_R
\end{pmatrix} = 
\begin{pmatrix}
L_0 &0 &\ldots &0 &0 \\
0   &0 &\ldots &0 &L_1 
\end{pmatrix}, \quad L_0, L_1 \in \mathbb{R}^{d\times d}.
\end{equation}
be the matrix representation of the characteristic boundary conditions for the hyperbolic system 
(\ref{eq:ibvp}).  $L_0$  and $L_1$ are the analytic boundary operators.

Up until this point no assumptions have been made on the structure of $H$.  Let $H$ be
a restricted full norm:
\begin{equation}
\label{eq:rfnorm}
H \equiv 
\begin{pmatrix}
h_{00}I \\
&\tilde{H} \\
&&h_{NN}I
\end{pmatrix}\in \mathbb{R}^{(N+1)d\times(N+1)d}.
\end{equation}
Then
\begin{equation}
\label{eq:lhpseudocommute}
LH = \bar{H} L \quad \Longleftrightarrow \quad L^T = H^{-1} L^T\bar{H},
\end{equation}
where
\begin{equation}
\label{eq:hgamma}
\nonumber
\bar{H} \equiv 
\begin{pmatrix}
h_{00}I \\
&h_{NN}I
\end{pmatrix}\in \mathbb{R}^{2d\times 2d}.
\end{equation}
Define $H_\Gamma=\bar{H}$. Thus, by (\ref{eq:psimplified2}):
\[
P = I - L^T\left[LL^T\right]^+L.
\]

\begin{prop}
\label{prop:simplify}
Let $L_i:\mathbb{R}^d\rightarrow \mathbb{R}^d$, $i=0,1$ be the analytic boundary operators
of (\ref{eq:ibvp}) with pseudoinverses $L^+_i$, $i=1,2$.  Define 
$L:V\rightarrow V_\Gamma$ as in (\ref{eq:charbc1}).  If $H$ is a restricted full 
norm (\ref{eq:rfnorm}), then
\begin{equation}
\label{eq:lpseoduinv}
\nonumber
L^+ =
\begin{pmatrix}
L^+_L &L^+_R
\end{pmatrix} \equiv
\begin{pmatrix}
L^+_0 &0 \\
0 &0 \\
\vdots &\vdots \\
0 &0 \\
0 &L^+_1
\end{pmatrix}
\end{equation}
for any
\begin{equation}
\label{eq:hgammageneral}
H_\Gamma \equiv 
\begin{pmatrix}
h_LI \\
&h_RI
\end{pmatrix} > 0, \quad I\in\mathbb{R}^{d\times d}.
\end{equation}
\end{prop}

\noindent
{\bf Proof}:  The four Moore-Penrose conditions
\begin{align*}
LL^+L &= L \\
 L^+LL^+ &= L^+ \\
H\left[L^+L\right] &= \left[L^+L\right]^TH \\ 
H_\Gamma\left[LL^+\right] &= \left[LL^+\right]^TH_\Gamma
\end{align*}
follow immediately from the structure of $L$, $L^+$, $H$ and $H_\Gamma$.  The details are left
as an exercise. \hfill$\Box$

\begin{rem}
\label{rem:simplified2}
According to Proposition~\ref{prop:simplify}, $LL^+$ is self-adjoint with respect to {\em any} 
norm $H_\Gamma$ of the form (\ref{eq:hgammageneral}).  It is true for
\[
H_\Gamma = 
\begin{pmatrix}
I \\
&I
\end{pmatrix}  \quad \mbox{and} \quad
H_\Gamma \equiv 
\begin{pmatrix}
h_{00}I \\
&h_{NN}I
\end{pmatrix}.
\]
The former expression represents the choice of scalar product in Definition~(\ref{def:binnerprod}).  
Proposition~\ref{prop:simplify} is a stronger result than (\ref{eq:psimplified2}), since $L^+$ is 
independent of $H$ and $H_\Gamma$, yet it fulfills the Moore-Penrose conditions for any 
$H_\Gamma$ (\ref{eq:hgammageneral}), not just for $H_\Gamma=\Bar{H}$.  For restricted full norms,
$H$ and $H_\Gamma$ are completely decoupled from one another as far as the Moore-Penrose
conditions are concerned.
\hfill $\Box$
\end{rem}

\subsection{Scalar advection equation}
\label{sec:advection}
Consider the hyperbolic model problem 
\begin{align}
\label{eq:standardmodel}
\nonumber
u_t + c(t)u_x & = 0,    \qquad 0 < x < 1, \quad t > 0 \\
      u(x, 0) & = f(x). \nonumber
\end{align}
The boundary conditions are defined as
\begin{align*}
\delta_0u(0, t) & = 0 \\
\delta_1u(1, t) & = 0,
\end{align*}
where
\[
\delta_0 = 
\left\{
\begin{matrix}
1 & \mbox{if\ } c(t) > 0 \\
0 & \mbox{if\ } c(t) \leq 0
\end{matrix}
\right.,
\qquad
\delta_1 = 
\left\{
\begin{matrix}
1 & \mbox{if\ } c(t) < 0 \\
0 & \mbox{if\ } c(t) \geq 0
\end{matrix}
\right..
\]
No boundary conditions are prescribed if $c(t) = 0$.  The model problem is 
discretized as
\begin{align}
\label{eq:semidiscrete}
v_t + c(t)PDPv & = 0, \quad t > 0 \\
          v(0) & = f, \nonumber
\end{align}
where $D:V\rightarrow V$ satisfies (\ref{eq:operatorpartsum}) for some restricted full norm $H$.  Define 
the discrete boundary operator $L:V \rightarrow V_\Gamma$:
\begin{equation}
\label{eq:bcop}
\nonumber
L \equiv 
\begin{pmatrix}
\delta_0 & 0 & \ldots & 0 & 0 \\
0        & 0 & \ldots & 0 & \delta_1    
\end{pmatrix}.
\end{equation}
In particular, if $c(t) = 0$, then
\[
\delta_0 = \delta_1 = 0,
\]
which means that $L=0$, i.~e., no boundary conditions are imposed in complete agreement with 
the analytic problem (\ref{eq:model}).  All prerequisites of Proposition~\ref{prop:simplify} are met.  Hence,
\[
L^+ = L^T(LL^T)^+.
\]
This leads to
\begin{equation}
\label{eq:bcoppseudoinv}
\nonumber
L^+ =
\begin{pmatrix}
\delta_0 & 0 \\
0        & 0 \\
\vdots   & \vdots \\
0        & 0 \\
0        & \delta_1   
\end{pmatrix}.
\end{equation}
Hence, by (\ref{eq:pseudoproj}):
\[
P = I -
\begin{pmatrix}
\delta_0 \\
& 0 \\
&& \ddots \\
&&& 0 \\
&&&& \delta_1
\end{pmatrix}.
\]
The projection $P$ thus varies with $t$ (piecewise constant in time) reflecting the time dependent 
boundary conditions.
This is in complete agreement with the analytic boundary conditions (prescription of 
ingoing characteristics).  Note that if $c(t) = 0$, then $P=I$ leaves the difference 
operator unchanged, that is to say, no boundary conditions are imposed.  This is consistent 
with the analytic problem, since $c(t) = 0$ implies that (\ref{eq:model}) reduces to
an ordinary differential equation.

Stability of (\ref{eq:semidiscrete}) is immediate:
\begin{align*}
\frac{d}{dt} \|v\|  & = 2(v, v_t) \\
					  & = -2c(t)(Pv, DPv) \\
					  & = c(t)\left[\left(Pv\right)^2_0 - \left(Pv\right)^2_N\right] \\
					  & \leq 0,
\end{align*}
where the last inequality is a consequence of the construction of $P$.  For this simple 
model problem one could of course have constructed the projection
directly.  The point is, (\ref{eq:pseudoproj}) is valid for any boundary operator
$L:V\rightarrow V_\Gamma$ for any scalar product $(\cdot, \cdot)$ that correpsponds to 
a restricted full norm.  The resulting projection will always lead to an energy estimate 
if one adheres to the pattern used in (\ref{eq:semidiscrete}).

\subsection{The heat equation}
\label{sec:heat}
The simplest parabolic example is furnished by the heat equation:
\begin{align}
\label{eq:heat}
\nonumber
u_t  & = u_{xx},    \qquad 0 < x < 1, \quad t > 0 \\
u(x, 0) & = f(x). \nonumber
\end{align}
The boundary conditions correspond to an adiabatic wall at $x=0$ and  $x=1$:
\begin{align*}
u_x(0, t) & = 0 \\
u_x(1, t) & = 0.
\end{align*}
Hence, the discrete boundary conditions become
\begin{align*}
(Dv)_0 & = 0 \\
(Dv)_N & = 0,
\end{align*}
where $D$ is defined by (\ref{eq:diffop}) and (\ref{eq:bstencil}) for some diagonal norm $H$.  Thus,
we define the boundary operator $L:V\rightarrow V_\Gamma$:
\begin{equation}
\label{eq:heatbcop}
\nonumber
L \equiv 
\begin{pmatrix}
d_0 &\ldots &d_{s-1} &0 &\ldots &0 &0        &\ldots &0\\
0   &\ldots &0       &0 &\ldots &0 &-d_{s-1} &\ldots &-d_0    
\end{pmatrix},
\end{equation}
where $H_\Gamma$ is taken to be the 2 by 2 identity matrix; 
$d_j \equiv d_{0j}$ for notational convenience.  The pseudoinverse becomes
\begin{equation}
\label{eq:heatpseudoinv}
\nonumber
L^+ =
\begin{pmatrix}
\hat{d}_0     &0 \\
\vdots        &\vdots \\
\hat{d}_{s-1} &0 \\    
0             &0 \\
\vdots        &\vdots \\ 
0             &0 \\
0             &-\hat{d}_{s-1} \\
\vdots        &\vdots \\ 
0             &-\hat{d}_0 
\end{pmatrix}, \qquad
\hat{d}_i \equiv \frac{d_i}{h_i}\left(\sum_{j=0}^{s-1}\frac{d_j^2}{h_j}\right)^{-1},
\end{equation}
where $h_j$ are the coefficients of the diagonal norm.  In this case $LL^+=I$, which means 
that the Penrose conditions reduce to the single requirement $HL^+L = (L^+L)^TH$.

\subsection{2 by 2 hyperbolic systems}
\label{sec:2by2}
Next we consider a diagonal 2 by 2 system, $u^T \equiv (u^{(1)} \ u^{(2)})$:
\begin{align*}
u_t + \Lambda u_x & = 0,    \qquad 0 < x < \infty, \quad t > 0, \\
          u(x, 0) & = f(x), \nonumber
\end{align*}
subject to characteristic boundary conditions
\begin{align*}
L_1u(0,t) &\equiv \delta_1\left(u^{(1)} -c_{12}(1-\delta_2)u^{(2)}\right) = 0 \\
L_2u(0,t) &\equiv \delta_2\left(-c_{21}(1-\delta_1)u^{(1)} + u^{(2)}\right) = 0,
\end{align*}
where
\[
\delta_j = 
\left\{
\begin{matrix}
1 & \mbox{if\ } \lambda_j > 0 \\
0 & \mbox{if\ } \lambda_j \leq 0
\end{matrix}
\right.,
\qquad
j=1,2.
\]
The analytic boundary conditions can now be expressed as
\[
L_0u = 0, \qquad
L_0 \equiv
\begin{pmatrix}
L_1 \\
L_2
\end{pmatrix} \in \mathbb{R}^{2\times 2}.
\]
Note that $L_iL_j^T = 0$ whenever $i \neq j$.  
This simplifies the construction of the pseudoinverse significantly:
\begin{align*}
L^+_0 = L^T_0(L_0L^T_0)^+ & = 
\begin{pmatrix} L_1^T &L_2^T \end{pmatrix}
\begin{pmatrix}
L_1L^T_1 & L_1L^T_2 \\
L_2L^T_1 & L_2L^T_2
\end{pmatrix}^+ \\
& =
\begin{pmatrix} L_1^T &L_2^T \end{pmatrix}
\begin{pmatrix}
\delta_1 \|L_1\|^{-2} \\
                             &\delta_2 \|L_2\|^{-2}
\end{pmatrix},
\end{align*}
where
\[
\|L_i\|^2 = 1 + c^2_{ij}(1-\delta_j)^2, \qquad i=1,2, \quad j\neq i.
\]
The Penrose conditions are satisfied ($\delta_j^2 = \delta_j, (1-\delta_j)^2 = 1-\delta_j$), 
in particular:
\[
L_0L^+_0 = 
\begin{pmatrix}
\delta_1 \\
        &\delta_2
\end{pmatrix}.
\]

If one were to consider the 2 by 2 system on the unit interval $0 < x < 1$ instead of the 
half space $0 < x < \infty$ one would get a boundary operator $L_N \in \mathbb{R}^{2\times 2}$ 
at $x=1$ completely analogous to $L_0$.  Thus, in the semidiscrete case the corresponding 
boundary operator $L:V\rightarrow V_\Gamma$ is given by
\begin{equation}
\label{eq:2by2bcop}
L \equiv 
\begin{pmatrix}
L_0 &0 &\ldots &0 \\
0   &0 &\ldots &L_N
\end{pmatrix}.
\end{equation}
Since the block rows of $L$ are orthogonal it follows that the pseudoinverse is given by
\begin{equation}
\label{eq:2by2pseudoinv}
L^+ =
\begin{pmatrix}
L^+_0  & 0 \\
0      & 0 \\
\vdots & \vdots \\
0      & L^+_N   
\end{pmatrix}.
\end{equation}
The verification of the Penrose conditions (\ref{eq:penrose1}) - (\ref{eq:penrose4}) is 
straightforward and is left as an exercise.

\subsection{3 by 3 hyperbolic systems}
\label{sec:3by3}
The 2 by 2 example from the previous section can be extended to the 3 by 3 case, 
$u^T \equiv (u^{(1)} \ u^{(2)} \ u^{(3)})$:
\begin{align*}
u_t + \Lambda u_x & = 0,    \qquad 0 < x < \infty, \quad t > 0 \\
          u(x, 0) & = f(x) \nonumber
\end{align*}
subject to characteristic boundary conditions
\begin{align*}
L_1u(0,t) &\equiv \delta_1\left(u^{(1)} - c_{12}(1-\delta_2)u^{(2)} - c_{13}(1-\delta_3)u^{(3)}\right) = 0 \\
L_2u(0,t) &\equiv \delta_2\left(-c_{21}(1-\delta_1)u^{(1)} + u^{(2)} - c_{23}(1-\delta_3)u^{(3)}\right) = 0 \\
L_3u(0,t) &\equiv \delta_3\left(-c_{31}(1-\delta_1)u^{(1)}  - c_{32}(1-\delta_2)u^{(2)} + u^{(3)}\right) = 0,
\end{align*}
where
\[
\delta_j = 
\left\{
\begin{matrix}
1 & \mbox{if\ } \lambda_j > 0 \\
0 & \mbox{if\ } \lambda_j \leq 0
\end{matrix}
\right.,
\qquad
j=1,2,3.
\]
Superficially this looks like the 2 by 2 example, but $L_1, L_2$ and $L_3$ are not
mutually orthogonal in general:
\begin{enumerate}
\item All $\lambda_i > 0$ implies all $\delta_i = 1$ (Supersonic inflow): $L_iL_j^T = 0$, $i\neq j$.
\item Two $\lambda_i > 0$ positive implies one $\delta_k = 0$ (Subsonic inflow): $L_iL_j^T \neq 0$ 
for $i,j \neq k$, $i\neq j$. 
\item One $\lambda_i > 0$ positive implies two $\delta_i = 0$ (Subsonic outflow): $L_iL_j^T = 0$, $i\neq j$.
\item All $\lambda_i \leq 0$ positive implies all $\delta_i = 0$ (Supersonic outflow): $L_iL_j^T = 0$, $i\neq j$. 
\end{enumerate}
Clearly, for subsonic inflow orthogonality of the boundary operators $L_i$ is lost.  The technique
in the 2 by 2 case for computing the pseudoinverse does not carry over to the 3 by 3 case.  We have the
following result, however \cite{tg:sapm,aa:rmpp}.

\begin{prop}
\label{prop:pseudoiter}
Let $L\in\mathbb{R}^{m \times n}$ be a given matrix partitioned as follows:
\[
L = 
\begin{pmatrix}
L_1 \\
\vdots \\
L_m
\end{pmatrix}, \quad
L_j\in\mathbb{R}^{1\times n}, \quad j=1,\ldots, m.
\]
Define $\tilde{L}_j\in\mathbb{R}^{j\times n}$:
\begin{equation}
\label{eq:pseudoiter1}
\nonumber
\tilde{L}_j \equiv 
\begin{pmatrix}
\tilde{L}_{j-1} \\
L_j
\end{pmatrix}, \quad j=2,\ldots, m, \quad \tilde{L}_1 \equiv L_1.
\end{equation}
Then
\begin{equation}
\label{eq:pseudoiter2}
\nonumber
\tilde{L}^+_j = 
\begin{pmatrix}
\left(I - K^T_jL_j\right)\tilde{L}^+_{j-1} &K^T_j
\end{pmatrix}, \quad j=2,\ldots, m, \quad \tilde{L}^+_1 = L^+_1,
\end{equation}
where
\[
K_j = \left\lbrace
\begin{array}{ll}
L_j\left(I-\tilde{L}^+_{j-1}\tilde{L}_{j-1}\right)/\lambda^2_j 
&\mbox{if } L_j \neq L_j\tilde{L}^+_{j-1}\tilde{L}_{j-1} \\ \\
L_j\tilde{L}^+_{j-1}\left[\tilde{L}^+_{j-1}\right]^T\!\!/(1 + \mu^2_j) 
&\mbox{otherwise}
\end{array}
\right.
\]
with
\[
\lambda_j \equiv \|L_j\left(I-\tilde{L}^+_{j-1}\tilde{L}_{j-1}\right)\|, \quad
\mu_j \equiv \|L_j\tilde{L}^+_{j-1}\|.
\]
\end{prop}

\noindent
{\bf Proof}: We note that $L=\tilde{L}_m$.  Hence, Proposition~\ref{prop:pseudoiter} provides
an iterative method for computing the pseudoinverse of an arbitrary matrix $L$.  We will
use induction over $j$ and the Moore-Penrose conditions to prove the claims.  

The result is obviously true for $j=1$:
\[
L^+_1 = \left(L^T_1L_1\right)^+L^T_1 = L^T_1\left(L_1L^T_1\right)^+.
\]
Before continuing, we will collect some results that will simplify the algebraic expressions
that follow.  If $L_j \neq L_j\tilde{L}^+_{j-1}\tilde{L}_{j-1}$, then
\begin{equation}
\label{eq:aux1}
\tilde{L}_{j-1}K^T_j = 0, \quad L_jK^T_j = 1,
\end{equation}
where we used definitions of $K_j$ and $\lambda_j$ together with the induction hypotheses
\begin{align*}
\tilde{L}_{j-1}\tilde{L}^+_{j-1}\tilde{L}_{j-1} &= \tilde{L}_{j-1}  \\
(I-\tilde{L}^+_{j-1}\tilde{L}_{j-1})^2 &= I-\tilde{L}^+_{j-1}\tilde{L}_{j-1} \\
(I-\tilde{L}^+_{j-1}\tilde{L}_{j-1})^T &= I-\tilde{L}^+_{j-1}\tilde{L}_{j-1}. 
\end{align*}
Similarly, if $L_j = L_j\tilde{L}^+_{j-1}\tilde{L}_{j-1}$, then by the definition of $K_j$ 
and the induction hypothesis $\tilde{L}^+_{j-1}\tilde{L}_{j-1}\tilde{L}^+_{j-1} = \tilde{L}^+_{j-1}$: 
\begin{equation}
\label{eq:aux2}
\tilde{L}^+_{j-1}\tilde{L}_{j-1}K^T_j = K^T_j.
\end{equation}
The following relations will be used frequently throughout the remainder of the proof:
\begin{align}
\tilde{L}^+_j\tilde{L}_j &= \left(I - K^T_jL_j\right)\tilde{L}^+_{j-1}\tilde{L}_{j-1} + K^T_jL_j \label{eq:aux3}\\
&= \tilde{L}^+_{j-1}\tilde{L}_{j-1} +  K^T_jL_j\left(I - \tilde{L}^+_{j-1}\tilde{L}_{j-1}\right). \nonumber
\end{align}

We begin by verifying the first Moore-Penrose condition $\tilde{L}_j\tilde{L}^+_j\tilde{L}_j=\tilde{L}_j$
when $L_j\neq L_j\tilde{L}^+_{j-1}\tilde{L}_{j-1}$:
\[
\tilde{L}_j\tilde{L}^+_j\tilde{L}_j = 
\begin{pmatrix}
\tilde{L}_{j-1} \\
L_j
\end{pmatrix}
\begin{pmatrix}
\tilde{L}^+_{j-1}\tilde{L}_{j-1} + K^T_jL_j\left[I- \tilde{L}^+_{j-1}\tilde{L}_{j-1}\right]
\end{pmatrix}  = \tilde{L}_j,
\]
where we used (\ref{eq:aux1}).  If $L_j = L_j\tilde{L}^+_{j-1}\tilde{L}_{j-1}$, then by (\ref{eq:aux3}):
\[
\tilde{L}^+_j\tilde{L}_j = \tilde{L}^+_{j-1}\tilde{L}_{j-1}.
\]
Hence, by the definition of $\tilde{L}_j$:
\[
\tilde{L}_j\tilde{L}^+_j\tilde{L}_j=
\begin{pmatrix}
\tilde{L}_{j-1} \\
L_j
\end{pmatrix}
\tilde{L}^+_{j-1}\tilde{L}_{j-1} = 
\tilde{L}_j,
\]
which concludes the verification of the first Moore-Penrose condition.

Next, consider the second Moore-Penrose condition $\tilde{L}^+_j\tilde{L}_j\tilde{L}^+_j$ subject
to the constraint $L_j\neq L_j\tilde{L}^+_{j-1}\tilde{L}_{j-1}$:
\begin{align*}
\tilde{L}^+_j\tilde{L}_j\tilde{L}^+_j &= 
\begin{pmatrix}
\left[I - K^T_jL_j\right]\tilde{L}^+_{j-1}\tilde{L}_{j-1} + K^T_jL_j 
\end{pmatrix}
\begin{pmatrix}
\left[I-K^T_jL_j\right]\tilde{L}^+_{j-1} &K^T_j
\end{pmatrix} \\
&= 
\begin{pmatrix}
\left[I - K^T_jL_j\right]\tilde{L}^+_{j-1}\tilde{L}_{j-1}\left[I-K^T_jL_j\right]\tilde{L}^+_{j-1} &K^T_j
\end{pmatrix} \\
&= 
\begin{pmatrix}
\left[I - K^T_jL_j\right]\tilde{L}^+_{j-1}\tilde{L}_{j-1}\tilde{L}^+_{j-1} &K^T_j
\end{pmatrix} \\
&= \tilde{L}^+_j,
\end{align*}
where we have used (\ref{eq:aux1}) repeatedly.  We have already shown that 
\[
\tilde{L}^+_j\tilde{L}_j = \tilde{L}^+_{j-1}\tilde{L}_{j-1},
\]
in case $L_j = L_j\tilde{L}^+_{j-1}\tilde{L}_{j-1}$.  Thus,
\begin{align*}
\tilde{L}^+_j\tilde{L}_j\tilde{L}^+_j &= \tilde{L}^+_{j-1}\tilde{L}_{j-1}
\begin{pmatrix}
\left[I-K^T_jL_j\right]\tilde{L}^+_{j-1} &K^T_j
\end{pmatrix} \\
&=
\begin{pmatrix}
\left[\tilde{L}^+_{j-1}\tilde{L}_{j-1}-K^T_jL_j\right]\tilde{L}^+_{j-1} &K^T_j
\end{pmatrix} \\
&= \tilde{L}^+_j,
\end{align*}
where we also used (\ref{eq:aux2}).  This shows that $\tilde{L}^+_j$ satisfies 
the second Moore-Penrose condition.

To verify the third Moore-Penrose condition under the assumption that 
$L_j\neq L_j\tilde{L}^+_{j-1}\tilde{L}_{j-1}$, we observe that (\ref{eq:aux3}) implies
\[
\tilde{L}^+_j\tilde{L}_j = 
\tilde{L}^+_{j-1}\tilde{L}_{j-1} + K^T_jL_j\left(I - \tilde{L}^+_{j-1}\tilde{L}_{j-1}\right) =
\tilde{L}^+_{j-1}\tilde{L}_{j-1} + \lambda^2_jK^T_jK_j,
\]
i.~e., $\tilde{L}^+_j\tilde{L}_j$ is symmetric.  As before, 
$\tilde{L}^+_j\tilde{L}_j = \tilde{L}^+_{j-1}\tilde{L}_{j-1}$ if $L_j = L_j\tilde{L}^+_{j-1}\tilde{L}_{j-1}$.
This completes the verification of the third Moore-Penrose condition.

The validation of the fourth Moore-Penrose is different, since we cannot leverage (\ref{eq:aux3})
anymore.  Instead, we form
\begin{equation}
\label{eq:aux4}
\nonumber
\tilde{L}_j\tilde{L}^+_j =
\begin{pmatrix}
\tilde{L}_{j-1}\left( I - K^T_jL_j \right)\tilde{L}^+_{j-1} &\tilde{L}_{j-1}K^T_j \\
L_j\left( I - K^T_jL_j \right)\tilde{L}^+_{j-1} &L_jK^T_j
\end{pmatrix} \equiv
\begin{pmatrix}
A &B \\
C &D
\end{pmatrix}.
\end{equation}
We will demonstrate that $A^T=A$, $B=C^T$ and $D^T=D$.  For $L_j\neq \tilde{L}_{j-1}\tilde{L}^+_{j-1}L_j$
these claims follow immediately from (\ref{eq:aux1}):
\[
\tilde{L}_j\tilde{L}^+_j = 
\begin{pmatrix}
\tilde{L}_{j-1}\tilde{L}^+_{j-1} &0^T \\
0 &1
\end{pmatrix}.
\]
Next, consider  $L_j = \tilde{L}_{j-1}\tilde{L}^+_{j-1}L_j$.  By the definition of $K_j$:
\begin{equation}
\label{eq:aux5}
D \equiv L_jK^T_j = \frac{1}{1+\mu^2_j}L_j\tilde{L}^+_{j-1}\left[ \tilde{L}^+_{j-1} \right]^TL^T_j 
= \frac{\mu^2_j}{1+\mu^2_j}.
\end{equation}
Thus,
\begin{equation}
\label{eq:aux6}
\nonumber
C \equiv L_j\left(I-K^T_jL_j\right)\tilde{L}^+_{j-1} = \frac{1}{1+\mu^2_j}L_j\tilde{L}^+_{j-1}.
\end{equation}
Also
\begin{align}
B \equiv \tilde{L}_{j-1}K^T_j &= 
\frac{1}{1+\mu^2_j}\tilde{L}_{j-1}\tilde{L}^+_{j-1}\left[ \tilde{L}^+_{j-1} \right]^TL^T_j \nonumber \\
&=
\frac{1}{1+\mu^2_j}\left[\tilde{L}_{j-1}\tilde{L}^+_{j-1}\right]^T\left[ \tilde{L}^+_{j-1} \right]^TL^T_j \nonumber \\
&=
\frac{1}{1+\mu^2_j}\left[ \tilde{L}^+_{j-1} \right]^TL^T_j = C^T. \label{eq:aux7}
\end{align}
Finally,
\begin{align}
A \equiv \tilde{L}_{j-1}\left( I - K^T_jL_j \right)\tilde{L}^+_{j-1} &=
\tilde{L}_{j-1}\tilde{L}^+_{j-1} - \tilde{L}_{j-1}K^T_jL_j\tilde{L}^+_{j-1} \nonumber \\
&=
\tilde{L}_{j-1}\tilde{L}^+_{j-1} - \frac{1}{1+\mu^2_j}\left[L_j\tilde{L}^+_{j-1}\right]^TL_j\tilde{L}^+_{j-1} \label{eq:aux8},
\end{align}
where we used (\ref{eq:aux7}) in the last equality.  From (\ref{eq:aux5}), (\ref{eq:aux7}) and
(\ref{eq:aux8}) it is clear that $A^T=A$, $B=C^T$ and $D^T=D$, whence the fourth Moore-Penrose condition 
holds.  This concludes the proof. \hfill $\Box$

After this detour we can return to the original problem of finding the pseudoinverse of the boundary operator
$L_0:\mathbb{R}^{3} \rightarrow \mathbb{R}^{3}$ representing the characteristic boundary conditions 
of the hyperbolic 3 by 3 system.  This is done by applying 
Proposition~\ref{prop:pseudoiter} to $L_0$ with $m=n=3$ thus resulting in $L^+_0=\tilde{L}^+_3$.  Hence, for 
restricted full norms we recover (\ref{eq:2by2bcop}) and (\ref{eq:2by2pseudoinv}) also in the 3 by 3 case.

\subsection{$d$ by $d$ hyperbolic systems}
\label{sec:dbyd}
The general one-dimensional hyperbolic IBVP is given by:
\begin{align*}
u_t + Au_x & = 0,    \qquad 0 < x < 1, \quad t > 0 \\
          u(x, 0) & = f(x),
\end{align*}
where $A\in \mathbb{R}^{d\times d}$ is assumed to be symmetric.  Furthermore,
the boundary conditions are assumed to be well-posed:
\begin{equation}
\label{eq:wellposed}
\nonumber
L_0u(0,t) = 0, \quad L_1u(1,t) = 0, \quad L_0, L_1 \in \mathbb{R}^{d \times d}.
\end{equation}
The rank of $L_0$, $L_1$ is $d_0 \leq d$ and $d_1 \leq d$; $d_0, d_1$ correspond 
to the number of ingoing characteristics at the respective boundary.
Except for special cases, one should not expect to find a closed formula for the 
pseudoinverses of $L_0$ and $L_1$.  Instead, one can use the singular value
decomposition \cite{gvl:matcomp} to obtain the pseudoinverse:
\begin{equation}
\label{eq:svd}
\nonumber
L_0 = U\Sigma V^T, \qquad U, V \in \mathbb{R}^{d\times d} \mbox{\ \ orthogonal}.
\end{equation}
The matrix $\Sigma$ is diagonal with $d_0$ positive elements.  The pseudoinverse can
then be expressed as
\begin{equation}
\label{eq:pseudosvd}
L^+_0 = V\Sigma^+U^T.
\end{equation}
Direct computations show that $L_0$ and $L^+_0$ satisfy the Penrose conditions.  Thus,
$L^+_0$ defined as above is indeed the pseudoinverse of $L_0$.  A similar factorization 
can be derived for $L_1$.

\begin{rem}
\label{rem:svd}
Since any matrix has an SVD decomposition, it follows that (\ref{eq:pseudosvd}) is 
valid for any $L_0$. \hfill$\Box$
\end{rem}

\noindent
For restricted full norms $H$, the discrete boundary operator 
$L:\mathbb{R}^{(N+1)d} \rightarrow \mathbb{R}^{2d}$ will have the same block structure
as (\ref{eq:2by2bcop}), which implies that $L^+$ is given by (\ref{eq:2by2pseudoinv}).

\subsubsection{Pseudoinverses and full norms}
\label{sec:pseudofullnorm}
We now shift focus to general norms $H$: 
\[
H = (h_{ij}I)_{0\leq i,j \leq N}, \qquad I \in \mathbb{R}^{d\times d}.
\]
The discrete boundary conditions $L:V \rightarrow V_\Gamma$ ($H_\Gamma=I\in\mathbb{R}^{2d\times 2d}$) can still be represented by a matrix as in (\ref{eq:2by2bcop}):
\[
L =
\begin{pmatrix}
L_0 &0 &\ldots &0 &0\\
0   &0 &\ldots &0 &L_1
\end{pmatrix} \in \mathbb{R}^{2d\times (N+1)d}.
\]
Since $H$ is SPD we can find a lower triangular Cholesky factor $G$ such that 
$GG^T = H$.  Let $U\Sigma V^T$ be the SVD of $L\left[G^{T}\right]^{-1}$, i.~e.:
\[
L = U\Sigma V^TG^T, \quad
U \in \mathbb{R}^{2d\times 2d}, \quad
\Sigma \in \mathbb{R}^{2d\times (N+1)d}, \quad
V \in \mathbb{R}^{(N+1)d\times (N+1)d}.
\]
Define
\[
S \equiv \left[G^{T}\right]^{-1}V\Sigma^+U^T \in \mathbb{R}^{(N+1)d\times 2d}.
\]
The generalized Moore-Penrose conditions ($H_2 = H_\Gamma = I$) 
\[
LSL = L, \quad SLS = S, \quad H(SL) = (SL)^TH, \quad LS = (LS)^T
\]
are readily established, which shows that $S$ is the pseudoinverse of $L$:
\begin{equation}
\label{eq:fullpseudo}
\nonumber
L^+ = \left[G^{T}\right]^{-1}V\Sigma^+U^T \in \mathbb{R}^{(N+1)d\times 2d}.
\end{equation}

\begin{rem}
\label{rem:general}
The technique described above carries
over verbatim to more general boundary conditions involving derivatives and function
values that may occur in conjunction with incompletely parabolic systems, for instance.  
In this case the boundary operator $L:V \rightarrow V_\Gamma$ 
can be expressed as
\begin{align*}
L(D) & \equiv
\begin{pmatrix}
L_L(D) \\
L_R(D)
\end{pmatrix} \\
& \equiv
\begin{pmatrix}
L_0(h) & \ldots & L_{s-1}(h) & 0 & \ldots & 0 & 0            & \ldots & 0 \\
0      & \ldots & 0          & 0 & \ldots & 0 & L_{N-s+1}(h) & \ldots & L_N(h)
\end{pmatrix},
\end{align*}
see \eqref{eq:bop}.  Apply the arguments presented earlier to $L(D)$ as defined above.  \hfill$\Box$
\end{rem}

\section{Multiblock stability}
\label{sec:multiblock} 
\begin{figure}
\centering
\includegraphics[width=7cm]{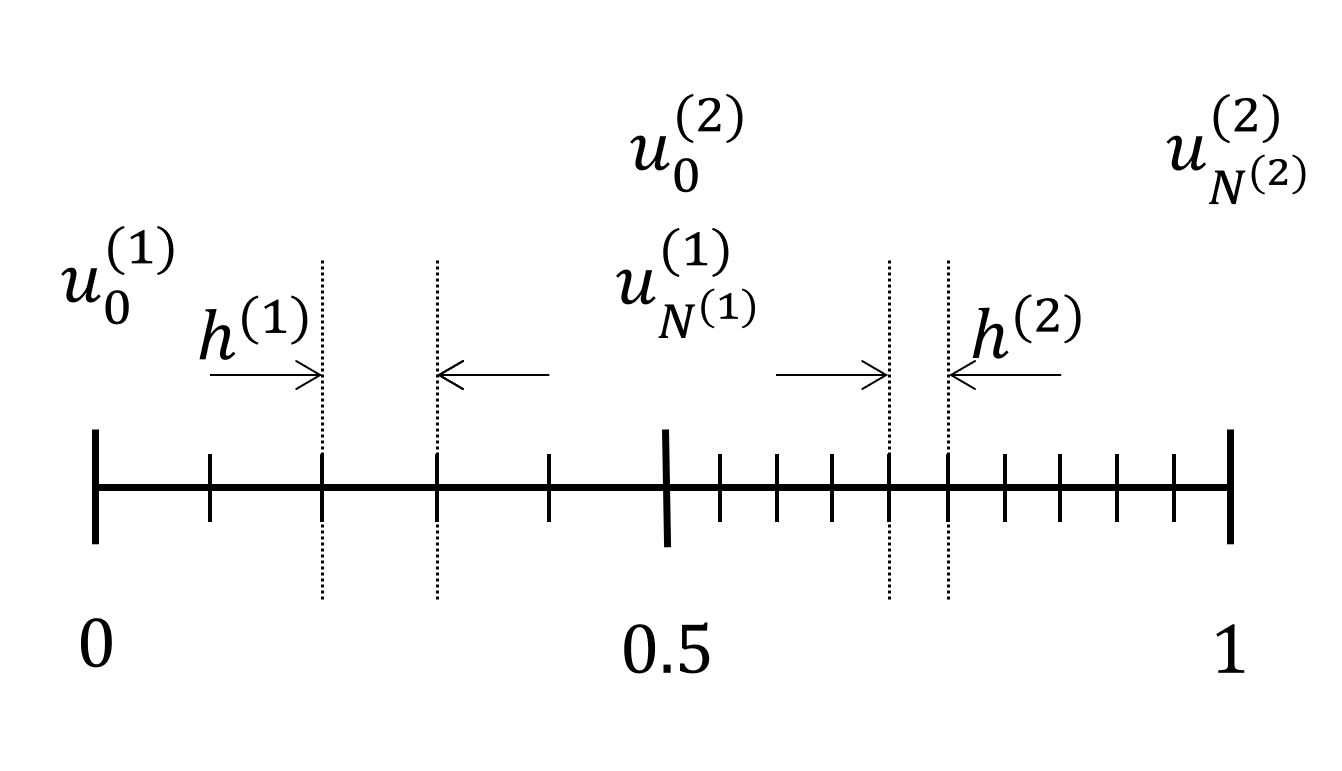}
\caption{Two uniform grids with states $u^{(1)}$ and $u^{(1)}$}
\label{fig:twogrids}
\end{figure}

We will study the $2$-grid problem in one space dimension, cf. Fig.~\ref{fig:twogrids}.
The mesh sizes are defined as
\begin{equation}
\label{eq:h1h2}
\nonumber
h^{(1)} = \frac{1}{2N^{(1)}}, \qquad
h^{(2)} = \frac{1}{2N^{(2)}}.
\end{equation}
Let
\begin{align}
\label{eq:twoint}
\Omega & \equiv \Omega_1 \cup \Omega_2 \nonumber \\
\Omega_1 & \equiv \left[0, 1/2\right] \\
\Omega_2 & \equiv \left[1/2, 1\right] \nonumber
\end{align}
with the grid points $x_j$:
\begin{align*}
\Omega_1 & : \left\{x_j = jh^{(1)}\right\}, \quad j=0,\ldots,N^{(1)} \\
\Omega_2 & : \left\{x_j = 1/2 + jh^{(2)}\right\}, \quad j=0,\ldots,N^{(2)}.
\end{align*}

To set the stage for the stability analysis that will follow, we introduce some additional 
terminology.  The notion of a {\em multiset} \cite{wb:mt} is a useful concept in many situations.  The canonical example is the representation of an integer in terms of its prime factors.  Contrary to a regular set, in a multiset an element can occur more than once.  Each element $x$ in a multiset $A$ is associated with a multiplicity $m(x)$.  Set
inclusion, union, intersection, and sum are extended to multisets through the multiplicity function:
\begin{itemize}
\item Inclusion: $A \subset B$ if $m_A(x) \leq m_B(x)$
\item Union: $C = A \cup B$ where $m_C(x) = \max(m_A(x), m_B(x))$
\item Intersection: $C = A \cap B$ where $m_C(x) = \min(m_A(x), m_B(x))$
\item Sum: $C = A + B \equiv A \cup B$ where $m_C(x) = m_A(x) + m_B(x)$
\end{itemize}

We can now define the multiset
\begin{equation}
\label{eq:iplus}
\nonumber
\Omega_+ \equiv \Omega_1 + \Omega_2.
\end{equation}
This implies that the multiplicity of $x=1/2$ is 2; the multiplicity of all other grid points is
one.  For future reference we define
\begin{equation}
\label{eq:n1n2}
\nonumber
N \equiv N^{(1)} + N^{(2)}.
\end{equation}
For each of the subintervals $\Omega_i, i=1,2$, there is an associated state
space $V_i$ in the sense of Definition~\ref{def:innerprod} with the corresponding
state vectors $u^{(i)},v^{(i)}$ and inner products $(\cdot,\cdot)_i$.

\subsection{Multiblock scalar products}
\label{sec:multiscalprod}
We begin with the following
\begin{define}
\label{def:auginnerprod}
Let the inner product space $V_+$ be a real vector space with the inner product
\[
(\cdot,\cdot)_+:V_+\times V_+ \rightarrow \mathbb{R}
\]
for all vectors $u^{(+)},v^{(+)}\in V_+$ and where
\begin{equation}
\label{eq:scalprodplus}
\nonumber
(u^{(+)},v^{(+)})_+ \equiv (u^{(1)},v^{(1)})_1 + (u^{(2)},v^{(2)})_2.
\end{equation}
The space $V_+$ will be referred to as the augmented state space. \hfill$\Box$
\end{define}

\begin{rem}
\label{rem:auginnerprod}
From the definition of $(u^{(+)},v^{(+)})_+$ it follows immediately that
\[
(u^{(+)},v^{(+)})_+ = \left[u^{(+)}\right]^TH^{(+)}v^{(+)},
\]
where
\begin{equation}
\label{eq:hplus}
H^{(+)} =
\overset{\raise3pt\hbox{$\scriptstyle N^{(1)}+1 \ \ N^{(2)}+1$}}{
\begin{pmatrix}
H^{(1)}\\
&H^{(2)}
\end{pmatrix}}
\begin{matrix}
{\scriptstyle N^{(1)}+1} \\
{\scriptstyle N^{(2)}+1}
\end{matrix}, 
\end{equation}
and where
\begin{equation}
\label{eq:uvplus}
\nonumber
u^{(+)} \equiv
\begin{pmatrix}
u^{(1)} \\
u^{(2)}
\end{pmatrix}, \quad
v^{(+)} \equiv
\begin{pmatrix}
v^{(1)} \\
v^{(2)}
\end{pmatrix} \in \mathbb{R}^{N+2}.
\end{equation}
It is clear that $u^{(+)}$ and $v^{(+)}$ are well-defined grid vectors on the 
multiset $\Omega_+$, since $x=1/2$ has multiplicity two corresponding to 
$u^{(1)}_{N^{(1)}}$ and $u^{(2)}_0$.  Hence, the augmented state space $V_+$ is
well defined.  Augmented state spaces will be defined for higher dimensions in
Section~\ref{sec:2dim}. \hfill$\Box$
\end{rem}

\subsubsection{The embedding operator $E$}
\label{sec:embedding}
Let $u, v\in\mathbb{R}^{N+1}$ be grid vectors on $\Omega=\Omega_1\cup\Omega_2$ (\ref{eq:twoint}):
\begin{equation}
\label{eq:uv1}
u \equiv 
\begin{pmatrix}
u_0 \\
\vdots \\
u_N
\end{pmatrix}, \quad
v \equiv 
\begin{pmatrix}
v_0 \\
\vdots \\
v_N
\end{pmatrix} \in \mathbb{R}^{N+1},
\end{equation}
where we recall that $N=N^{(1)}+N^{(2)}$.  Note the formal similarity of this definition 
with (\ref{eq:uv}).  It should be pointed out that every point in $\Omega$ has
multiplicity one.

Next, define a mapping $E:\mathbb{R}^{N+1} \rightarrow \mathbb{R}^{N+2}$:
\begin{equation}
\label{eq:epart}
E = 
\begin{pmatrix}
E^{(1)} \\
E^{(2)}
\end{pmatrix},
\end{equation}
where
\begin{align}
E^{(1)} &\equiv 
\overset{\raise3pt\hbox{$\scriptstyle \ 0 \ \quad \quad \ N^{(1)} \ N^{(1)}+1 \quad N$}}
{\left(
\begin{array}{rcc|cccc}
1  &\ &\ &0  &\ldots &0 \\
&\ddots &\ &\vdots  &&\vdots \\
&&1 &0 &\ldots &0 \\
\end{array}
\right)}
\begin{array}{c}
\scriptstyle 0 \\ \\ \\
\scriptstyle N^{(1)}
\end{array} \quad \in\mathbb{R}^{(N^{(1)}+1)\times (N+1)} \label{eq:e1embed} \\
\nonumber \\
E^{(2)} &\equiv 
\overset{\raise3pt\hbox{$\scriptstyle 0 \ \ \quad \quad N^{(1)} \ N^{(1)}+1 \quad N$}}
{\left(
\begin{array}{rcc|cccc}
0 &\ldots &1 &0 \\
0 &\ldots &0 &1 \\
\vdots &&\vdots &&\ddots \\
0 &\ldots &0 &&&1
\end{array}
\right)}
\begin{array}{c}
\scriptstyle 0 \\
\scriptstyle 1 \\ \\ \\
\scriptstyle N^{(2)}
\end{array} \quad \in\mathbb{R}^{(N^{(2)}+1)\times (N+1)}. \label{eq:e2embed}
\end{align} 
Partition $I^{(i)}\in\mathbb{R}^{(N^{(i)}+1)\times(N^{(i)}+1)}$:
\begin{equation}
\label{eq:idpart}
\nonumber
I^{(1)} = 
\begin{pmatrix}
\tilde{I}^{(1)} \\
&1
\end{pmatrix}, \quad
I^{(2)} = 
\begin{pmatrix}
1 \\
&\tilde{I}^{(2)}
\end{pmatrix}, \quad \tilde{I}^{(i)} \in \mathbb{R}^{N^{(i)}\times N^{(i)}}, \quad i=1,2.
\end{equation}
The embedding $E$ can then be alternatively written as:
\begin{equation}
\label{eq:eexplicit}
\nonumber
E = 
\begin{pmatrix}
\tilde{I}^{(1)} &0 &0 \\
0 &1 &0 \\
0 &1 &0 \\
0 &0 &\tilde{I}^{(2)}
\end{pmatrix}.
\end{equation}
This partition will be helpful when discussing the structure of the adjoint and pseudoinverse 
of the embedding operator $E$ in later sections.

Define an inner product on $\mathbb{R}^{N+1}\times\mathbb{R}^{N+1}$, $N=N^{(1)}+N^{(2)}$, as follows:
\begin{equation}
\label{eq:innerprodembed}
(u,v)_H \equiv (Eu)^TH^{(+)}Ev.
\end{equation}

\begin{prop}
\label{prop:scalprodembed}
The scalar product (\ref{eq:innerprodembed}) is well defined.
\end{prop}

\noindent
{\bf Proof}:
Bilinearity and symmetry are obvious.  Positivity is also a straightforward consequence:
\[
0 = (u,u)_H = (Eu)^TH^{(+)}Eu \Longleftrightarrow Eu = 0
\Longleftrightarrow u = 0.
\]
The first equivalence follows from $H^{(+)}>0$; the second equivalence holds since $E$ 
is one-to-one. \hfill$\Box$

We are now ready to define the state space $V$ for the grid vectors $u, v$ (\ref{eq:uv1})
defined on $\Omega=\Omega_1\cup\Omega_2$:

\begin{define}
\label{def:innerprodembed}
Let the inner product space $V$ be a real vector space with the inner product
\[
(\cdot,\cdot):V\times V \rightarrow \mathbb{R}
\]
for all vectors $u,v\in V$ given by (\ref{eq:uv1}) and where
\[
(u,v) \equiv (u,v)_H.
\]
The inner product $(\cdot,\cdot)_H$ is defined in (\ref{eq:innerprodembed}). \hfill$\Box$
\end{define}

\begin{rem}
\label{rem:innerprodembed}
The single-domain definition~\ref{def:innerprod} carries over {\em verbatim} to the multi-domain 
definition~\ref{def:innerprodembed}.  Only the inner products $H$ differ as defined by 
(\ref{eq:innerprod}) and (\ref{eq:innerprodembed}).  The inner product on $V$ is obtained by 
restricting $(\cdot,\cdot)_+$ of Definition~\ref{def:auginnerprod} to
${\cal R}(E) \subset \mathbb{R}^{N+2}$.  The vectors
\[
u^{(e)} = Eu, \quad v^{(e)} = Ev \in \mathbb{R}^{N+2},
\]
are the embeddings of $u, v \in \mathbb{R}^{N+1}$ into $\mathbb{R}^{N+2}$.  The 
embeddings $u^{(e)}, v^{(e)}$ satisfy
\[
u^{(e)}_{N^{(1)}}=u^{(e)}_{N^{(1)}+1}, \quad v^{(e)}_{N^{(1)}}=v^{(e)}_{N^{(1)}+1},
\]
by construction.  \hfill$\Box$
\end{rem}

Completely analogous to interpreting $D$ as a linear operator $D:V\rightarrow V$, we will
from now on consider the {\em embedding} $E$ as an operator $E:V\rightarrow V_+$.  The following results are direct consequences of the operator definition of $E$.

\begin{prop}
\label{prop:isometry}
The embedding operator $E:V\rightarrow V_+$ is an isometry.
\end{prop}

\noindent
{\bf Proof}:
Follows immediately from Definition~\ref{def:innerprodembed} and (\ref{eq:innerprodembed}).\hfill$\Box$

\begin{prop}
\label{prop:embedadjoint}
Given $E:V\rightarrow V_+$, then $E^* = H^{-1}E^TH^{(+)}$.
\end{prop}

\noindent
{\bf Proof}:
Follows immediately from Definition~\ref{def:adjoint}.\hfill$\Box$

\begin{prop}
\label{prop:embedpseudo}
Given $E:V\rightarrow V_+$, then $E^+ = E^*$.
\end{prop}

\noindent
{\bf Proof}:
By Proposition~\ref{prop:embedadjoint} and (\ref{eq:innerprodembed}):
\[
E^*E = \left[H^{-1}E^TH^{(+)}\right]E = H^{-1}\left[E^TH^{(+)}E\right] = I.
\]
Thus, the Moore-Penrose conditions reduce to the single requirement
\[
\left[EE^*\right]^* = EE^*,
\]
which is trivially satisfied. \hfill$\Box$

\subsubsection{Structure of $H$}
\label{sec:hstructure}
From Definition~\ref{def:innerprodembed} and (\ref{eq:innerprodembed}):
\begin{equation}
\label{eq:h}
H = E^TH^{(+)}E.
\end{equation}
Partition $E$:
\begin{equation}
\label{eq:epartalt}
\nonumber
E = \left(
\begin{array}{c|c}
I &0 \\
\hline
J & \Pi
\end{array}
\right), \qquad
\begin{array}{l}
I \in \mathbb{R}^{(N^{(1)}+1)\times(N^{(1)}+1)} \\
0 \in \mathbb{R}^{(N^{(1)}+1)\times N^{(2)}} \\
J \in \mathbb{R}^{(N^{(2)}+1)\times(N^{(1)}+1)} \\
\Pi \in \mathbb{R}^{(N^{(2)}+1)\times N^{(2)}}.
\end{array}
\end{equation}
$J$ has a single nonzero element $J_{0N^{(1)}}=1$; $\Pi$ is a shift operator:
\[
\Pi\begin{pmatrix}
u_{N^{(1)}+1} \\
\vdots \\
u_{N^{(1)}+N^{(2)}}
\end{pmatrix}
=
\begin{pmatrix}
0 \\
u_{N^{(1)}+1} \\
\vdots \\
u_{N^{(1)}+N^{(2)}}
\end{pmatrix}.
\]
Hence, by (\ref{eq:hplus}) and (\ref{eq:h}):
\begin{equation}
\label{eq:hpart}
H =
\begin{pmatrix}
H^{(1)} + J^TH^{(2)}J & J^TH^{(2)}\Pi \\
\Pi^TH^{(2)}J & \Pi^TH^{(2)}\Pi
\end{pmatrix},
\end{equation}
where
\begin{align*}
J^TH^{(2)}J &= 
\begin{pmatrix}
0 &\ldots &0 \\
\vdots &&\vdots \\
0 &\ldots &H^{(2)}_{00}
\end{pmatrix} \\
J^TH^{(2)}\Pi &= 
\begin{pmatrix}
0 &\ldots &0 \\
\vdots &&\vdots \\
H^{(2)}_{01} &\ldots &H^{(2)}_{0N^{(2)}}
\end{pmatrix} \\
\Pi^TH^{(2)}J &= 
\begin{pmatrix}
0 &\ldots &H^{(2)}_{10} \\
\vdots &&\vdots \\
0 &\ldots &H^{(2)}_{N^{(2)}0}
\end{pmatrix} \\
\Pi^TH^{(2)}\Pi &= 
\begin{pmatrix}
H^{(2)}_{11} &\ldots &H^{(2)}_{1N^{(2)}} \\
\vdots &&\vdots \\
H^{(2)}_{N^{(2)}1} &\ldots &H^{(2)}_{N^{2)}N^{(2)}}
\end{pmatrix}.
\end{align*}
Substituting the above expressions in (\ref{eq:hpart}) allows us to symbolically express 
$H\in\mathbb{R}^{(N+1)\times(N+1)}$ as
\begin{equation}
\label{eq:hsymb}
H = 
\begin{pmatrix}
H^{(1)} \\
&0
\end{pmatrix} + 
\begin{pmatrix}
0 \\
&H^{(2)}
\end{pmatrix}, \qquad H^{(i)} \in \mathbb{R}^{N^{(i)}+1}, \quad i=1,2.
\end{equation}
Note that block operations are not defined since the diagonal zero block of the first matrix is 
$N^{(2)}\times N^{(2)}$;  the zero block of the second matrix is $N^{(1)}\times N^{(1)}$.  Addition is still 
defined at the element level since both block matrices are $(N+1)\times (N+1)$.  Heuristically,
$H$ is obtained from $H^{(+)}$ by shifting $H^{(2)}$ one step to the NW along the main diagonal of $H^{(+)}$
and adding overlapping elements together:
\begin{equation}
\label{eq:hoverlap1}
H=\left(
\begin{array}{ccc|ccc}
h_{00}^{(1)} &\ldots &h_{0N^{(1)}}^{(1)} \\
\vdots &&\vdots \\
h_{N^{(1)}0}^{(1)} &\ldots &\tilde{h}_{N^{(1)}N^{(1)}}^{(1)} &h_{01}^{(2)} &\ldots &h_{0N^{(2)}}^{(2)} \\
\hline
&&h_{10}^{(2)} &h_{11}^{(2)} &\ldots &h_{1N^{(2)}}^{(2)} \\
&&\vdots &\vdots && \vdots \\
&&h_{N^{(2)}0}^{(2)} &h_{N^{(2)}1}^{(2)} &\ldots &h_{N^{(2)}N^{(2)}}^{(2)}
\end{array}
\right).
\end{equation}
Alternatively, we can block $H$ based on the second matrix of (\ref{eq:hsymb}):
\begin{equation}
\label{eq:hoverlap2}
H=\left(
\begin{array}{ccc|ccc}
h_{00}^{(1)} &\ldots &h_{0,N^{(1)}-1}^{(1)} &h_{0N^{(1)}}^{(1)}\\
\vdots &&\vdots &\vdots \\
h_{N^{(1)}-1,0}^{(1)} &\ldots &h_{N^{(1)}-1,N^{(1)}-1}^{(1)} &h_{N^{(1)}-1,N^{(1)}}^{(1)} \\
\hline
h_{N^{(1)}0}^{(1)} &\ldots &h_{N^{(1)},N^{(1)}-1}^{(1)} &\tilde{h}_{00}^{(2)} &\ldots &h_{0N^{(2)}}^{(2)}   \\
&&&\vdots &&\vdots \\
&&&h_{N^{(2)}0}^{(2)} &\ldots &h_{N^{(2)}N^{(2)}}^{(2)}
\end{array}
\right).
\end{equation}
Pictorially, $H^{(1)}$ is shifted one step towards SE along the main diagonal of $H^{(+)}$:
\[ 
\tilde{h}_{N^{(1)}N^{(1)}}^{(1)} = \tilde{h}_{00}^{(2)} = h_{N^{(1)}N^{(1)}}^{(1)} + h_{00}^{(2)}.
\]

\subsubsection{Structure of $H^{-1}$}
\label{sec:hblockinv}
The inverse of $H$ will be needed when discussing difference operators that have a summation-by-parts
property.  To this end, we will derive an explicit expression for the block inverse of $H$ and then
apply the formula to (\ref{eq:hoverlap1}) and (\ref{eq:hoverlap2}).  This is a well-known result in linear algebra, but we include the derivation presented below for completeness.

Let
\begin{equation}
\label{eq:hblock}
\nonumber
H =
\begin{pmatrix}
H_{11} &H_{12} \\
H_{21} &H_{22}
\end{pmatrix}
\end{equation}
be an arbitrary blocking of $H$ with square diagonal blocks $H_{ii}$.  Since $H$ is SPD it follows 
immediately that the blocks $H_{ii}$ inherit this property.  In particular, $H_{ii}$ are nonsingular.
Make the ansatz
\[
\tilde{H} \equiv
\begin{pmatrix}
H_{11}^{-1} &-H_{11}^{-1}H_{12}H_{22}^{-1} \\
-H_{22}^{-1}H_{21}H_{11}^{-1} &H_{22}^{-1}
\end{pmatrix}.
\]
Straightforward matrix multiplication yields
\begin{align*}
H\tilde{H} &=
\begin{pmatrix}
I-H_{12}H_{22}^{-1}H_{21}H_{11}^{-1} &0 \\
0 &I-H_{21}H_{11}^{-1}H_{12}H_{22}^{-1} 
\end{pmatrix} \\
\tilde{H}H &=
\begin{pmatrix}
I-H_{11}^{-1}H_{12}H_{22}^{-1}H_{21} &0 \\
0 &I-H_{22}^{-1}H_{21}H_{11}^{-1}H_{12} 
\end{pmatrix}.
\end{align*}
We thus arrive at two equivalent expressions for $H^{-1}$:
\begin{align*}
H^{-1} &\! = \!
\begin{pmatrix}
\left[H_{11}-H_{12}H_{22}^{-1}H_{21}\right]^{-1} &\hspace{-10pt} -H_{11}^{-1}H_{12}\left[H_{22}-H_{21}H_{11}^{-1}H_{12}\right]^{-1} \\
-H_{22}^{-1}H_{21}\left[H_{11}-H_{12}H_{22}^{-1}H_{21}\right]^{-1} &\hspace{-10pt} \left[H_{22}-H_{21}H_{11}^{-1}H_{12}\right]^{-1}
\end{pmatrix} \\
H^{-1} &\! = \!
\begin{pmatrix}
\left[H_{11}-H_{12}H_{22}^{-1}H_{21}\right]^{-1} &\hspace{-10pt} -\! \left[H_{11} -H_{12}H_{22}^{-1}H_{21}\right]^{-1}\!\! H_{12}H_{22}^{-1} \\
-\left[H_{22}-H_{21}H_{11}^{-1}H_{12}\right]^{-1}H_{21}H_{11}^{-1} &\hspace{-10pt} \left[H_{22} -H_{21}H_{11}^{-1}H_{12}\right]^{-1}
\end{pmatrix}
\hspace{-3pt}
.
\end{align*}
Equivalence follows from the identities
\begin{align}
H_{11}^{-1}H_{12}\left[H_{22}-H_{21}H_{11}^{-1}H_{12}\right]^{-1} &= \left[H_{11}-H_{12}H_{22}^{-1}H_{21}\right]^{-1}H_{12}H_{22}^{-1} \label{eq:schur1} \\
H_{22}^{-1}H_{21}\left[H_{11}-H_{12}H_{22}^{-1}H_{21}\right]^{-1} &= \left[H_{22}-H_{21}H_{11}^{-1}H_{12}\right]^{-1}H_{21}H_{11}^{-1}, \label{eq:schur2} 
\end{align}
which are established by pre- and post-multiplication of $H_{22}-H_{21}H_{11}^{-1}H_{12}$ and
$H_{11}-H_{12}H_{22}^{-1}H_{21}$.  These factors are known as the Schur complements of $H_{11}$ and
$H_{22}$.  Note that both Schur complements exist since $H$ (and thus $H_{ii}$), is nonsingular.

We conclude the block inverse discussion by showing that one Schur complement can be expressed
in terms of the other.  Let $S$ be the Schur complement of $H_{11}$:
\begin{equation}
\label{eq:schur}
S \equiv H_{22}-H_{21}H_{11}^{-1}H_{12}.
\end{equation}
Hence, $\left[H_{11}-H_{12}H_{22}^{-1}H_{21}\right]^{-1}=$
\begin{align*}
& H_{11}^{-1}\left[\left(H_{11}-H_{12}H_{22}^{-1}H_{21}\right) + H_{12}H_{22}^{-1}H_{21}\right]\left[H_{11}-H_{12}H_{22}^{-1}H_{21}\right]^{-1} = \\
& H_{11}^{-1}\left[I + H_{12}H_{22}^{-1}H_{21}\left(H_{11}-H_{12}H_{22}^{-1}H_{21}\right)^{-1}\right] = \qquad [(\ref{eq:schur2}), (\ref{eq:schur})] \\
& H_{11}^{-1} + H_{11}^{-1}H_{12}S^{-1}H_{21}H_{11}^{-1}.
\end{align*}
We have thus arrived at the final representation
\begin{equation}
\label{eq:hblockinv}
H^{-1} =
\begin{pmatrix}
H_{11}^{-1} + H_{11}^{-1}H_{12}S^{-1}H_{21}H_{11}^{-1} &-H_{11}^{-1}H_{12}S^{-1} \\
-S^{-1}H_{21}H_{11}^{-1} &S^{-1}
\end{pmatrix},
\end{equation}
where the Schur complement $S$ is given by (\ref{eq:schur}).  The identities
\[
HH^{-1} = H^{-1}H = I
\]
follow immediately from (\ref{eq:hblockinv}).  Note the structural similarity between the ansatz $\tilde{H}$
and the block inverse $H^{-1}$.  The latter is formally obtained from the former by replacing $H_{22}$ with 
$S$ and by adding $H_{11}^{-1}H_{12}S^{-1}H_{21}H_{11}^{-1}$ to $H_{11}^{-1}$.

Apply (\ref{eq:hblockinv}) where the blocks $H_{ij}$ are given by (\ref{eq:hoverlap1}).  Since $H_{11}$ and its Schur complement $S$ are nonsingular, $H^{-1}$ is block diagonal iff $H_{21}$ and $H_{12}=H_{21}^T$ vanish identically.  But this happens iff
\[
h_{0j}^{(2)} = 0, \quad j=1, \ldots, N^{(2)}.
\]
Since the partition of $H$ is arbitrary, we could just as well have started with (\ref{eq:hoverlap2}) instead,
which would result in
\[
h_{N^{(1)}j}^{(1)} = 0, \quad j=0, \ldots, N^{(1)}-1.
\]
Hence, both sets of constraints must be fulfilled.  Thus,
\begin{equation}
\label{eq:rfnorm2}
H^{(2)} =
\begin{pmatrix}
h^{(2)}_{00} \\
&\tilde{H}^{(2)}
\end{pmatrix}, \quad  \tilde{H}^{(2)}\in \mathbb{R}^{N^{(2)}\times N^{(2)}},
\end{equation}
and 
\begin{equation}
\label{eq:rfnorm1}
H^{(1)} =
\begin{pmatrix}
\tilde{H}^{(1)} \\
&h^{(1)}_{N^{(1)}N^{(1)}}
\end{pmatrix},
\quad \tilde{H}^{(1)}\in \mathbb{R}^{N^{(1)}\times N^{(1)}}.
\end{equation}
Consequently, by (\ref{eq:hsymb}):
\begin{align}
H &=
\begin{pmatrix}
\tilde{H}^{(1)} \\
&h^{(1)}_{N^{(1)}N^{(1)}} + h^{(2)}_{00} \\
&&\tilde{H}^{(2)}
\end{pmatrix} \label{eq:hstruct} \\
H^{-1} &=
\begin{pmatrix}
\left[\tilde{H}^{(1)}\right]^{-1} \\
&\left[h^{(1)}_{N^{(1)}N^{(1)}} + h^{(2)}_{00}\right]^{-1} \\
&&\left[\tilde{H}^{(2)}\right]^{-1}
\end{pmatrix}. \label{eq:hinvstruct} 
\end{align}
In other words, $H^{-1}$ is block diagonal iff the norms $H^{(i)}$ are restricted full norms
\cite{po:hodmdi}.  

\subsubsection{Structure of $E^*$}
\label{sec:eadjointstructure}
Since $E^*$ involves $H^{-1}$ it will be assumed that $H^{(1)}$ and $H^{(2)}$ are
restricted full norms (\ref{eq:rfnorm1}) and (\ref{eq:rfnorm2}).  

By (\ref{eq:epart}), (\ref{eq:e1embed}) and (\ref{eq:e2embed}):
\[
E^TH^{(+)} = 
\begin{pmatrix}
\tilde{H}^{(1)} &0                &0            &0 \\
0               &h^{(1)}_{N^{(1)}N^{(1)}} &h^{(2)}_{00} &0 \\
0               &0                &0            &\tilde{H}^{(2)}
\end{pmatrix}
\]
and hence, by Proposition~\ref{prop:embedadjoint} and (\ref{eq:hinvstruct}):
\begin{equation}
\label{eq:adjstructure}
E^* = 
\begin{pmatrix}
\tilde{I}^{(1)} &0    &0      &0 \\
0   &\chi &1-\chi &0 \\
0   &0    &0      &\tilde{I}^{(2)}
\end{pmatrix}, \quad \chi \equiv \frac{h^{(1)}_{N^{(1)}N^{(1)}}}{h^{(1)}_{N^{(1)}N^{(1)}} + h^{(2)}_{00}}.
\end{equation}

\begin{rem}
\label{rem:eadjointstructure}
Suppose that the restricted full norms $H^{(1)}$ and $H^{(2)}$ satisfy the structural 
requirements of (\ref{eq:norm}), which is a very common situation in practical 
computations.  Then
\[
h^{(1)}_{N^{(1)}N^{(1)}} = \mu h^{(1)},\quad h^{(2)}_{00} = \mu h^{(2)},
\]
for some $\mu > 0$; $h^{(1)}$ and $h^{(2)}$ are the mesh sizes in the respective domains, whence
\[
\chi = \frac{h^{(1)}}{h^{(1)}+h^{(2)}}.
\]
\ \hfill$\Box$
\end{rem}

\subsection{Multiblock difference operators}
\label{sec:multidiffop}
Analogous to (\ref{eq:hplus}), define
\begin{equation}
\label{eq:dplus}
D^{(+)} \equiv
\overset{\raise3pt\hbox{$\scriptstyle N^{(1)}+1 \ \ N^{(2)}+1$}}{
\begin{pmatrix}
D^{(1)}\\
&D^{(2)}
\end{pmatrix}}
\begin{matrix}
{\scriptstyle N^{(1)}+1} \\
{\scriptstyle N^{(2)}+1}
\end{matrix},
\end{equation}
where $D^{(i)}$ satisfies the summation-by-parts rule (\ref{eq:operatorpartsum}) with respect to 
$(\cdot,\cdot)_i = (\cdot,\cdot)_{H^{(i)}}, i=1,2$.  It follows that $D:V_+ \rightarrow V_+$.
Define the multiblock difference operator $D$:  

\begin{define}
\label{def:diffopembed}
Given the inner product space $V$ as in Definition~\ref{def:innerprodembed}, the difference
operator $D:V\rightarrow V$ is defined as
\begin{equation}
\label{eq:diffopembed}
D \equiv H^{-1}E^TH^{(+)}D^{(+)}E.
\end{equation}
\ \hfill$\Box$
\end{define}
The main result can be formulated as 

\begin{prop}
\label{prop:diffopembed}
Let $D:V\rightarrow V$ be as in Definition~\ref{def:diffopembed}.  Then
\begin{itemize}
\item[(i)] $D$ satisfies summation by parts with respect to the inner product 
$(\cdot,\cdot)$ of Definition~\ref{def:innerprodembed}.
\item[(ii)] $D$ is a consistent approximation of $\partial/\partial x$.
\item[(iii)] $D=E^+D^{(+)}E$.
\end{itemize}
\end{prop}

\noindent
{\bf Proof}:  The first statement follows, since
\begin{align*}
(u,Dv) &= \left[u^{(1)}\right]^TH^{(1)}D^{(1)}v^{(1)} + \left[u^{(2)}\right]^TH^{(2)}D^{(2)}v^{(2)} \\
       &= u^{(2)}_{N^{(2)}}v^{(2)}_{N^{(2)}} - u^{(1)}_0v^{(1)}_0 - \left[D^{(1)}u^{(1)}\right]^TH^{(1)}v^{(1)} 
	                                                      - \left[D^{(2)}u^{(2)}\right]^TH^{(2)}v^{(2)},
\end{align*}
where we used $u^{(2)}_0 = u^{(1)}_{N^{(1)}}$ with a similar constraint for $v$ thanks to the 
embedding $E$:
\[
u^{(e)} = Eu =
\begin{pmatrix}
u^{(1)} \\
u^{(2)}
\end{pmatrix}.
\]
Hence
\[
(u,Dv) = u_Nv_N - u_0v_0 - (Du,v), \quad N=N_1+N_2,
\]
which proves the first assertion.

To prove the second assertion, let $u(x)$ be a smooth function on $[0,1]$.  Define
\[
u \equiv
\begin{pmatrix}
u(x_0) \\
\vdots \\
u(x_N)
\end{pmatrix}, \qquad
x_0=0, x_{N}=1.
\]
Then
\begin{align*}
D^{(+)}Eu = D^{(+)}u^{(e)} 
&= 
\begin{pmatrix}
D^{(1)}u^{(1)} \\
D^{(2)}u^{(2)} 
\end{pmatrix} \\
&=
\begin{pmatrix}
u^{(1)}_x \\
u^{(2)}_x
\end{pmatrix} + {\cal O}(h^p) \qquad \left[u^{(1)}_x[N^{(1)}] = u^{(2)}_x[0] + {\cal O}(h^p)\right] \\
&= u^{(e)}_x + {\cal O}(h^p) \\
&= Eu_x + {\cal O}(h^p).
\end{align*}
By Definition~\ref{def:diffopembed}:
\[
Du = H^{-1}E^TH^{(+)}D^{(+)}Eu = H^{-1}E^TH^{(+)}Eu_x + {\cal O}(h^p) = u_x + {\cal O}(h^p),
\]
which finishes the proof of the second assertion.

The third statement, finally, follows immediately from Propositions~\ref{prop:embedadjoint},~\ref{prop:embedpseudo}. 
This concludes the proof.\hfill$\Box$

\subsubsection{Structure of $D$}
\label{sec:dstructure}
From (\ref{eq:diffopembed}) it is clear that constructing $D$ involves computing
the inverse of $H$.  For general norms $H^{(i)}$, in particular implicit ones, this
can be a costly operation, cf.~(\ref{eq:hblockinv}) for the explicit formula for $H^{-1}$.
For explicit norms involving small blocks close to the boundaries the inverse can be 
precomputed using symbolic tools from which the numeric coefficients can be generated
reliably.

There are important classes of norms, where the computation of the leading factor $H^{-1}E^TH^{(+)}=E^*$ can be significantly simplified.  
Suppose that $H^{(i)}$ are restricted full norms.  Then $D$ as defined in (\ref{eq:diffopembed}) can be written as
\begin{equation}
\label{eq:dembedstructure}
\begin{pmatrix}
D_{00}^{(1)} &\!\! \ldots &\!\! D_{0N^{(1)}}^{(1)} &\!\! 0 &\!\! \ldots &\!\! 0 \\
\vdots &&\!\! \vdots &\!\! \vdots &&\!\! \vdots \\
D_{N^{(1)}-1,0}^{(1)} &\!\!\ldots &\!\! D_{N^{(1)}-1,N^{(1)}}^{(1)} &\!\! 0 &\!\! \ldots &\!\! 0 \\
\chi D_{N^{(1)}0}^{(1)} &\!\!\ldots &\!\! \chi D_{N^{(1)}N^{(1)}}^{(1)} + (1-\chi) D_{00}^{(2)} &\!\! (1-\chi) D_{01}^{(2)} &\!\! \ldots &\!\! (1-\chi) D_{0N^{(2)}}^{(2)}   \\
0 &\!\!\ldots &\!\! D_{10}^{(2)} &\!\! D_{11}^{(2)} &\!\! \ldots &\!\! D_{1N^{(2)}}^{(2)}  \\
\vdots &&\!\! \vdots  &\!\! \vdots &&\!\! \vdots  \\
0 &\!\!\ldots &\!\! D_{N^{(2)}0}^{(2)} &\!\! D_{N^{(2)}1}^{(2)} &\!\! \ldots &\!\! D_{N^{(2)}N^{(2)}}^{(2)}
\end{pmatrix}
\end{equation}
where we used (\ref{eq:adjstructure}).  In practice, $Du$ is computed by evaluating $D^{(i)}u^{(i)}$ in each subinterval followed by computing the weighted mean
\[
\chi\left[D^{(1)}u^{(1)}\right]_{N_1} + (1-\chi)\left[D^{(2)}u^{(2)}\right]_0.
\]
The arithmetic overhead is negligible compared to computing the difference stencils
in each subinterval.

A very common situation is that the boundary stencils of $D^{(1)}$ and $D^{(2)}$ are the "anti-reflections" of
one another, cf. (\ref{eq:antireflected}).  If this is the case, then $\chi$ becomes
\[
\chi = \frac{h^{(1)}}{h^{(1)} + h^{(2)}},
\]
see Remark~\ref{rem:eadjointstructure}.  Hence,
\[
\chi D_{N^{(1)},N^{(1)}-j}^{(1)} = \frac{1}{h_1 + h_2}d_{N^{(1)},N^{(1)}-j} = -\frac{1}{h^{(1)} + h^{(2)}}d_{0j} = -(1-\chi) D_{0j}^{(2)}.
\]
In particular, for $j=0$:
\[
\chi D_{N^{(1)}N^{(1)}}^{(1)} + (1-\chi) D_{00}^{(2)} = 0.
\]
For $j>s$, where $s$ is some positive constant (independent of $N^{(1)}$  and $N^{(2)}$), the stencil coefficients 
$d_{0j}$ are zero.  Hence, the middle row of $D$ represents an anti-symmetric difference stencil
corresponding to different mesh sizes $h^{(i)}$ to the left and right of the center point $x=1/2$.  
For $h^{(1)}=h^{(2)}=h$ the traditional centered anti-symmetric difference stencil is recovered.

\subsection{A one-dimensional model example}
\label{sec:1dmodel}
In this section we will apply the previous technique to the one-dimensional advection equation in
skew symmetric form:
\begin{align*}
u_t + \frac{1}{2}(c(x,t)u)_x + \frac{1}{2}c(x,t)u_x & = 0,    \qquad 0 < x < 1, \quad t > 0 \\
                                          u(x, 0) & = f(x).
\end{align*}
The boundary conditions are defined as
\begin{align*}
\delta_1u(0, t) & = 0 \\
\delta_2u(1, t) & = 0,
\end{align*}
where
\[
\delta_1 = 
\left\{
\begin{matrix}
1 & \mbox{if\ } c(0,t) > 0 \\
0 & \mbox{if\ } c(0,t) \leq 0
\end{matrix}
\right.,
\qquad
\delta_2 = 
\left\{
\begin{matrix}
1 & \mbox{if\ } c(1,t) < 0 \\
0 & \mbox{if\ } c(1,t) \geq 0
\end{matrix}
\right..
\]

The unit interval $[0,1]$ will be split into two halves $\Omega_1$ and $\Omega_2$ as in (\ref{eq:twoint}).
For each interval we define diagonal norms 
\begin{equation}
\label{eq:H1H2}
\nonumber
H^{(1)} =
\begin{pmatrix}
h^{(1)}_{00} \\
&\ddots \\
&&h^{(1)}_{N^{(1)}N^{(1)}} 
\end{pmatrix}, \qquad
H^{(2)} =
\begin{pmatrix}
h^{(2)}_{00} \\
&\ddots \\
&&h^{(2)}_{N^{(2)}N^{(2)}} 
\end{pmatrix},
\end{equation}
and the associated SBP operators $D^{(1)}$ and $D^{(2)}$.  Hence, by (\ref{eq:hstruct}):
\[
H = 
\begin{pmatrix}
h^{(1)}_{00} \\
&\ddots \\
&&h^{(1)}_{N^{(1)}N^{(1)}} + h^{(2)}_{00}  \\
&&&\ddots \\
&&&&h^{(2)}_{N^{(2)}N^{(2)}} 
\end{pmatrix} \in \mathbb{R}^{(N+1)\times (N+1)}.
\]
The corresponding difference operator $D$ is given by (\ref{eq:dembedstructure}).
The coefficient matrix $C$ is defined as (recall that $x_{N} = 1$)
\begin{equation}
\label{eq:c}
\nonumber
C =
\begin{pmatrix}
c(0,t) \\
&\ddots \\
&&c(x_{N},t) 
\end{pmatrix} \in\mathbb{R}^{(N+1)\times(N+1)}.
\end{equation}
The boundary conditions, finally, can be expressed as
\begin{equation}
\label{eq:bcic}
\nonumber
Lu = 0, \qquad u = 
\begin{pmatrix} u_0 &\ldots &u_{N}\end{pmatrix}^T,
\end{equation}
where
\begin{equation}
\label{eq:bcicop}
\nonumber
L = 
\begin{pmatrix}
\delta_1&\ldots &0 \\
0 &\ldots &\delta_2
\end{pmatrix} \in\mathbb{R}^{2\times(N+1)}.
\end{equation}
As usual, the projection operator $P$ becomes
\[
P = I-L^+L = I-L^TL.
\]
The semidiscrete system can then be expressed as 
\begin{align*}
u_t + P\left(\frac{1}{2}DCPu  + \frac{1}{2}CDPu\right)& = 0, \quad t > 0 \\
          u(0) & = f.
\end{align*}
Since $D$ satisfies summation by parts with respect $H$, and since $C$ and $P$ are self-adjoint with respect
to $H$ ($HP=P^TH$, $HC=C^TH$), an energy estimate follows.  The arithmetic operations of the single domain 
case carry over {\em verbatim} to the multidomain case.

\section{Two space dimensions, single-block case}
\label{sec:2dim}
Consider the unit square $\Omega = [0,1]\times[0,1]$ with grid points
\begin{equation}
\label{eq:gridpoints}
\nonumber
(x_i, y_j) \equiv (ih_1, jh_2), \quad
h_1 = 1/N_1,\quad h_2 = 1/N_2. 
\end{equation}
For future reference we define the discrete boundaries corresponding to $y=0$, $x=1$,
$y=1$ and $x=0$:
\begin{align}
\label{eq:boundarysegment}
\Gamma_1 &\equiv  \left\{ (ih_1, 0) \right\} \nonumber \\
\Gamma_2 &\equiv  \left\{ (1, jh_2) \right\} \\
\Gamma_3 &\equiv  \left\{ (ih_1, 1) \right\} \nonumber \\
\Gamma_4 &\equiv  \left\{ (0, jh_2) \right\}. \nonumber
\end{align}
The ordering of the boundary segments $\Gamma_i$ corresponds to traversing the 
boundary $\Gamma$ in the positive (counterclockwise) direction:
\begin{equation}
\label{eq:boundary}
\nonumber
\Gamma \equiv \cup_{i=1}^4\Gamma_i.
\end{equation}
Internally, the boundary segments $\Gamma_1$ and $\Gamma_2$ are ordered according 
to increasing $i$ and $j$; the ordering of $\Gamma_3$ and $\Gamma_4$ corresponds to 
decreasing $i$ and $j$.
\begin{figure}
\centering
\includegraphics[width=0.7\columnwidth]{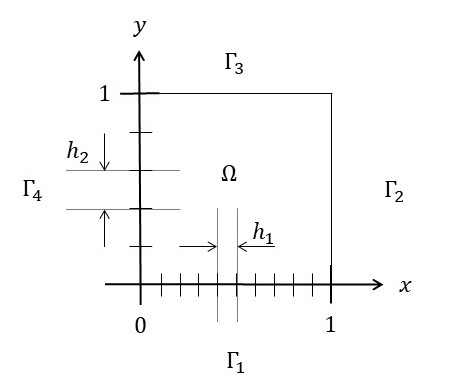}
\caption{Unit square with mesh size $h_1$ and $h_2$}
\label{fig:unitsquare}
\end{figure}
Each grid point (including the four corner points) 
\begin{equation}
\label{eq:cornerpoints}
\nonumber
(0, 0),\quad (1, 0),\quad (1, 1),\quad (0, 1),
\end{equation}
has multiplicity one.  

Analogous to the one-dimensional case we introduce
\begin{equation}
\label{eq:msetboundary}
\nonumber
\Gamma_+ = \sum^4_{i=1}\Gamma_i.
\end{equation}
Hence, the multiplicity of the four corner points is two; all other points 
have multiplicity one.  Just like the one-dimensional case, the increased multiplicity at the corner 
points is triggered by partial summation occurring twice at the same grid point.  In the one-dimensional
multidomain formulation this situation occurred at the interface between the two computational domains.
In the two-dimensional formulation the boundary state representing the corner points is 
needed twice:  once for partial summation in the $x$-direction and once in the $y$-direction.

\subsection{The solution state space $V$}
\label{sec:2dstate}
At each point $(x_i, y_j)$ we define a state variable $u_{ij}(t)$.  Arrange
the state variables into a column vector:
\begin{equation}
\label{eq:2dstate}
u \equiv
\begin{pmatrix}
u_0 \\
\vdots \\
u_{N_2}
\end{pmatrix}, \quad
u_j \equiv
\begin{pmatrix}
u_{0j} \\
\vdots \\
u_{N_1j}
\end{pmatrix} \in \mathbb{R}^{N_1 + 1},  \quad
0\leq j\leq N_2.
\end{equation}
This block structure corresponds to the usual row ordering of the state variables $u_{ij}$ 
with $i$ being the inner index and $j$ the outer one. 

Summation by parts is simplified if we also introduce an alternate representation of the 
grid function $u$ corresponding to column ordering of $u_{ij}$, where $j$ is the inner index:
\begin{equation}
\label{eq:2daltstate}
u \equiv
\begin{pmatrix}
u^{0} \\
\vdots \\
u^{N_1}
\end{pmatrix}, \quad
u^{i} \equiv
\begin{pmatrix}
u_{i0} \\
\vdots \\
u_{iN_2}
\end{pmatrix} \in \mathbb{R}^{N_2 + 1},   \quad
0\leq i\leq N_1,
\end{equation}
where we have adopted the convention of using superscripts for column reference.

The two-dimensional norm is constructed as
\begin{equation}
\label{eq:h2d}
H \equiv H_2 \otimes H_1 \in \mathbb{R}^{(N_1+1)(N_2+1)\times (N_1+1)(N_2+1)},
\end{equation}
where $H_i, i=1,2$, represent one-dimensional scalar products (\ref{eq:innerprod}).
Use $H$ to define an inner product on $\mathbb{R}^{(N_1+1)(N_2+1)\times (N_1+1)(N_2+1)}$:
\begin{equation}
\label{eq:2dinnerprod}
(u,v)_H \equiv u^THv.
\end{equation}

\begin{define}
\label{def:2dinnerprod}
Let the inner product space $V$ be a real vector space with the inner product
\[
(\cdot,\cdot):V\times V \rightarrow \mathbb{R}
\]
for all vectors $u,v\in V$ given by (\ref{eq:2dstate}) and where
\[
(u,v) \equiv (u,v)_H.
\]
The inner product $(\cdot,\cdot)_H$ is defined in (\ref{eq:2dinnerprod}). \hfill$\Box$
\end{define}

\begin{rem}
\label{rem:2dinnerprod}
Definition~\ref{def:2dinnerprod} together with definitions (\ref{eq:2dstate}) and 
(\ref{eq:2daltstate}) imply that
\[
u_j\in V_1, \quad u^i\in V_2,
\]
where $V_1$ and $V_2$ correspond to Definition~\ref{def:innerprod} with scalar
products $(u_j,v_j)_1 \equiv u_j^TH_1v_j$ and $(u^i,v^i)_2 \equiv \left[u^i\right]^TH_2v^i$.
The vector space $V$ of Definition~\ref{def:2dinnerprod} will be referred to as the 
two-dimensional solution state space.\hfill$\Box$
\end{rem}

\subsection{Summation-by-parts operators $D_x$ and $D_y$}
\label{sec:2dsbpops}
Define the two-dimensional difference operators $D_x, D_y:V\rightarrow V$
\begin{align}
D_x & \equiv I_2 \otimes D_1 \label{eq:dxdy}\\
D_y & \equiv D_2 \otimes I_1,\nonumber
\end{align}
where 
\begin{align*}
I_1, D_1&:V_1 \rightarrow V_1 \\
I_2, D_2&:V_2 \rightarrow V_2.
\end{align*}
$D_i, i=1,2,$ satisfies summation by parts with respect to $(\cdot,\cdot)_i$; 
$\otimes$ denotes the Kronecker product of two matrices.

The following notation along with some technical lemmas concerning Kronecker 
products will prove useful:

\begin{lemma}
\label{lemma:kron1}
Let $A$, $B$, $C$, $D$ be matrices such that the ordinary matrix multiplications 
$AB$ and $CD$ exist, then the mixed product satisfies
\[
(A\otimes B)(C\otimes D) = (AC)\otimes (BD).
\]
\end{lemma}

\begin{lemma}
\label{lemma:kron2}
The Kronecker product $A\otimes B$ is invertible iff $A$ and $B$ are invertible,
in which case
\[
(A\otimes B)^{-1} = A^{-1}\otimes B^{-1}.
\]
This property applies to pseudoinverses as well:
\[
(A\otimes B)^+ = A^+\otimes B^+.
\]
\end{lemma}

\begin{lemma}
\label{lemma:kron3}
Transposition is distributive over the Kronecker product:
\[
(A\otimes B)^T = A^T\otimes B^T.
\]
\end{lemma}

Let
\begin{align}
H_x &\equiv I_2 \otimes H_1 = 
\begin{pmatrix}
H_1 \\
&\ddots \\
&&H_1
\end{pmatrix} \label{eq:hxhydef} \\ 
H_y &\equiv H_2 \otimes I_1 =
\begin{pmatrix}
h^{(2)}_{00}I_1 & \ldots & h^{(2)}_{0N_2}I_1 \\
\vdots          &        & \vdots \\
h^{(2)}_{N_20}I_1 & \ldots & h^{(2)}_{N_2N_2}I_1
\end{pmatrix}. \nonumber
\end{align}
The next lemma is an immediate consequence of Lemma ~\ref{lemma:kron1}.

\begin{lemma}
\label{lemma:hxhydxdy}
Let $H_x$, $H_y$, $D_x$, $D_y$ be defined by (\ref{eq:hxhydef}) and (\ref{eq:dxdy}).
Then
\begin{align}
H_xH_y &= H_yH_x = H \nonumber \\
D_xH_y &= H_yD_x  \label{eq:hxhy} \\
D_yH_x &= H_xD_y. \nonumber
\end{align}
\end{lemma}

\noindent
{\bf Proof}: Applying Lemma~\ref{lemma:kron1} twice:
\[
H_xH_y = \left[I_2\otimes H_1\right]\left[H_2\otimes I_1\right] = H_2\otimes H_1 =
\left[H_2\otimes I_1\right]\left[I_2\otimes H_1\right] = H_yH_x.
\]
The remaining statements are proved in a similar way. \hfill$\Box$

\begin{thm}
\label{thm:2dpartsum}
Let $D_x,D_y:V\rightarrow V$ be given by (\ref{eq:dxdy}) and $(\cdot,\cdot)$ by
Definition~\ref{def:2dinnerprod}.  Then
\begin{align*}
(u,D_xv) &= (u^{N_1},v^{N_1})_2 - (u^{0},v^{0})_2 - (D_xu,v) \\
(u,D_yv) &= (u_{N_2},v_{N_2})_1 - (u_{0},v_{0})_1 - (D_yu,v).
\end{align*}
\end{thm}

\noindent
{\bf Proof:}
We have
\begin{align*}
(u,D_xv) &= u^THD_xv = u^TH_yH_xD_xv \\
           &= \sum^{N_2}_{k,l=0}h^{(2)}_{kl}\left[u^T_kH_1D_1v_l\right] = 
		      \sum^{N_2}_{k,l=0}h^{(2)}_{kl}(u_k,D_1v_l)_1 \\
		   &= \sum^{N_2}_{k,l=0}h^{(2)}_{kl}\left[u_{N_1k}v_{N_1l} - u_{0k}v_{0l} - (D_1u_k,v_l)_1\right] \\
		   &= \left[u^{N_1}\right]^TH_2v^{N_1} - \left[u^{0}\right]^TH_2v^{0}
		   - \sum^{N_2}_{k,l=0}h^{(2)}_{kl}\left[(D_1u_k)^TH_1v_l\right] \\
		   &= (u^{N_1},v^{N_1})_2 - (u^{0},v^{0})_2 
		   - \sum^{N_2}_{k,l=0}(D_xu)_k^T\left[h^{(2)}_{kl}I_1\right](H_xv)_l \\
		   &= (u^{N_1},v^{N_1})_2 - (u^{0},v^{0})_2 - (D_xu,v).
\end{align*}
Analogously,
\begin{align*}
(u,D_yv) &= u^TH_yH_xD_yv = u^TH_yD_yH_xv = (u,D_yH_xv)_{H_y} \\
         &= u^T_{N_2}H_1v_{N_2} - u^T_{0}H_1v_{0} - (D_yu,H_xv)_{H_y} \\
		 &= (u_{N_2},v_{N_2})_1 - (u_0,v_0)_1 - (D_yu)^TH_yH_xv \\
		 &= (u_{N_2},v_{N_2})_1 - (u_0,v_0)_1 - (D_yu,v),
\end{align*}
which completes the proof. \hfill$\Box$

\begin{rem}
\label{rem:2dpartsum}
Theorem~\ref{thm:2dpartsum} holds for any one-dimensional norms $H_1$ and $H_2$,
i.~e., any scalar products $(\cdot,\cdot)_1$ and $(\cdot,\cdot)_2$.  This 
includes implicit norms as long as the summation-by-parts rule 
(\ref{eq:partsum}) holds. \hfill$\Box$
\end{rem}

\subsection{The boundary state $V_\Gamma$}
\label{sec:2dbstate}
For each boundary segment $\Gamma_i$, define the grid vectors $u_{\Gamma_i}$:
\begin{align}
u_{\Gamma_1} &\equiv u_0 \nonumber \\
u_{\Gamma_2} &\equiv u^{N_1} \label{eq:segmentbstate}\\
u_{\Gamma_3} &\equiv J_1u_{N_2} \nonumber \\
u_{\Gamma_4} &\equiv J_2u^{0},\nonumber 
\end{align}
where $J_1\in\mathbb{R}^{(N_1+1)\times(N_1+1)}$ and $J_2\in\mathbb{R}^{(N_2+1)\times(N_2+1)}$ are
anti-diagonal permutation matrices, cf. (\ref{eq:adiag}).  The reason for including the permutation 
matrices is to ensure that the definition of the boundary state corresponds to the positive 
orientation of the boundary $\Gamma$ (\ref{eq:boundarysegment}).  Let
\begin{equation}
\label{eq:2dbstate}
u^{(e)}_\Gamma \equiv
\begin{pmatrix}
u_{\Gamma_1} \\
u_{\Gamma_2} \\
u_{\Gamma_3} \\
u_{\Gamma_4}
\end{pmatrix} \in \mathbb{R}^{2(N+2)}, \quad
u_\Gamma \equiv
\begin{pmatrix}
u_{\Gamma_1}[:\!\! N_1 ] \\
u_{\Gamma_2}[:\!\! N_2 ] \\
u_{\Gamma_3}[:\!\! N_1 ] \\
u_{\Gamma_4}[:\!\! N_2 ]
\end{pmatrix} \in \mathbb{R}^{2N \times 2N},
\end{equation}
represent the grid vectors on $\Gamma_+$ and $\Gamma$, $N=N_1+N_2$;  $u_{\Gamma_i}[:\!\! N_i ]$ refers
to the standard Python notation for extracting all but the last element of $u_{\Gamma_i}$.  Since
\begin{align}
u_{\Gamma_1}[N_1] &= u_{\Gamma_2}[0] = u_{N_10} \nonumber \\
u_{\Gamma_2}[N_2] &= u_{\Gamma_3}[0] = u_{N_1N_2}\label{eq:uboundaryconstraint} \\
u_{\Gamma_3}[N_1] &= u_{\Gamma_4}[0] = u_{0N_2} \nonumber \\
u_{\Gamma_4}[N_2] &= u_{\Gamma_1}[0] = u_{00} \nonumber,
\end{align}
it follows that $ u^{(e)}_\Gamma$ is indeed an embedding of $u_\Gamma$, and conversely, 
$u_\Gamma$ is the restriction of $u_{\Gamma}^{(e)}$.  Hence,
\begin{equation}
\label{eq:uboundaryembed}
\nonumber
u^{(e)}_\Gamma = Eu_\Gamma
\end{equation}
in complete analogy with Remark~\ref{rem:innerprodembed}.  The matrix representation of 
the embedding operator $E\in\mathbb{R}^{2(N+2)\times 2N}$ is given by
\begin{equation}
\label{eq:eembedboundary}
E = 
\begin{blockarray}{ccc|cc|cc|cc}
  &{\scriptstyle 0} &c_1 & &c_2 & &c_3 & &c_4 \\
  \begin{block}{c(cc|cc|cc|cc)}
  {\scriptstyle 0}     
      &1 &                 &  &                 &  &                 &  &                \\
      &  &\tilde{I}^{(1)}  &  &                 &  &                 &  &                \\
  r_1 &  &0                &1 &                 &  &                 &  &                \\
  \cline{1-9}
      &  &                 &1 &                 &  &                 &  &                \\
      &  &                 &  &\tilde{I}^{(2)}  &  &                 &  &                \\
  r_2 &  &                 &  &0                &1 &                 &  &                \\
  \cline{1-9}
      &  &                 &  &                 &1 &                 &  &                \\
      &  &                 &  &                 &  &\tilde{I}^{(1)}  &  &                \\
  r_3 &  &                 &  &                 &  &0                &1 &                \\
  \cline{1-9}
      &  &                 &  &                 &  &                 &1 &                \\
      &  &                 &  &                 &  &                 &  &\tilde{I}^{(2)} \\
  r_4 &1 &                 &  &                 &  &                 &  &0               \\
  \end{block}
\end{blockarray}\,,
\end{equation}
where $\tilde{I}^{(i)} \in \mathbb{R}^{(N^{(i)}-1)\times (N^{(i)}-1)}$.  The indices 
$r_j$ and $c_j$ correspond to the last element of each block in $u^{(e)}_\Gamma$ and $u_\Gamma$:
\begin{align*}
r_1 &= N_1 \\
r_2 &= N_1+N_2+1 \\
r_3 &= 2N_1+N_2+2 \\
r_4 &= 2N_1+2N_2+3
\end{align*}
and
\begin{align*}
c_1 &= N_1-1 \\
c_2 &= N_1+N_2-1 \\
c_3 &= 2N_1+N_2-1 \\
c_4 &= 2N_1+2N_2-1.
\end{align*}

For arbitrary $2(N+2)$-dimensional grid vectors $u_{\Gamma_+}$, $v_{\Gamma_+}$, 
not necessarily subject to the constraints (\ref{eq:uboundaryconstraint}), we define the scalar 
product
\begin{equation}
\label{eq:bstateinnerprodplus}
\langle u, v \rangle_+  \equiv u^T_{\Gamma_+}H^{(+)}_\Gamma v^T_{\Gamma_+},
\end{equation}
where
\begin{equation}
\label{eq:hgamma+}
H^{(+)}_\Gamma  =
\begin{pmatrix}
H_1 \\
&H_2 \\
&&J_1H_1J_1 \\
&&&J_2H_2J_2
\end{pmatrix} \in \mathbb{R}^{2(N+2)\times 2(N+2)}
\end{equation}
with $H_i$ being one-dimensional norms and $J_i$ anti-diagonal permutation matrices~\eqref{eq:adiag}.
Note that $H^{(+)}_\Gamma $ has the same structure as $H^{(+)}$ in (\ref{eq:hplus}).  
Following the same pattern as before, we define a scalar product on 
$\mathbb{R}^{2N}\times  \mathbb{R}^{2N}$:
\begin{equation}
\label{eq:bstateinnerprod}
\langle u_{\Gamma}, v_{\Gamma} \rangle  \equiv \left(u^{(e)}_\Gamma\right)^TH^{(+)}_\Gamma u^{(e)}_\Gamma = 
\left(Eu_\Gamma\right)^TH^{(+)}_\Gamma \left(Eu_\Gamma\right).
\end{equation}
It follows immediately from the definition 
that $\langle \cdot, \cdot \rangle$ is a well-defined scalar product, 
cf.~(\ref {eq:innerprodembed}).  In particular, for grid vectors given by 
(\ref{eq:segmentbstate}) and (\ref{eq:2dbstate}):
\begin{align*}
\langle u_{\Gamma}, v_{\Gamma} \rangle 
  &\equiv \langle u_{\Gamma_1}, v_{\Gamma_1}\rangle_1 + \langle u_{\Gamma_2}, v_{\Gamma_2}\rangle_2 
   + \langle u_{\Gamma_3}, v_{\Gamma_3}\rangle_3 + \langle u_{\Gamma_4}, v_{\Gamma_4}\rangle_4 \\
  &= (u_{0}, v_{0})_1 + (u^{N_1}, v^{N_1})_2
   + (u_{N_2}, v_{N_2})_1 + (u^{0}, v^{0})_2
\end{align*}
and thus
\[
\langle 1, 1 \rangle = 4,
\]
since $(1, 1)_i=1$, see~(\ref{eq:constraint}).

\begin{define}
\label{def:b2dinnerprodembed}
Let the inner product space $V_\Gamma$ be a real vector space with the inner product
\[
\langle\cdot,\cdot\rangle:V_\Gamma\times V_\Gamma \rightarrow \mathbb{R}
\]
for all vectors $u_\Gamma,v_\Gamma\in V_\Gamma$ given by (\ref{eq:2dbstate}) and where
$\langle\cdot,\cdot\rangle$ is defined in (\ref{eq:bstateinnerprod}). \hfill$\Box$
\end{define}

\begin{rem}
\label{rem:2dpartsum1}
Given Definition~\ref{def:b2dinnerprodembed}, Theorem~\ref{thm:2dpartsum} can be expressed as
\begin{align*}
(u,D_xv) &= \langle u,v\rangle_2 - \langle u,v\rangle_4 - (D_xu,v) \\
(u,D_yv) &= \langle u,v\rangle_3 - \langle u,v\rangle_1 - (D_yu,v),
\end{align*}
where we have dropped the subscripts of $u_{\Gamma_i}, v_{\Gamma_i}$ in the inner products
$\langle \cdot,\cdot\rangle_i$ for each boundary segment $\Gamma_i$.  The restriction
of $u$ and $v$ to a particular boundary $\Gamma_i$ is implied by the notation. \hfill$\Box$
\end{rem}

\subsection{Energy estimates}
\label{sec:energy}
Consider a symmetric hyperbolic system in two dimensions ($A,B\in\mathbb{R}^{d\times d}$):
\begin{align}
u_t + Au_x + Bu_y &= 0, \quad (x,y) \in \Omega, \quad t>0 \label{eq:hyper2d}\\
          u(x,y,0)&= f(x,y) \nonumber
\end{align}
with characteristic boundary conditions:
\begin{equation}
\label{eq:2dcharbc}
L_iu \equiv \left(\left[Q^{(i)}_I\right]^T - S_i\left[Q^{(i)}_{II}\right]^T\right)u=0, 
\quad (x,y) \in \Gamma_i, \quad t>0.
\end{equation}
The $d$ by $d$ matrices $Q^{(1)}=Q^{(3)}$ diagonalize $B^T=B$.  Similarly,
$Q^{(2)}=Q^{(4)}$ diagonalize $A^T=A$:
\begin{equation}
\label{eq:2dbdiag}
\lambda^{(1,3)} \equiv \left[Q^{^{(1,3)}}\right]^TBQ^{(1,3)} , \quad 
\lambda^{(2,4)} \equiv \left[Q^{^{(2,4)}}\right]^TAQ^{(2,4)} .
\end{equation}
There are $d^{(i)}_1$ (locally) ingoing characteristics represented by the eigenvectors $Q^{(i)}_I$ 
for each boundary segment $\Gamma_i$.  Similarly, $Q^{(i)}_{II}$ corresponds to the outgoing ones.

The semidiscrete approximation of (\ref{eq:hyper2d}) is written as
\begin{align}
v_t + P(AD_x + BD_y)Pv &= 0, \quad t>0 \label{eq:semihyper2d} \\
          v(0)&= f, \quad f=Pf ,\nonumber
\end{align}
where $A,B:V \rightarrow V$ are defined by
\[
A \equiv \mbox{diag}(A_{ij}), B \equiv \mbox{diag}(B_{ij})
\in \mathbb{R}^{(N_1+1)(N_2+1)d\times (N_1+1)(N_2+1)d}
\]
with $A_j=A\in\mathbb{R}^{d\times d}$, $B_j=B\in\mathbb{R}^{d\times d}$.  We 
have deliberately used the same symbols for the coefficient matrices in the 
analytic and semidiscrete formulations.  

As usual,
\begin{equation}
\label{eq:2dcharbcproj}
\nonumber
P = I - L^+L.
\end{equation}
The boundary operator $L:V \rightarrow V_{\Gamma_+}$ is represented as a block
matrix
\begin{equation}
\label{eq:2dcharbcdiscrete}
L \equiv 
\begin{pmatrix}
L_1 \\
L_2 \\
L_3 \\
L_4
\end{pmatrix},
\end{equation}
where
\begin{align}
L_1 &= I_2[0,\!:] \otimes I_1 \otimes \left(\left[Q^{(1)}_I\right]^T - S_1\left[Q^{(1)}_{II}\right]^T \right) \nonumber \\
L_2 &= I_2 \otimes I_1[N_1,\!:] \otimes \left(\left[Q^{(2)}_I\right]^T - S_2\left[Q^{(2)}_{II}\right]^T \right) \label{eq:2dcharbcdiscrete1} \\
L_3 &= I_2[N_2,\!:] \otimes I_1 \otimes \left(\left[Q^{(3)}_I\right]^T - S_3\left[Q^{(3)}_{II}\right]^T \right) \nonumber \\
L_4 &= I_2 \otimes I_1[0,\!:] \otimes \left(\left[Q^{(4)}_I\right]^T - S_4\left[Q^{(4)}_{II}\right]^T \right) \nonumber
\end{align}
and
\[
Q^{(i)}_{I}\in\mathbb{R}^{d \times d^{(i)}_1},\quad
Q^{(i)}_{II}\in\mathbb{R}^{d \times (d-d^{(i)}_1)},\quad
S_i\in\mathbb{R}^{d^{(i)}_1\times (d-d^{(i)}_1)}.
\]

\begin{prop}
\label{prop:hyper2dsemistability}
The semidiscrete system (\ref{eq:semihyper2d}) - (\ref{eq:2dcharbcdiscrete1}) is
a strictly stable approximation of (\ref{eq:hyper2d}) - (\ref{eq:2dcharbc}).
\end{prop}

\noindent
{\bf Proof}:
Applying the energy method to (\ref{eq:semihyper2d}) will result in one-dimensional boundary terms like
\[
(u^{0},Au^{0})_2, \quad (u^{N_1},Au^{N_1})_2, \quad
(u_0,Bu_0)_1, \quad (u_{N_2},Bu_{N_2})_1,
\]
where the one-dimensional block diagonal matrices $A$ and $B$ are defined as
\begin{align*}
A &\equiv \mbox{diag}(A_j) \in \mathbb{R}^{(N_2+1)d\times (N_2+1)d} \\
B &\equiv \mbox{diag}(B_j) \in \mathbb{R}^{(N_1+1)d\times (N_1+1)d}
\end{align*}
with $A_j=A\in \mathbb{R}^{d\times d}, B_j=B\in \mathbb{R}^{d\times d}$.  Similarly, let
\[
Q_i \equiv \mbox{diag}\left(Q^{(i)}\right) \in \mathbb{R}^{(N_i+1)d\times (N_i+1)d}, 
\quad i=1,2,
\]
where the eigenvectors $Q^{(i)}\in\mathbb{R}^{d\times d}$ are defined by 
(\ref{eq:2dbdiag}).  Thus,  $Q_1$ and $Q_2$ correspond to the 
eigenvectors of $B\in\mathbb{R}^{(N_1+1)d\times (N_1+1)d}$ and $A\in\mathbb{R}^{(N_2+1)d\times (N_2+1)d}$.  
The block diagonal structure of $Q_i$ implies that
\begin{equation}
\label{eq:eigencommute}
\nonumber
H_iQ_i=Q_iH_i,  \quad i=1,2.
\end{equation}
Furthermore, 
\[
Q_3=Q_1, \quad Q_4=Q_2,
\]
since $Q^{(3)} = Q^{(1)}$ and $Q^{(4)} = Q^{(2)}$.  All four boundary terms can thus be 
reduced to
\begin{equation}
\label{eq:standardform}
\nonumber
(v,\Lambda_i v)_i, \qquad
\Lambda_i \equiv \mbox{diag}(\lambda^{(i)})\in \mathbb{R}^{(N_i+1)d\times (N_i+1)d};
\end{equation}
$\lambda^{(i)}$ is given by (\ref{eq:2dbdiag}) and  
$v$ represents the characteristic boundary state corresponding to 
the state vector $u_{\Gamma_i}$, $i=1,\ldots, 4$.  

Since $u=Pu$, where $u$ is the two-dimensional state vector (\ref{eq:2dstate}), 
(\ref{eq:2daltstate}) and where $P$ is the projection representing (\ref{eq:2dcharbcdiscrete1}), 
it follows that each component $v_j$ satisfies
\begin{equation}
\label{eq:charbc2d}
L_iv_j=0, \qquad
L_i=
\overset{\begin{matrix}{\scriptstyle d^{(i)}_1} &{\scriptstyle d^{(i)}_2} &{\scriptstyle d^{(i)}_3}\end{matrix}}{
\begin{pmatrix}
I &\ 0 &-S_i
\end{pmatrix}}\begin{matrix}{\scriptstyle d^{(i)}_1}\end{matrix}, \quad d^{(i)}_1 + d^{(i)}_2 + d^{(i)}_3 = d.
\end{equation}
Without loss of generality, we can drop the dependence on the boundary segment $\Gamma_i$.  Thus, it suffices to analyze boundary terms where the coefficient matrices $A$ and $B$ are diagonal:
\[
(u,\Lambda u)
\]
and where each boundary point satisfies the analytic boundary condition
\[
Lu_j = 0, \qquad j=0,\ldots, N.
\]

Going forward, we will partition $\lambda\in\mathbb{R}^{d\times d}$ as
\[
\lambda = 
\begin{pmatrix}
\lambda_+ \\
&0 \\
&&\lambda_-
\end{pmatrix}
\begin{matrix}
{\scriptstyle d_1} \\
{\scriptstyle d_2} \\
{\scriptstyle d_3}
\end{matrix},
\]
where $\lambda_+$ and $\lambda_-$ are the strictly positive and negative parts
of $\lambda$.  Define auxiliary $d$ by $d$ diagonal matrices
\[
\lambda^{(+)} \equiv
\begin{pmatrix}
\lambda_+ \\
&0 \\
&&0
\end{pmatrix}, \quad
\lambda^{(-)} \equiv
\begin{pmatrix}
0 \\
&0 \\
&&\lambda_-
\end{pmatrix},
\]
i.~e., $\lambda = \lambda^{(+)} + \lambda^{(-)}$.  Next, we construct
\[
\Lambda_+ \equiv \mbox{diag}\left(\lambda^{(+)}\right), \qquad 
\Lambda_- \equiv \mbox{diag}\left(\lambda^{(-)}\right) \in \mathbb{R}^{(N+1)d\times(N+1)d} .
\]
Hence,
\[
(u,\Lambda u) = (u,\Lambda_+u) + (u,\Lambda_-u) = (u,\Lambda_+u) - (u,|\Lambda_-|u).
\]
We use the notation $|\cdot|$ to indicate the absolute value of the nonzero 
elements of $\Lambda_-$.  This expression can be bounded from above by zero
as the following argument will show.  We first observe that
\[
 H\Lambda_+ = \left[\Lambda_+\right]^{1/2}H\left[\Lambda_+\right]^{1/2},
\]
since the $d$ by $d$ diagonal blocks of $\Lambda_+$ all equal to $\lambda_+$ and
since the blocks of $H$ are given by $h_{ij}I$, where $I$ is the $d$ by $d$ identity 
matrix.  The matrix square root is meaningful since 
$\Lambda_+$ is diagonal (not just block diagonal) with nonnegative elements.  Define
\[
v \equiv \left[\Lambda_+\right]^{1/2}u.
\]
Thus,
\[
(u,\Lambda_+u)_H = (v,v)_H.
\]
Similarly,
\[
(u,|\Lambda_-|u)_H = (w,w)_H,
\]
where
\[
w \equiv \left[|\Lambda_-|\right]^{1/2}u.
\]
But all $H$-norms  are equivalent, i.~e., we can find constants
$0 < c_- \leq c_+$ such that
\[
c_-(u,u)_2 \leq (u,u) \leq c_+(u,u)_2 
\]
holds for all $u$ where
\[
(u,u)_2 \equiv h\frac{1}{2}u^T_0u_0 + h\sum_{j=1}^{N-1}u^T_ju_j + h\frac{1}{2}u^T_Nu_N.
\]
Hence,
\begin{align*}
(v,v) &\leq c_+(v,v)_2 \\
(w,w) &\geq c_-(w,w)_2 .
\end{align*}
The following inequality has thus been established:
\begin{equation}
\label{eq:eestimate}
(u,\Lambda u) \leq c_+(u,\Lambda_+u)_2 - c_-(u,|\Lambda_-|u)_2.
\end{equation}
The first term of the right-hand side is made up of terms like
(ignoring the mesh size $h$)
\[
c_+u^T_j\lambda^{(+)}u_j.
\]
Dropping the spatial script $j$ we notice that
\[
c_+u^T\lambda^{(+)}u = c_+u^T_I\lambda_+u_I.
\]
But
\[
Lu = 0 \quad \Longleftrightarrow \quad u_I = Su_{II}.
\]
Hence,
\[
c_+u^T\lambda^{(+)}u = c_+u^T_{II}S^T\lambda_+Su_{II}.
\]
This expression can be balanced by the corresponding term from second
scalar product of (\ref{eq:eestimate}):
\[
c_-u^T|\lambda^{(-)}|u = c_-u^T_{II}|\lambda_-|u_{II},
\]
whence,
\[
c_+u^T\lambda^{(+)}u - c_-u^T|\lambda^{(-)}|u = u^T_{II}\left[c_+S^T\lambda_+S - c_-|\lambda_-| \right]u_{II}
\leq 0
\]
if $S$ is sufficiently small. Since this estimate holds for each boundary state $u_j$, an energy
estimate follows:
\[
(u,\Lambda u) \leq 0,
\]
which proves the claim \hfill$\Box$

\subsection{Structure of the pseudoinverse}
\label{sec:2dcharbc}
In section~\ref{sec:bc} we derived conditions that rendered a particularly simple
expression for the pseudoinverse $L^+$.  We will now
show that this result can be applied to characteristic boundary conditions 
(\ref{eq:2dcharbc}) in two space dimensions provided $H$ is a restricted full norm.  
This will be done by proving that (\ref{eq:lhpseudocommute}) holds.  

Combining (\ref{eq:2dcharbc}) and (\ref{eq:charbc2d}):
\begin{equation}
\label{eq:2dcharbcmatrix}
\nonumber
L_i=\begin{matrix}{\scriptstyle d^{(i)}_1}\end{matrix}
\overset{\begin{matrix}{\scriptstyle d^{(i)}_1} &{\scriptstyle d^{(i)}_2} &{\scriptstyle d^{(i)}_3}\end{matrix}}{
\begin{pmatrix}
I &\ 0 &-S_i
\end{pmatrix}}\left[Q^{(i)}\right]^T, \quad d^{(i)}_1 + d^{(i)}_2 + d^{(i)}_3 = d, \quad i=1,\ldots,4.
\end{equation}
As opposed to the earlier energy analysis of the boundary conditions, the subsequent
analysis does not require the detailed structure of the above expression.  To
simplify the notation, we pad the leading matrix with zero blocks:
\[
\overset{\begin{matrix}{\scriptstyle d^{(i)}_1} &{\scriptstyle d^{(i)}_2} &{\scriptstyle d^{(i)}_3}\end{matrix}}{
\begin{pmatrix}
I &\ 0 &-S_i \\
0 &0 &0 \\
0 &0 &0
\end{pmatrix}}\begin{matrix}
{\scriptstyle d^{(i)}_1} \\
{\scriptstyle d^{(i)}_2} \\
{\scriptstyle d^{(i)}_3}
\end{matrix}.
\]
Note that $L^+$ is well defined no matter if $L$ has full rank or not.  From now on $L_i$
will be regarded as a $d$ by $d$ matrix for each boundary point $(x,y)\in\Gamma_i$.  

Define the auxiliary boundary operators
\begin{equation}
\label{eq:2dcharbcaux1}
\nonumber
\tilde{L}_i \equiv
\begin{pmatrix}
L_i \\
&\ddots \\
&&L_i
\end{pmatrix} \in \mathbb{R}^{(N_1+1)d \times (N_1+1)d}, 
\quad i=1,3,
\end{equation}
and
\begin{align*}
\tilde{L}_2 &=
\begin{pmatrix}
0 &\ldots &0 &L_2
\end{pmatrix} \in \mathbb{R}^{d \times (N_1+1)d} \\
\tilde{L}_4 &=
\begin{pmatrix}
L_4 &0 &\ldots &0
\end{pmatrix} \in \mathbb{R}^{d \times (N_1+1)d}.
\end{align*}
The boundary conditions (\ref{eq:2dcharbcdiscrete1}) can thus be expressed as
\begin{align}
\tilde{L}_1u_0 &= 0 \nonumber \\
\tilde{L}_2u_j &= 0, \quad j=0,\ldots,N_2 \label{eq:bc1} \\
\tilde{L}_3u_{N_2} &= 0 \nonumber \\
\tilde{L}_4u_j &= 0, \quad j=0,\ldots,N_2. \nonumber
\end{align}
Hence,
\begin{align*}
L_1 &= 
\begin{pmatrix}
\tilde{L}_1 &0 &\ldots &0
\end{pmatrix} \in \mathbb{R}^{(N_1+1)d \times (N_1+1)(N_2+1)d} \\
L_3 &= 
\begin{pmatrix}
0 &\ldots &0 &\tilde{L}_3
\end{pmatrix} \in \mathbb{R}^{(N_1+1)d \times (N_1+1)(N_2+1)d} \\
L_i &= 
\begin{pmatrix}
\tilde{L}_i \\
&\ddots \\
&&\tilde{L}_i
\end{pmatrix} \in \mathbb{R}^{(N_2+1)d \times (N_1+1)(N_2+1)d}, 
\quad i=2,4,
\end{align*}
where $L_i$ constitute the blocks of the boundary operator $L:V\rightarrow V_{\Gamma_+}$ (\ref{eq:2dcharbcdiscrete}).
We will show that
\[
L_iH = \bar{H}_iL_i, \quad i=1,\ldots 4,
\]
whence
\[
LH = \bar{H}L, \quad
\bar{H} = 
\begin{pmatrix}
\bar{H}_1 \\
&\bar{H}_2 \\
&&\bar{H}_3 \\
&&&\bar{H}_4
\end{pmatrix}
\begin{matrix}
{\scriptstyle (N_1+1)d} \\
{\scriptstyle (N_2+1)d} \\
{\scriptstyle (N_1+1)d} \\
{\scriptstyle (N_2+1)d}
\end{matrix}.
\]

We begin with $L_1H$ and $L_3H$.  By (\ref{eq:h2d}) and (\ref{eq:hxhy}):
\[
H=H_2\otimes H_1 = H_xH_y = H_yH_x.
\]
Hence,
\begin{align*}
L_1H_x &= 
\begin{pmatrix}
\tilde{L}_1H_1 &0 &\ldots &0
\end{pmatrix} = 
H_1\begin{pmatrix}
\tilde{L}_1 &0 &\ldots &0
\end{pmatrix} = H_1L_1 \\
L_3H_x &= 
\begin{pmatrix}
0 &\ldots &0 &\tilde{L}_3H_1
\end{pmatrix} = 
H_1\begin{pmatrix}
0 &\ldots &0 &\tilde{L}_3
\end{pmatrix} = H_1L_3,
\end{align*}
and
\begin{align*}
L_1H_y &=
\begin{pmatrix}
h^{(2)}_{00}\tilde{L}_1 &\ldots &h^{(2)}_{0N_2}\tilde{L}_1
\end{pmatrix} \\
L_3H_y &=
\begin{pmatrix}
h^{(2)}_{N_20}\tilde{L}_3 &\ldots &h^{(2)}_{N_2N_2}\tilde{L}_3
\end{pmatrix}.
\end{align*}
If we require that
\[
h^{(2)}_{0j} = 0, \quad 0<j\leq N_2 \qquad \mbox{and} \qquad 
h^{(2)}_{N_2j} = 0, \quad 0\leq j< N_2,
\]
then
\begin{align*}
L_1H_y &= h^{(2)}_{00}L_1 \\
L_3H_y &= h^{(2)}_{N_2N_2}L_3
\end{align*}
and thus
\begin{align*}
L_1H &= \bar{H}_1L_1, \qquad \bar{H}_1 \equiv h^{(2)}_{00}H_1 \\
L_3H &= \bar{H}_3L_3, \qquad \bar{H}_3 \equiv h^{(2)}_{N_2N_2}H_1.
\end{align*}

Next, we turn to $L_iH$, $i=2,4$:
\[
L_iH_x = 
\begin{pmatrix}
\tilde{L}_iH_1 \\
&\ddots \\
&&\tilde{L}_iH_1
\end{pmatrix},
\quad
\]
where
\begin{align*}
\tilde{L}_2H_1 &= 
\begin{pmatrix}
h^{(1)}_{N_10}L_2 &\ldots &h^{(1)}_{N_1N_1}L_2
\end{pmatrix} \\
\tilde{L}_4H_1 &= 
\begin{pmatrix}
h^{(1)}_{00}L_4 &\ldots &h^{(1)}_{0N_1}L_4
\end{pmatrix}.
\end{align*}
By requiring
\[
h^{(1)}_{0j} = 0, \quad 0<j\leq N_1 \qquad \mbox{and} \qquad 
h^{(1)}_{N_1j} = 0, \quad 0\leq j< N_1,
\]
we arrive at
\begin{align*}
\tilde{L}_2H_1 &= h^{(1)}_{N_1N_1}\tilde{L}_2 \\
\tilde{L}_4H_1 &= h^{(1)}_{00}\tilde{L}_4,
\end{align*}
i.~e.,
\begin{align*}
L_2H_x &= h^{(1)}_{N_1N_1}L_2 \\
L_4H_x &= h^{(1)}_{00}L_4.
\end{align*}
Furthermore,
\[
L_iH_y = 
\begin{pmatrix}
h^{(2)}_{00}\tilde{L}_i &\ldots &h^{(2)}_{0N_2}\tilde{L}_i \\
\vdots &&\vdots \\
h^{(2)}_{N_20}\tilde{L}_i &\ldots &h^{(2)}_{N_2N_2}\tilde{L}_i
\end{pmatrix} = H_2L_i, \quad i=2,4,
\]
which imposes no new constraints on $H_1$ and $H_2$.  Hence,
\begin{align*}
L_2H &= \bar{H}_2L_2, \qquad \bar{H}_2 \equiv h^{(1)}_{N_1N_1}H_2 \\
L_4H &= \bar{H}_4L_4, \qquad \bar{H}_4 \equiv h^{(1)}_{00}H_2.
\end{align*}
Summing up, if we take $H_1$ and $H_2$ to be restricted full norms, then 
\begin{equation}
\label{eq:pseudocommute}
LH = \bar{H}L,
\end{equation}
where $L$ represents the characteristic boundary conditions (\ref{eq:2dcharbc}),
$H$ is given by (\ref{eq:h2d}), and where
\begin{equation}
\label{eq:charbcnorm}
\bar{H} = 
\begin{pmatrix}
\bar{H}_1 \\
&\bar{H}_2 \\
&&\bar{H}_3 \\
&&&\bar{H}_4
\end{pmatrix}
\begin{matrix}
{\scriptstyle (N_1+1)d} \\
{\scriptstyle (N_2+1)d} \\
{\scriptstyle (N_1+1)d} \\
{\scriptstyle (N_2+1)d}
\end{matrix}.
\end{equation}

\begin{rem}
\label{rem:agnostic}
We have tacitly assumed that the boundary state is represented by the four state vectors
\[
u_0, \ u_{N_2} \in\mathbb{R}^{(N_1+1)d}, \qquad u^{0}, \ u^{N_1} \in\mathbb{R}^{(N_2+1)d}.
\]
Using $u_{\Gamma_i}$ instead simply amounts to reordering (\ref{eq:bc1}):
\[
J_1L_3u=0, \quad J_2L_4u=0.
\]
But $J^2_i=I$.  Thus
\begin{align*}
J_1L_3H &= J_1\bar{H}_3L_3 = \left[J_1\bar{H}_3J_1\right]J_1L_3\\
J_2L_4H &= J_2\bar{H}_4L_4 = \left[J_2\bar{H}_4J_2\right]J_2L_4,
\end{align*}
which implies (\ref{eq:pseudocommute}) but with $J_1\bar{H}_3J_1$,  
$J_2\bar{H}_4J_2$ substituted for $\bar{H}_3$ and $\bar{H}_4$.  Hence, 
we can safely ignore the permutation matrices $J_1$ and $J_2$ when
establishing $LH=\bar{H}L$.  Finally, it should be noted that if 
$H_1$ and $H_2$ satisfy (\ref{eq:norm}), then
\[
J_1\bar{H}_3J_1 = \bar{H}_3, \quad 
J_2\bar{H}_4J_2 = \bar{H}_4.
\]
\hfill$\Box$
\end{rem}

\subsubsection{The simplified projection $P$ revisited}
\label{sec:simpleprojrevisited}
We showed in the previous section that
\[
LH = \bar{H}L,
\]
if $L$ represents the characteristic boundary conditions (\ref{eq:2dcharbcdiscrete}), (\ref{eq:2dcharbcdiscrete1}) of the two-dimensional hyperbolic system (\ref{eq:hyper2d}); the scalar product in $V$ is represented by
\[
H = H_2\otimes H_1,
\]
where $H_1$ and $H_2$ are restricted full norms (\ref{eq:rfnorm}); $\bar{H}$ is defined by (\ref{eq:charbcnorm}):
\begin{equation}
\label{eq:charbcnorm1}
\bar{H} = 
\begin{pmatrix}
h^{(2)}_{00}H_1 \\
&h^{(1)}_{N_1N_1}H_2 \\
&&h^{(2)}_{N_2N_2}H_1 \\
&&&h^{(1)}_{00}H_2
\end{pmatrix}.
\end{equation}
Analogous to Section~\ref{sec:psimplified}, we then regard $L$ as a mapping $L:V\rightarrow V_{\Gamma_+}$ with $H^{(+)}_\Gamma $ defined by (\ref{eq:charbcnorm1}) instead of (\ref{eq:hgamma+}).  Hence,
\[
L^* = L^T
\]
and thus
\[
P = I - L^T\left[LL^T\right]^+L
\]
in complete agreement with (\ref{eq:psimplified2}).

\begin{rem}
\label{rem:simpleprojrevisited}
In sections~\ref{sec:energy} and~\ref{sec:simpleprojrevisited} it was assumed that $L$ is a mapping 
$L:V\rightarrow V_{\Gamma_+}$ and not $L:V\rightarrow V_\Gamma$.  The underlying Euclidean spaces of 
$V_{\Gamma_+}$ and $V_\Gamma$ are $\mathbb{R}^{2(N+2)d}$ and $\mathbb{R}^{2Nd}$.  It is largely a matter of convenience which boundary state space to choose.  Both lead to projection operators that satisfy 
\[
Pv=v, \quad HP = P^TH,
\]
which is what the stability analysis requires.  As a general rule, algebraic manipulation become easier if one chooses 
$L:V\rightarrow V_{\Gamma_+}$.  In our case, all state vectors in $V_{\Gamma_+}$ satisfy the constraints 
(\ref{eq:uboundaryconstraint}) by construction.  The case $L:V\rightarrow V_\Gamma$ corresponds to redefining 
(\ref{eq:2dcharbcdiscrete1}) to account for the corner points only once. \hfill$\Box$
\end{rem}

\section{Two space dimensions, multiblock case}
\label{sec:2dimmb}
The embedding operator $E:V\rightarrow V_+$ introduced in Section~\ref{sec:embedding} played a crucial role
when establishing difference operators $D:V\rightarrow V$ that satisfy summation by parts in a multiblock 
scenario.  This was done by embedding a lower-dimensional vector space as a manifold in a higher-dimensional 
space.  This technique will now be generalized to domains $\Omega$ in two space dimensions.  We will
restrict ourselves to the case where grid lines match at the subdomain interfaces.

\subsection{Two-block difference operators, case 1}
\label{sec:twoblock1diffop}
Theorem~\ref{thm:2dpartsum} will be generalized to two-dimensional difference operators defined on 
$\Omega=\Omega_1\cup\Omega_2$, $\Omega_1=[0,1/2]\times[0,1]$, $\Omega_2=[1/2,1]\times[0,1]$. 
The grid points are defined as
\begin{align*}
\Omega_1: \quad (x_i, y_j) &= (ih^{(1)}_1, jh_2), \quad 0\leq i\leq N^{(1)}_1, 0\leq j \leq N_2 \\
\Omega_2: \quad (x_i, y_j) &= (0.5 + ih^{(2)}_1, jh_2), \quad 0\leq i\leq N^{(2)}_1, 0\leq j \leq N_2,
\end{align*}
where the mesh sizes are defined as
\[
h^{(i)}_1 = \frac{1}{2N^{(i)}_1}, \ i=1,2, \quad
h^{(1)}_2 = h^{(2)}_2 = h_2 = \frac{1}{N_2}.
\]
It should be observed that the grid spacings $h^{(j)}_2, j=1,2$, in the $y$-dimension are assumed to 
be the same across the interface between $\Omega_1$ and $\Omega_2$.  On each domain we 
define grid vectors (\ref{eq:2dstate}), scalar products (\ref{eq:h2d}), 
(\ref{eq:2dinnerprod}) and difference operators (\ref{eq:dxdy}):
\begin{equation}
\label{eq:opdefi}
\Omega_i: \quad u^{(i)},v^{(i)}, \quad H^{(i)} = H^{(i)}_xH^{(i)}_y, \quad D^{(i)}_x, D^{(i)}_y.
\end{equation}
\begin{figure}
\centering
\includegraphics[width=7cm]{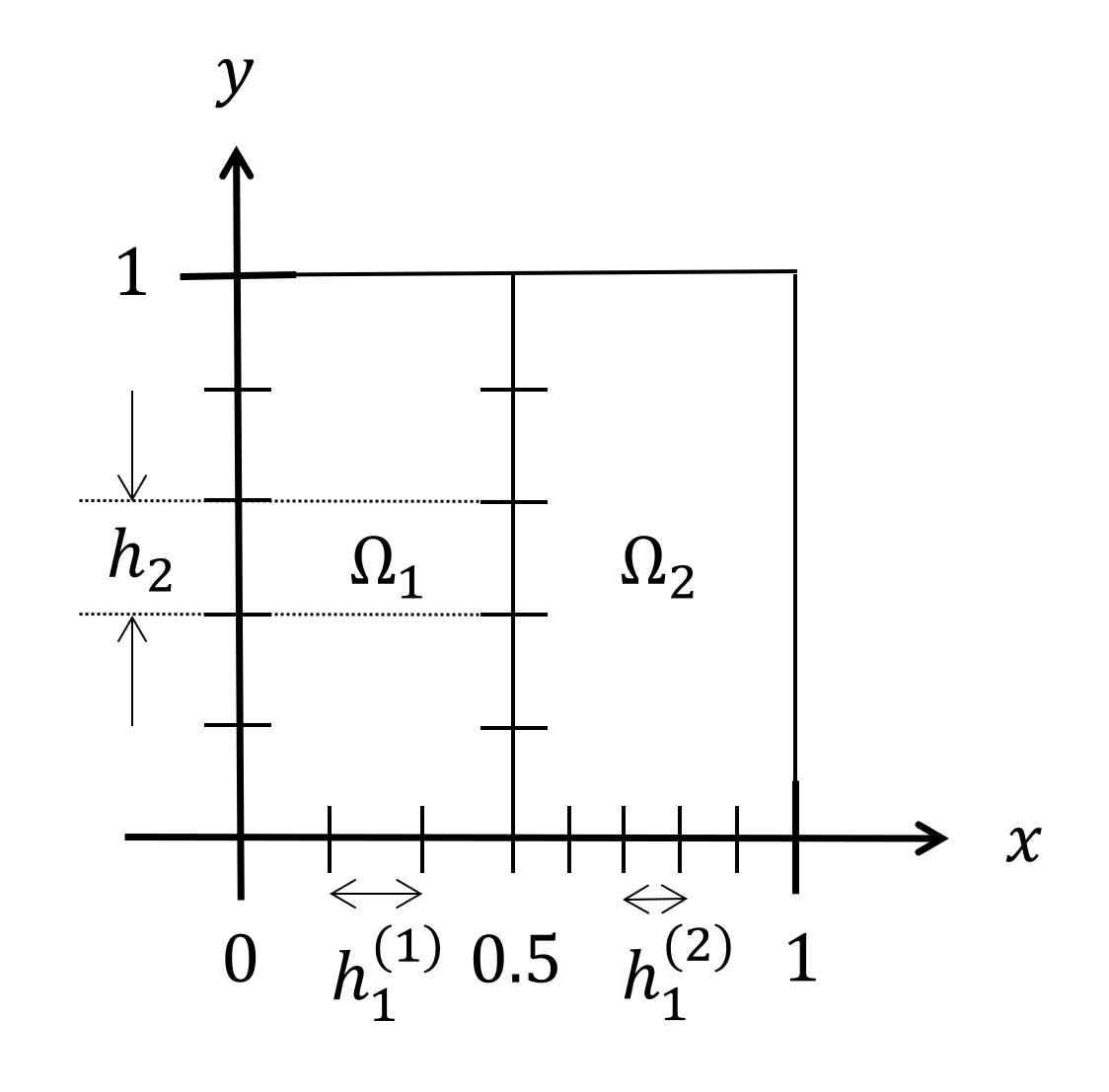}
\caption{Two blocks, case 1}
\label{fig:twoblockscase1}
\end{figure}

Analogous to the one-dimensional case, cf. Remark~\ref{rem:auginnerprod}, we define grid
vectors on $\Omega_+=\Omega_1+\Omega_2$:
\begin{equation}
\label{eq:2duvplus}
\nonumber
u^{(+)} \equiv
\begin{pmatrix}
u^{(1)} \\
u^{(2)}
\end{pmatrix}
\begin{matrix}
{\scriptstyle r_1} \\
{\scriptstyle r_2}
\end{matrix}, \quad
v^{(+)} \equiv
\begin{pmatrix}
v^{(1)} \\
v^{(2)}
\end{pmatrix}
\begin{matrix}
{\scriptstyle r_1} \\
{\scriptstyle r_2}
\end{matrix}.
\end{equation}
The augmented state space $V_+$ is defined exactly as in 
Definition~\ref{def:auginnerprod}, whence: 
\begin{equation}
\label{eq:2dhplus}
(u^{(+)},v^{(+)})_+ = \left[u^{(+)}\right]^TH^{(+)}v^{(+)},
 \quad
H^{(+)} =
\overset{\raise3pt\hbox{$\scriptstyle c_1 \ \ \ \ \ \ c_2$}}{
\begin{pmatrix}
H^{(1)} \\
&H^{(2)}
\end{pmatrix}}
\begin{matrix}
{\scriptstyle r_1} \\
{\scriptstyle r_2}
\end{matrix},
\end{equation}
where
\[
c_i = r_i = (N^{(i)}_1+1)(N_2+1), \ i=1,2.
\]

\subsubsection{The embedding operator $E_x$}
\label{sec:exembedding}
Let $u, v$ be grid vectors on $\Omega=\Omega_1\cup\Omega_2$ as defined in (\ref{eq:2dstate}):
\[
u, v \in \mathbb{R}^{(N_1+1)(N_2+1)}, \quad N_1 \equiv N^{(1)}_1 + N^{(2)}_1.
\]
Note that we traverse all of $\Omega$ horizontally and then vertically, that is
\[
u_j, v_j \in \mathbb{R}^{N_1+1}, \quad j=0, \ldots, N_2.
\]

As in Section~\ref{sec:embedding}, define a mapping $E_x:\mathbb{R}^{(N_1+1)(N_2+1)}\rightarrow
\mathbb{R}^{(N_1+2)(N_2+1)}$:
\begin{equation}
\label{eq:expart}
E_x = 
\begin{pmatrix}
E^{(1)}_x \\
E^{(2)}_x
\end{pmatrix}, \quad  E^{(i)}_x \equiv I_2 \otimes E^{(i)}_1, \quad 
I_2\in\mathbb{R}^{(N_2+1)\times(N_2+1)},
\end{equation}
and where $E^{(1)}_1$ and $E^{(2)}_1$ are defined by (\ref{eq:e1embed}), (\ref{eq:e2embed}) replacing
$N^{(1)}\rightarrow N^{(1)}_1$, $N^{(2)}\rightarrow N^{(2)}_1$ and $N \rightarrow N_1$, which
reflects the fact that $\Omega_1$ and $\Omega_2$ are joined in the $x$-direction.

Define an inner product on $\mathbb{R}^{(N_1+1)(N_2+1)}\times\mathbb{R}^{(N_1+1)(N_2+1)}$:
\begin{equation}
\label{eq:2dinnerprodembedx}
(u,v)_H \equiv (E_xu)^TH^{(+)}E_xv \quad \Longleftrightarrow \quad
H = E^T_xH^{(+)}E_x.
\end{equation}
Definition~\ref{def:innerprodembed} carries over with (\ref{eq:innerprodembed}) replaced 
by (\ref{eq:2dinnerprodembedx}).
Hence, $E_x$ is a mapping between two inner product spaces $V$ and $V_+$, formally written as
$E_x:V\rightarrow V_+$.  The adjoint $E^*_x$ is given by
\[
E^*_x = H^{-1}E^T_xH^{(+)} \quad \Longrightarrow \quad
E^*_xE_x = I.
\]
The Moore-Penrose conditions thus imply that
\begin{equation}
\label{eq:pseudoembedx}
E^+_x=E^*_x
\end{equation}
just like in the one-dimensional case.

\subsubsection{Multiblock difference operators $D_x$ and $D_y$}
\label{sec:2dmultiblockx}
Similar to (\ref{eq:dplus}), define $D^{(+)}_x,D^{(+)}_y:V_+ \rightarrow V_+$:
\begin{equation}
\label{eq:2ddplus}
D^{(+)}_x =
\overset{\raise3pt\hbox{$\scriptstyle c_1 \ \ \ \ \ \ c_2$}}{
\begin{pmatrix}
D^{(1)}_x \\
&D^{(2)}_x
\end{pmatrix}}
\begin{matrix}
{\scriptstyle r_1} \\
{\scriptstyle r_2}
\end{matrix}, \quad
D^{(+)}_y =
\overset{\raise3pt\hbox{$\scriptstyle c_1 \ \ \ \ \ \ c_2$}}{
\begin{pmatrix}
D^{(1)}_y \\
&D^{(2)}_y
\end{pmatrix}}
\begin{matrix}
{\scriptstyle r_1} \\
{\scriptstyle r_2}
\end{matrix}.
\end{equation}

\begin{define}
\label{def:2ddiffopembedx}
Given the inner product space $V$ with inner product (\ref{eq:2dinnerprodembedx}), the difference
operators $D_x,D_y:V\rightarrow V$ are defined as
\begin{align}
\label{eq:2ddiffopembedx}
D_x &\equiv H^{-1}E^T_xH^{(+)}D^{(+)}_xE_x \\
D_y &\equiv H^{-1}E^T_xH^{(+)}D^{(+)}_yE_x. \nonumber
\end{align}
\ \hfill$\Box$
\end{define}

As in the single domain case, we introduce an alternate column ordering of the (restricted) state vector $u$ on 
$\Omega = \Omega_1\cup\Omega_2$ to allow for a convenient notation when discussing summation by parts in two dimensions:
\begin{equation}
\label{eq:2durestrictedalt}
\nonumber
u =
\begin{pmatrix}
u^0 \\
\vdots \\
u^N
\end{pmatrix}, \quad
u^i \equiv 
\begin{pmatrix}
u_{i0} \\
\vdots \\
u_{iN_2}
\end{pmatrix},\quad  0\leq i \leq N_1.
\end{equation}
The two-dimensional equivalent of Proposition~(\ref{prop:diffopembed}) can be formulated as:

\begin{prop}
\label{prop:2ddiffopembedx}
Let $D_x,D_y:V\rightarrow V$ be as in Definition~\ref{def:2ddiffopembedx}.  Then
\begin{itemize}
\item[(i)] $D_x,D_y$ satify summation by parts with respect to the inner product 
(\ref{eq:2dinnerprodembedx}):
\begin{align*}
(u,D_xv) &= \langle u,v\rangle_2 - \langle u,v\rangle_4 - (D_xu,v) \\
(u,D_yv) &= \langle u,v\rangle_3 - \langle u,v\rangle_1 - (D_yu,v),
\end{align*}
iff $H^{(1)}_2 = H^{(2)}_2 = H_2 \in \mathbb{R}^{(N_2+1)\times (N_2+1)}$, where
the one-dimensional norm $H_2$ is that of (\ref{eq:h2d}); $H_1$ corresponds to 
(\ref{eq:h}).
\item[(ii)] $D_x,D_y$ are consistent approximations of $\partial/\partial x$ and $\partial/\partial y$.
\item[(iii)] $D_x=E^+_xD^{(+)}_xE_x$ and $D_y=E^+_xD^{(+)}_yE_x$.
\end{itemize}
\end{prop}

\noindent
{\bf Proof:}
Define
\begin{equation}
\label{eq:uembed2dx}
u^{(e)} \equiv E_xu =
\begin{pmatrix}
E^{(1)}_xu \\
E^{(2)}_xu
\end{pmatrix}
\equiv
\begin{pmatrix}
u^{(1)} \\
u^{(2)}
\end{pmatrix}.
\end{equation}
This corresponds to row-wise ordering of the elements of $u^{(1)}$ and $u^{(2)}$ (note that $u$ is row-ordered as well):
\[
u^{(k)} =
\begin{pmatrix}
u^{(k)}_0 \\
\vdots \\
u^{(k)}_{N_2} 
\end{pmatrix} \in \mathbb{R}^{(N^{(k)}_1+1)(N_2+1)}, \quad
u^{(k)}_j =
\begin{pmatrix}
u^{(k)}_{0j} \\
\vdots \\
u^{(k)}_{N^{(k)}_1j} 
\end{pmatrix} \in \mathbb{R}^{N^{(k)}_1+1}.
\]
By construction, the following compatibility conditions are fulfilled:
\begin{equation}
\label{eq:2dconstraintx}
\nonumber
u^{(1)}[N^{(1)}_1,:] \equiv 
\begin{pmatrix}
u^{(1)}_{N^{(1)}_10} \\
\vdots \\
u^{(1)}_{N^{(1)}_1N_2}
\end{pmatrix} =
\begin{pmatrix}
u_{N^{(1)}_10} \\
\vdots \\
u_{N^{(1)}_1N_2}
\end{pmatrix} = 
\begin{pmatrix}
u^{(2)}_{00} \\
\vdots \\
u^{(2)}_{0N_2}
\end{pmatrix} \equiv
u^{(2)}[0,:].
\end{equation}
From definitions (\ref{eq:2ddiffopembedx}), (\ref{eq:2dinnerprodembedx}) and (\ref{eq:uembed2dx}):
\[
(u,D_xv) = (u^{(1)},D^{(1)}_xv^{(1)})_{H^{(1)}} + (u^{(2)},D^{(2)}_xv^{(2)})_{H^{(2)}},
\]
where $u^{(i)}, v^{(i)}$, $i=1,2$, satisfy the above constraints.  Apply 
Theorem~\ref{thm:2dpartsum} to $D^{(i)}_x$ defined on $\Omega_i$:
\begin{align*}
(u^{(i)},D^{(i)}_xv^{(i)})_{H^{(i)}} &= (u^{(i)}[N^{(i)}_1,:],v^{(i)}[N^{(i)}_1,:])_{H^{(i)}_2} - (u^{(i)}[0,:],v^{(i)}[0,:])_{H^{(i)}_2} \\
                                     &- (D^{(i)}_xu^{(i)},v^{(i)})_{H^{(i)}}.
\end{align*}
Adding the two equations yields
\begin{align*}
(u,D_xv)   &= (u^{(2)}[N^{(2)}_1,:],v^{(2)}[N^{(2)}_1,:])_{H^{(2)}_2} - (u^{(2)}[0,:],v^{(2)}[0,:])_{H^{(2)}_2} \\
           &+ (u^{(1)}[N^{(1)}_1,:],v^{(1)}[N^{(1)}_1,:])_{H^{(1)}_2} - (u^{(1)}[0,:],v^{(1)}[0,:])_{H^{(1)}_2} \\
		   &- (D_xu,v).
\end{align*}
Thus
\[
(u^{(1)}[N^{(1)}_1,:],v^{(1)}[N^{(i)}_1,:])_{H^{(1)}_2} - (u^{(2)}[0,:],v^{(2)}[0,:])_{H^{(2)}_2} = 0
\]
independently of $u^{(i)}$, $v^{(i)}$, $i=1,2$, iff $H^{(1)}_2 = H^{(2)}_2 = H_2$, 
which proves the first part of the first assertion (we used (\ref{eq:2durestrictedalt}) in the boundary 
integrals).

The second part of the first claim follows more or less directly from
definitions (\ref{eq:2ddiffopembedx}), (\ref{eq:2dinnerprodembedx}) and (\ref{eq:uembed2dx}):
\begin{align*}
(u,D_yv)   &= u^TE^T_xH^{(+)}D^{(+)}_yE_xv \\
           &= (u^{(1)},D^{(1)}_yv^{(1)})_{H^{(1)}} + (u^{(2)},D^{(2)}_yv^{(2)})_{H^{(2)}}.
\end{align*}
According to Theorem~\ref{thm:2dpartsum}
\begin{align*}
(u,D_yv)   &= u^TE^T_xH^{(+)}D^{(+)}_yE_xv \\
           &= (u^{(1)}_{N_2},v^{(1)}_{N_2})_{H^{(1)}_1} + (u^{(2)}_{N_2},v^{(2)}_{N_2})_{H^{(2)}_1} \\
		   &- (u^{(1)}_0,v^{(1)}_0)_{H^{(1)}_1} - (u^{(2)}_0,v^{(2)}_0)_{H^{(2)}_1} \\
		   &- (D^{(1)}_yu^{(1)},v^{(1)})_{H^{(1)}} - (D^{(2)}_yu^{(2)},v^{(2)})_{H^{(2)}}.
\end{align*}
But
\begin{align*}
(u^{(1)}_{N_2},v^{(1)}_{N_2})_{H^{(1)}_1} + (u^{(2)}_{N_2},v^{(2)}_{N_2})_{H^{(2)}_1} &=
\begin{pmatrix}
u^{(1)}_{N_2} \\
u^{(2)}_{N_2}
\end{pmatrix}^T
\begin{pmatrix}
H^{(1)}_1 \\
&H^{(2)}_1
\end{pmatrix}
\begin{pmatrix}
v^{(1)}_{N_2} \\
v^{(2)}_{N_2}
\end{pmatrix} \\
&= (u_{N_2},v_{N_2})_{H_1},
\end{align*}
where the second equality follows from
\[
u^{(1)}_{N_2} = E^{(1)}_1u_{N_2}, \quad u^{(2)}_{N_2} = E^{(2)}_1u_{N_2},
\]
and from (\ref{eq:innerprodembed}).  The remaining terms in the right member of the above expression 
for $(u,D_yv)$ are treated in similar way, which concludes the proof of the first claim.

The second assertion can be proved exactly as in Proposition (\ref{prop:diffopembed}).  The third claim, finally, 
is an immediate consequence of (\ref{eq:pseudoembedx}). \hfill$\Box$

\subsubsection{Structure of $H$, $D_x$, $D_y$ and $E^*_x$}
\label{sec:structurex}
Below, we have gathered some results pertaining to the matrix representation of the
operators $H$, $D_x$, $D_y$ and $E^*_x$. 

\begin{prop}
\label{prop:2dhstructx}
Let $H$ be as in (\ref{eq:2dinnerprodembedx}).  If $H^{(1)}_2 = H^{(2)}_2 \equiv H_2$,
then
\begin{align*}
H   &= H_xH_y = H_yH_x  \\
H_x &= I_2 \otimes H_1  \\
H_y &= H_2\otimes I_1,
\end{align*}
where $H_1$ is the one-dimensional norm defined in (\ref{eq:h}).
\end{prop}

\noindent
{\bf Proof}:
From the definition of $H$:
\[
H = E^T_xH^{(+)}E_x = \left[E^{(1)}_x\right]^TH^{(1)}E^{(1)}_x + \left[E^{(2)}_x\right]^TH^{(2)}E^{(2)}_x.
\]
By Lemma~\ref{lemma:hxhydxdy}:
\begin{equation}
\label{eq:htemp}
H = \left[H^{(1)}_xE^{(1)}_x\right]^TH^{(1)}_yE^{(1)}_x + \left[H^{(2)}_xE^{(2)}_x\right]^TH^{(2)}_yE^{(2)}_x.
\end{equation}
But ($H^{(i)}_2 = H_2, i=1,2$)
\[
H^{(i)}_yE^{(i)}_x = \left[H_2\otimes I^{(i)}_1\right]\left[I_2\otimes E^{(i)}_1\right].
\]
Applying Lemma~\ref{lemma:kron1} twice:
\[
\left[H_2\otimes I^{(i)}_1\right]\left[I_2\otimes E^{(i)}_1\right] = H_2\otimes E^{(i)}_1 =
\left[I_2\otimes E^{(i)}_1\right]\left[H_2\otimes I_1\right].
\]
Define
\[
H_y \equiv H_2\otimes I_1.
\]
Hence,
\[
H^{(i)}_yE^{(i)}_x = E^{(i)}_xH_y, \quad i=1,2.
\]
Using this relation in (\ref{eq:htemp}) yields
\[
H = \left(\left[E^{(1)}_x\right]^TH^{(1)}_xE^{(1)}_x + \left[E^{(2)}_x\right]^TH^{(2)}_xE^{(2)}_x\right)H_y.
\]
But
\begin{align*}
H^{(i)}_x &= I_2\otimes H^{(i)}_1 \\
E^{(i)}_x &= I_2\otimes E^{(i)}_1.
\end{align*}
Thus,
\[
H = \left[I_2\otimes H_1 \right]H_y,
\]
where
\begin{equation}
\label{eq:h1x}
\nonumber
H_1 \equiv \left[E^{(1)}_1\right]^TH^{(1)}_1E^{(1)}_1 + \left[E^{(2)}_1\right]^TH^{(2)}_1E^{(2)}_1 \in \mathbb{R}^{(N_1+1)\times (N_1+1)}.
\end{equation}
Clearly, $H_1$ is the one-dimensional norm defined by (\ref{eq:h}).  Consequently, 
all of the results pertaining to $H$ in sections~\ref{sec:hstructure},~\ref{sec:hblockinv} 
also apply to $H_1$. 
 
Let
\[
H_x \equiv I_2\otimes H_1.
\]
Then
\[
H=H_xH_y.
\]
Furthermore,
\begin{align*}
H_xH_y &= \left[I_2\otimes H_1\right]\left[H_2\otimes I_1\right] \\
       &= H_2\otimes H_1 \\
	   &= \left[H_2\otimes I_1\right]\left[I_2\otimes H_1\right] \\
	   &= H_yH_x,
\end{align*}
which shows that the norm $H$ (\ref{eq:2dinnerprodembedx})
inherits the structure and behavior of the norms $H^{(i)}$ defined on the
subdomains $\Omega_i$. \hfill$\Box$

\begin{prop}
\label{prop:2ddxdystructx}
Let $D_x$ and $D_y$ be as in Definition~\ref{def:2ddiffopembedx}.  If $D^{(1)}_2 = D^{(2)}_2 \equiv D_2$,
then
\begin{align*}
D_x &= I_2 \otimes D_1 \\
D_y &= D_2 \otimes I_1 \\
D_xH_y &= H_yD_x \\
D_yH_x &= H_xD_y,
\end{align*}
where $D_1$ is the one-dimensional difference operator defined in (\ref{eq:diffopembed}).
\end{prop}

\noindent{\bf Proof}:
By the definition of $D_x$:
\begin{align*}
D_x &= H^{-1}E^T_xH^{(+)}D^{(+)}_xE_x \\
    &= H^{-1}
\left(
\left[H^{(1)}_yE^{(1)}_x\right]^TH^{(1)}_xD^{(1)}_xE^{(1)}_x + 
\left[H^{(2)}_yE^{(2)}_x\right]^TH^{(2)}_xD^{(2)}_xE^{(2)}_x  
\right).
\end{align*}
But $H^{(1)}_2 = H^{(2)}_2 = H_2$ according to Proposition~\ref{prop:2ddiffopembedx}.  Thus
\[
\left[H^{(i)}_yE^{(i)}_x\right]^T = \left[E^{(i)}_xH_y\right]^T = H_y\left[E^{(i)}_x\right]^T,\quad i=1,2,
\]
which implies ($H=H_xH_y=H_yH_x$)
\begin{align*}
D_x &= H^{-1}_x
\begin{pmatrix}
\left[E^{(1)}_x\right]^TH^{(1)}_xD^{(1)}_xE^{(1)}_x + \left[E^{(2)}_x\right]^TH^{(2)}_xD^{(2)}_xE^{(2)}_x 
\end{pmatrix} \\
     &= I_2\otimes D_1 \qquad \left[(\ref{eq:dxdy}), (\ref{eq:hxhydef}), (\ref{eq:expart})\right],
\end{align*} 
where 
\begin{align*}
D_1 &\equiv H^{-1}_1\left(\left[E^{(1)}_1\right]^TH^{(1)}_1D^{(1)}_1E^{(1)}_1 + 
                         \left[E^{(2)}_1\right]^TH^{(2)}_1D^{(2)}_1E^{(2)}_1\right) \\
    &=  H^{-1}_1E^T_1H^{(+)}_1D^{(+)}_1E_1,
\end{align*}
which is the corresponding one-dimensional difference operator (\ref{eq:diffopembed}).  Hence,
we have recovered the structure of (\ref{eq:dxdy}).

Similarly, for $D_y$:
\[
D_y = H^{-1}E^T_xH^{(+)}D^{(+)}_yE_x
= H^{-1}E^T_x
\begin{pmatrix}
H^{(1)}D^{(1)}_yE^{(1)}_x  \\
H^{(2)}D^{(2)}_yE^{(2)}_x 
\end{pmatrix}.
\]
But $H^{(1)}_2 = H^{(2)}_2 = H_2$ by necessity.  It is therefore natural to also require that 
$D^{(1)}_2 = D^{(2)}_2 = D_2$.  Hence,
\[
D^{(i)}_y = D_2 \otimes I^{(i)}_1,
\]
and so
\[
D^{(i)}_yE^{(i)}_x = E^{(i)}_x\left[D_2 \otimes I_1\right].
\]
Thus,
\[
D_y = H^{-1}E^T_x
\begin{pmatrix}
H^{(1)}E^{(1)}_x  \\
H^{(2)}E^{(2)}_x 
\end{pmatrix}D_2 \otimes I_1=
H^{-1}\left[E^T_xH^{(+)}E_x\right]D_2 \otimes I = D_2 \otimes I_1
\]
in complete agreement with (\ref{eq:dxdy}).  Finally, by Lemma~\ref{lemma:hxhydxdy}:
\begin{align*}
D_xH_y &= H_yD_x \\
D_yH_x &= H_xD_y.
\end{align*}
All claims have thus been established. \hfill$\Box$

\begin{prop}
\label{prop:2dembedxstructx}
Let $E_x:V\rightarrow V_+$. Then
\[
E^+_x=E^*_x = \left(I_2\otimes \left[E^{(1)}_1\right]^*\ I_2\otimes \left[E^{(2)}_1\right]^*\right).
\]
\end{prop}

\noindent{\bf Proof:}
In Proposition~\ref{prop:2dhstructx} it was shown that
\[
H^{(i)}_yE^{(i)}_x
=
E^{(i)}_xH_y ,\quad i=1,2,
\]
where
\begin{align*}
H^{(i)}_y &= H_2 \otimes I^{(i)}_1  \\
H_y &= H_2 \otimes I_1.
\end{align*}
Hence , from (\ref{eq:pseudoembedx}):
\begin{align*}
E^+_x = E^*_x = H^{-1}E^T_xH^{(+)} &= H^{-1}_xE^T_xH^{(+)}_x \\
                                   &=
\begin{pmatrix}
H^{-1}_x\left[E^{(1)}_x\right]^TH^{(1)}_x &H^{-1}_x\left[E^{(2)}_x\right]^TH^{(2)}_x
\end{pmatrix}.
\end{align*}
But
\begin{align*}
H_x &= I_2\otimes H_1 \\
H^{(i)}_x &= I_2\otimes H^{(i)}_1 \\
E^{(i)}_x &= I_2\otimes E^{(i)}_1.
\end{align*}
By Lemmas~\ref{lemma:kron1},~\ref{lemma:kron2},~\ref{lemma:kron3}:
\[
H^{-1}_x\left[E^{(i)}_x\right]^TH^{(i)}_x = I_2\otimes H^{-1}_1\left[E^{(i)}_1\right]^TH^{(i)}_1.
\]
Finally, we observe that $E^{(i)}:V\rightarrow V_{(i)}$, where $V$ and $V_{(i)}$ are vector
spaces with inner products represented by $H_1$ and $H^{(i)}_1$, i.~e.,
\[
H^{-1}_1\left[E^{(i)}_1\right]^TH^{(i)}_1 = \left[E^{(i)}_1\right]^*,
\]
which concludes the proof.  \hfill$\Box$

\begin{rem}
\label{rem:partialadjx}
For restricted full norms $H^{(1)}_1$ and $H^{(2)}_1$ one has
\begin{equation}
\label{eq:e1e2adjx}
\nonumber
\left[E^{(1)}_1\right]^* =
\begin{pmatrix}
\tilde{I}^{(1)}_1 &0    \\
0   &\chi \\
0   &0
\end{pmatrix},
\quad
\left[E^{(2)}_1\right]^* =
\begin{pmatrix}
0      &0 \\
1-\chi &0 \\
0      &\tilde{I}^{(2)}_1
\end{pmatrix},
\end{equation}
in complete agreement with (\ref{eq:adjstructure}); 
$\tilde{I}^{(i)}_1\in\mathbb{R}^{N^{(i)}_1\times N^{(i)}_1}, i=1,2$.  Hence,
the adjoint of the embedding operator $E:V\rightarrow V_+$, can be viewed as 
an averaging operator.  The one-dimensional expression (\ref{eq:adjstructure}) 
is formally recovered by setting $N_2=0$.  But
\[
\left[E^{(i)}_1\right]^*E^{(i)}_1 \neq I^{(i)}_1 \quad  \Longrightarrow \quad \left[E^{(i)}_1\right]^+ \neq \left[E^{(i)}_1\right]^*.
\]
As a final remark, it should be noted that $D_1$ is given by (\ref{eq:dembedstructure}) if $H^{(i)}, i=1,2$, 
are restricted full norms.  \hfill$\Box$
\end{rem}

\subsection{Two-block difference operators, case 2}
\label{sec:twoblock2diffop}
In this section we will construct difference operators for $\Omega=\Omega_1\cup\Omega_2$ where
\[
\Omega_1 = [0,1]\times[0,1/2], \quad
\Omega_2 = [0,1]\times[1/2,1].
\]
Hence, the interface between the domains is located at $y=0.5$.  The grid points are 
defined as
\begin{align*}
\Omega_1: \quad (x_i, y_j) &= (ih_1, jh^{(1)}_2), \quad 0\leq i\leq N_1, 0\leq j \leq N^{(1)}_2 \\
\Omega_2: \quad (x_i, y_j) &= (ih_1,0.5 + jh^{(2)}_2), \quad 0\leq i\leq N_1, 0\leq j \leq N^{(2)}_2,
\end{align*}
where the mesh sizes are given by
\[
h^{(1)}_1 = h^{(2)}_1 \equiv h_1 = \frac{1}{N_1}, \quad
h^{(j)}_2 = \frac{1}{2N^{(j)}_2}, \quad
j=1,2.
\]
This time, the mesh sizes $h^{(i)}_1$ are the same across the interface $y=0.5$.  
\begin{figure}
\centering
\includegraphics[width=7cm]{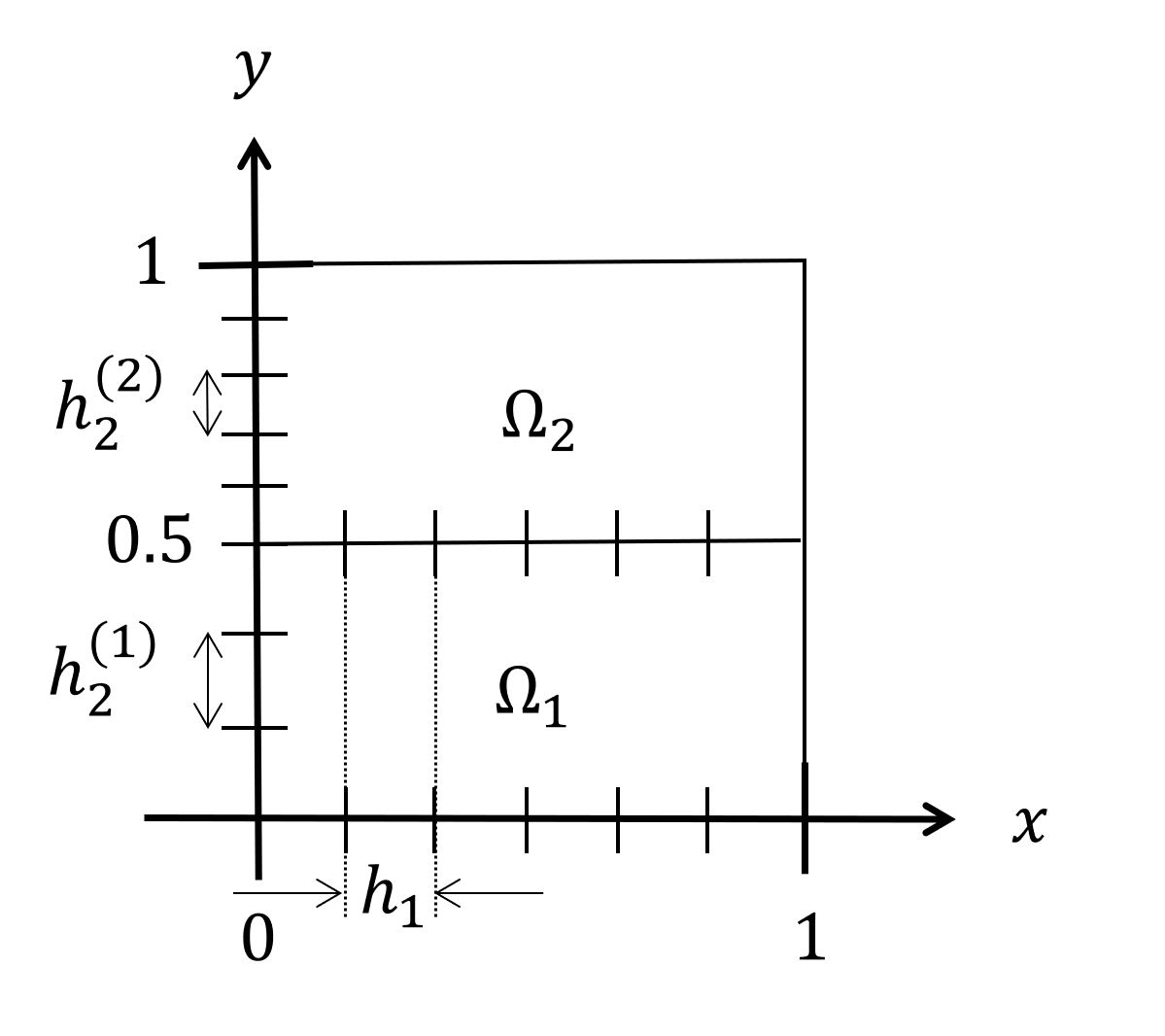}
\caption{Two blocks, case 2}
\label{fig:twoblockscase2}
\end{figure}
Definitions (\ref{eq:opdefi}) - (\ref{eq:2dhplus}) remain unchanged with
\[
c_j = r_j = (N_1+1)(N^{(j)}_2+1).
\]

\subsubsection{The embedding operator $E_y$}
\label{sec:eyembedding}
Let $u, v$ be grid vectors on $\Omega=\Omega_1\cup\Omega_2$ as defined in (\ref{eq:2dstate}):
\[
u, v \in \mathbb{R}^{(N_1+1)(N_2+1)}, \quad N_2 \equiv N^{(1)}_2 + N^{(2)}_2.
\]
We still traverse $\Omega$ horizontally and then vertically, that is
\[
u_j, v_j \in \mathbb{R}^{N_1+1}, \quad j=0, \ldots, N_2.
\]
The embedding $E_y:\mathbb{R}^{(N_1+1)(N_2+1)}\rightarrow \mathbb{R}^{(N_1+1)(N_2+2)}$ will 
be different, however:
\begin{equation}
\label{eq:eypart}
E_y = 
\begin{pmatrix}
E^{(1)}_y \\
E^{(2)}_y
\end{pmatrix}, \quad  E^{(i)}_y \equiv E^{(i)}_2 \otimes I_1, \quad 
I_1\in\mathbb{R}^{(N_1+1)\times(N_1+1)},
\end{equation}
where $E^{(1)}_2$ and $E^{(2)}_2$ are defined by (\ref{eq:e1embed}), (\ref{eq:e2embed}) replacing
$N^{(1)}\rightarrow N^{(1)}_2$, $N^{(2)}\rightarrow N^{(2)}_2$ and $N \rightarrow N_2$ since 
$\Omega_1$ and $\Omega_2$ are joined in the $y$-direction.  The inner product of 
Definition~\ref{def:innerprodembed} is given by
\begin{equation}
\label{eq:2dinnerprodembedy}
(u,v)_H \equiv (E_yu)^TH^{(+)}E_yv \quad \Longleftrightarrow \quad H = E^T_yH^{(+)}E_y.
\end{equation}
Hence, $E_y$ is a mapping between two inner product spaces $V$ and $V_+$ and thus
\[
E^*_y = H^{-1}E^T_yH^{(+)},
\]
i.~e.,
\begin{equation}
\label{eq:pseudoembedy}
E^*_yE_y = I \quad \Longrightarrow \quad E^+_y=E^*_y.
\end{equation}

\subsubsection{Multiblock difference operators $D_x$ and $D_y$}
\label{sec:2dmultiblocky}
Let $D^{(+)}_x,D^{(+)}_y:V_+ \rightarrow V_+$ be as in Definition~\ref{def:2ddiffopembedx}:

\begin{define}
\label{def:2ddiffopembedy}
Given the inner product space $V$ with inner product (\ref{eq:2dinnerprodembedy}), the difference
operators $D_x,D_y:V\rightarrow V$ are defined as
\begin{align*}
D_x &\equiv H^{-1}E^T_yH^{(+)}D^{(+)}_xE_y \\
D_y &\equiv H^{-1}E^T_yH^{(+)}D^{(+)}_yE_y.
\end{align*}
\ \hfill$\Box$
\end{define}
The two-dimensional version of Proposition~(\ref{prop:diffopembed}) when the domains $\Omega_1$ 
and $\Omega_2$ are joined along $y=0.5$ reads:

\begin{prop}
\label{prop:2ddiffopembedy}
Let $D_x,D_y:V\rightarrow V$ be as in Definition~\ref{def:2ddiffopembedy}.  Then
\begin{itemize}
\item[(i)] $D_x,D_y$ satify summation by parts with respect to the inner product 
(\ref{eq:2dinnerprodembedy}):
\begin{align*}
(u,D_xv) &= \langle u,v\rangle_2 - \langle u,v\rangle_4 - (D_xu,v) \\
(u,D_yv) &= \langle u,v\rangle_3 - \langle u,v\rangle_1 - (D_yu,v),
\end{align*}
iff $H^{(1)}_1 = H^{(2)}_1 = H_1 \in \mathbb{R}^{(N_1+1)\times (N_1+1)}$, where
the one-dimensional norm $H_1$ is that of (\ref{eq:h2d}); $H_2$ corresponds to 
(\ref{eq:h}).
\item[(ii)] $D_x,D_y$ are consistent approximations of $\partial/\partial x$ and $\partial/\partial y$.
\item[(iii)] $D_x=E^+_yD^{(+)}_xE_y$ and $D_y=E^+_yD^{(+)}_yE_y$.
\end{itemize}
\end{prop}

\noindent
{\bf Proof:}
Define
\begin{equation}
\label{eq:uembed2dy}
\nonumber
u^{(e)} \equiv E_yu =
\begin{pmatrix}
E^{(1)}_yu \\
E^{(2)}_yu
\end{pmatrix}
\equiv
\begin{pmatrix}
u^{(1)} \\
u^{(2)}
\end{pmatrix}.
\end{equation}
This is also a row-ordered embedding of $u$, but it is different from (\ref{eq:uembed2dx}):
\begin{equation}
\label{eq:2dconstrainty}
\nonumber
u^{(1)} =
\begin{pmatrix}
u_0 \\
\vdots \\
u_{N^{(1)}_2} 
\end{pmatrix} \in \mathbb{R}^{(N_1+1)(N^{(1)}_2+1)}, \quad
u^{(2)} =
\begin{pmatrix}
u_{N^{(1)}_2} \\
u_{N^{(1)}_2+1} \\
\vdots \\
u_{N_2} 
\end{pmatrix} \in \mathbb{R}^{(N_1+1)(N^{(2)}_2+1)}.
\end{equation}
As in the previous case, the operators $D^{(i)}_x$ and  $D^{(i)}_y$ satisfy summation by parts in their
respective domains:
\begin{align*}
(u^{(i)},D^{(i)}_xv^{(i)})_{H^{(i)}} &= (u^{(i)}[N_1,:],v^{(i)}[N_1,:])_{H^{(i)}_2} 
                                      - (u^{(i)}[0,:],v^{(i)}[0,:])_{H^{(i)}_2} \\
                                     &- (D^{(i)}_xu^{(i)},v^{(i)})_{H^{(i)}},  \quad i=1,2.
\end{align*}
Adding the two equations and using the definition of $D_x$:
\begin{align*}
(u,D_xv)   &= (u^{(1)}[N_1,:],v^{(1)}[N_1,:])_{H^{(1)}_2} + (u^{(2)}[N_1,:],v^{(2)}[N_1,:])_{H^{(2)}_2} \\
           &- (u^{(1)}[0,:],v^{(1)}[0,:])_{H^{(1)}_2} - (u^{(2)}[0,:],v^{(2)}[0,:])_{H^{(2)}_2} \\
           &- (D_xu,v).
\end{align*}
By construction, the $i$th column of $u^{(e)}$ satisfies 
\[
\begin{pmatrix}
u^{(1)} \\
u^{(2)}
\end{pmatrix}[i,\!:] = E_2u^i, \quad E_2 = 
\begin{pmatrix}
E^{(1)}_2 \\
E^{(2)}_2
\end{pmatrix}, \quad 0\leq i \leq N_1,
\]
whence, by (\ref{eq:h}):
\[
(u^{(1)}[i,:],v^{(1)}[i,:])_{H^{(1)}_2} + (u^{(2)}[i,:],v^{(2)}[i,:])_{H^{(2)}_2} =
(u^i,u^i)_{H_2}.
\]
This proves the first assertion.

To prove the second claim, we note that
\begin{align*}
(u,D_yv)   &= (u^{(2)}_{N^{(2)}_2},v^{(2)}_{N^{(2)}_2})_{H^{(2)}_1} - (u^{(2)}_0,v^{(2)}_0)_{H^{(2)}_1} \\ 
           &+ (u^{(1)}_{N^{(1)}_2},v^{(1)}_{N^{(1)}_2})_{H^{(1)}_1} - (u^{(1)}_0,v^{(1)}_0)_{H^{(1)}_1} \\
           &- (D_yu,v).
\end{align*}
But
\[
u^{(1)}_{N^{(1)}_2} = u^{(2)}_0, \quad v^{(1)}_{N^{(1)}_2} = v^{(2)}_0,
\]
whence the middle scalar products cancel out iff $H^{(1)}_1=H^{(2)}_1=H_1$.  But
\[
u^{(1)}_0 = u_0, \quad u^{(2)}_{N^{(2)}_2} = u_{N_2}.
\]
The grid vectors $v^{(1)}$ and $v^{(2)}$ satisfy identical relations.  This proves the second claim. 

The second assertion of the proposition can be proved in the same manner as Proposition (\ref{prop:diffopembed}).  
The third claim follows directly from the definition of $D_y$ and (\ref{eq:pseudoembedy}).
\hfill$\Box$

\subsubsection{Structure of $H$, $D_x$, $D_y$ and $E^*_y$}
\label{sec:structurey}
This section summarizes some structural results for $H$, $D_x$, $D_y$ and $E^*_y$.

\begin{prop}
\label{prop:2dhstructy}
Let $H$ be as in (\ref{eq:2dinnerprodembedy}).  If $H^{(1)}_1 = H^{(2)}_1 \equiv H_1$,
then
\begin{align*}
H   &= H_xH_y = H_yH_x \\
H_x &= I_2 \otimes H_1 \\
H_y &= H_2\otimes I_1,
\end{align*}
where $H_2$ is the one-dimensional norm defined in (\ref{eq:h}).
\end{prop}

\noindent
{\bf Proof}:
By (\ref{eq:hxhydef}):
\[
H^{(i)}_x = I^{(i)}_2\otimes H_1, \quad
H^{(i)}_y = H^{(i)}_2\otimes I_1,
\]
where we used
\[
H^{(1)}_1 = H^{(2)}_1 = H_1 \in \mathbb{R}^{(N_1+1)\times (N_1+1)}.
\]
This is a necessary condition for summation by parts to hold in the $y$-direction.
Let $E^{(i)}_y$ be defined as in (\ref{eq:eypart}):
\[
E^{(i)} = E^{(i)}_2\otimes I_1, \quad i=1,2.
\]
Then
\begin{align*}
H^{(i)}_xE^{(i)}_y &= \left[I^{(i)}_2\otimes H_1\right]\left[E^{(i)}_2\otimes I_1\right] \\
                 &= E^{(i)}_2\otimes H_1 \\
                 &= \left[E^{(i)}_2\otimes I_1\right]\left[I_2\otimes H_1\right] \\
				 &= E^{(i)}_yH_x, \quad i=1,2,
\end{align*}
where
\[
H_x \equiv I^{(2)}\otimes H_1.
\]
Thus,
\begin{align*}
H &= 
\left[E^{(1)}_y \right]^TH^{(1)}_yH^{(1)}_xE^{(1)}_y + 
\left[E^{(2)}_y \right]^TH^{(2)}_yH^{(2)}_xE^{(2)}_y \\
  &=
\left(\left[E^{(1)}_y \right]^TH^{(1)}_yE^{(1)}_y + 
\left[E^{(2)}_y \right]^TH^{(2)}_yE^{(2)}_y \right)H_x \\
  &= \left[\left(\left[E^{(1)}_2\right]^TH^{(1)}_2E^{(1)}_2 + \left[E^{(2)}_2\right]^TH^{(2)}_2E^{(2)}_2\right)\otimes I_1\right]H_x \\
  &=\left[H_2 \otimes I_1\right]H_x,
\end{align*}
where
\[
H_2 \equiv \left[E^{(1)}_2\right]^TH^{(1)}_2E^{(1)}_2 + \left[E^{(2)}_2\right]^TH^{(2)}_2E^{(2)}_2
\]
is identical to the one-dimensional inner product (\ref{eq:h}).  Hence, we define
\[
H_y \equiv H_2 \otimes I_1.
\]
Finally,
\[
H_xH_y = H_yH_x,
\]
which follows immediately from the definitions of $H_x$, $H_y$ and 
Lemma~\ref{lemma:kron1}.  \hfill$\Box$

\begin{prop}
\label{prop:2ddxdystructy}
Let $D_x$ and $D_y$ be as in Definition~\ref{def:2ddiffopembedy}.  If $D^{(1)}_1 = D^{(2)}_1 \equiv D_1$,
then
\begin{align*}
D_x &= I_2 \otimes D_1 \\
D_y &= D_2 \otimes I_1 \\
D_xH_y &= H_yD_x  \\
D_yH_x &= H_xD_y ,
\end{align*}
where $D_2$ is the one-dimensional difference operator defined in (\ref{eq:diffopembed}).
\end{prop}

\noindent
{\bf Proof}: By Definition~\ref{def:2ddiffopembedy}
\[
D_x = H^{-1}E^T_yH^{(+)}D^{(+)}_xE_y = 
H^{-1}E^T_y
\begin{pmatrix}
H^{(1)}D^{(1)}_xE^{(1)}_y \\
H^{(2)}D^{(2)}_xE^{(2)}_y
\end{pmatrix},
\]
$E_y$ is given by (\ref{eq:eypart}).  The operator
$D_x^{(i)}$ is defined as $[$(\ref{eq:dxdy})$]$:
\[
D^{(i)}_x =  I^{[i)}_2 \otimes D^{(i)}_1 =  I^{[i)}_2  \otimes D_1,
\]
where we used the assumption $D^{(1)}_1=D^{(2)}_1=D_1$.  This is not very restrictive, since 
$H^{(1)}_1 = H^{(2)}_1 = H_1$ by necessity.  Hence,
\begin{align*}
D^{(i)}_xE^{(i)}_y &= \left[I^{(i)}_2  \otimes D_1\right]\left[E^{(i)}_2 \otimes I_1\right] \\
                   &= E^{(i)}_2 \otimes D_1 \\
			   	   &= \left[E^{(i)}_2 \otimes I_1\right]\left[I_2\otimes D_1\right] \\
				   &= E^{(i)}_y\left[I_2\otimes D_1\right], \quad i=1,2. 
\end{align*}
Lemma~\ref{lemma:kron1} was implicitly invoked in the previous calculations.  Thus
\[
D_x = 
H^{-1}E^T_y
\begin{pmatrix}
H^{(1)}E^{(1)}_y \\
H^{(2)}E^{(2)}_y
\end{pmatrix}\left[I_2\otimes D_1\right] =
I_2\otimes D_1.
\]

Next, consider
\begin{align*}
D_y &= H^{-1}E^T_yH^{(+)}D^{(+)}_yE_y \\
    &= H^{-1}
\left(
\left[H^{(1)}_xE^{(1)}_y\right]^TH^{(1)}_yD^{(1)}_yE^{(1)}_y +
\left[H^{(2)}_xE^{(2)}_y\right]^TH^{(2)}_yD^{(2)}_yE^{(2)}_y
\right),
\end{align*}
where we used $H^{(i)}=H^{(i)}_xH^{(i)}_y=H^{(i)}_yH^{(i)}_x$.  Since $H^{(i)}_1=H_1$:
\begin{align*}
H^{(i)}_xE^{(i)}_y &= \left[I^{(i)}_2 \otimes H_1\right]\left[E^{(i)}_2 \otimes I_1\right] \\
                 &= E^{(i)}_2\otimes H_1 \\
				 &= \left[E^{(i)}_2 \otimes I_1\right]\left[I_2 \otimes H_1\right] \\
				 &= E^{(i)}_yH_x.
\end{align*}
Thus,
\begin{align}
D_y &= H^{-1}
\left(
\left[E^{(1)}_yH_x\right]^TH^{(1)}_yD^{(1)}_yE^{(1)}_y +
\left[E^{(2)}_yH_x\right]^TH^{(2)}_yD^{(2)}_yE^{(2)}_y
\right) \nonumber \\
   &=H^{-1}_y
\left(
\left[E^{(1)}_y\right]^TH^{(1)}_yD^{(1)}_yE^{(1)}_y +
\left[E^{(2)}_y\right]^TH^{(2)}_yD^{(2)}_yE^{(2)}_y
\right), \label{eq:dytemp}
\end{align}
since $H=H_xH_y$ by Proposition~\ref{prop:2dhstructy}.  But
\[
H^{(i)}_y = H^{(i)}_2\otimes I_1, \quad D^{(i)}_y = D^{(i)}_2\otimes I_1, \quad E^{(i)}_y= E^{(i)}_2\otimes I_1,
\]
by definition.  Furthermore, by Proposition~\ref{prop:2dhstructy} and Lemma~\ref{lemma:kron2}:
\[
H^{-1}_y = H^{-1}_2\otimes I_1.
\]
Hence, applying Lemmas~\ref{lemma:kron1},\ref{lemma:kron3} to (\ref{eq:dytemp}):
\[
D_y = D_2\otimes I_1, \quad 
D_2 \equiv H^{-1}_2\left(\left[E^{(1)}_2\right]^TH^{(1)}_2D^{(1)}_2E^{(1)}_2 + 
                         \left[E^{(2)}_2\right]^TH^{(2)}_2D^{(2)}_2E^{(2)}_2\right).
\]
Obviously, $D_2$ is identical to the one-dimensional expression (\ref{eq:diffopembed}).  To
conclude this discussion, we observe that
\begin{align*}
D_xH_y &= H_yD_x \\
D_yH_x &= H_xD_y. 
\end{align*}
\hfill$\Box$

\begin{prop}
\label{prop:2dembedystructy}
Let $E_y:V\rightarrow V_+$. Then
\[
E^+_y=E^*_y = \left(\left[E^{(1)}_2\right]^* \otimes I_1 \ \left[E^{(2)}_2\right]^* \otimes I_1\right).
\]
\end{prop}

\noindent{\bf Proof}:
From the proof of Proposition~\ref{prop:2dhstructy} we recall
\[
H^{(i)}_xE^{(i)}_y = E^{(i)}_yH_x, \quad H_x = I_2\otimes H_1, \quad i=1,2,
\]
where $E_y$ is defined by (\ref{eq:eypart}).  Thus, $[$(\ref{eq:pseudoembedy})$]$:
\begin{align*}
E^+_y &= E^*_y  \\
      &=
\begin{pmatrix}
H^{-1}\left[E^{(1)}_y\right]^TH^{(1)} &H^{-1}\left[E^{(2)}_y\right]^TH^{(2)}
\end{pmatrix} \\
	  &=
\begin{pmatrix}
H^{-1}_y\left[E^{(1)}_y\right]^TH^{(1)}_y &H^{-1}_y\left[E^{(2)}_y\right]^TH^{(2)}_y
\end{pmatrix}.
\end{align*}
But
\begin{align*}
H^{-1}_y  &= H^{-1}_2\otimes I_1 \\
H^{(i)}_y &= H^{(i)}_2\otimes I_1 \\
E^{(i)}_y &= E^{(i)}_2\otimes I_1.
\end{align*}
Hence,
\[
H^{-1}_y\left[E^{(i)}_y\right]^TH^{(i)}_y = \left[H^{-1}_2\left[E^{(i)}_2\right]^TH^{(i)}_2\right]\otimes I_1.
\]
As in the proof of Proposition~\ref{prop:2dembedxstructx}, we have $E^{(i)}:V\rightarrow V_{(i)}$.  This time,
$V$ and $V_{(i)}$ have inner products represented by $H_2$ and $H^{(i)}_2$.  Thus,
\[
H^{-1}_2\left[E^{(i)}_2\right]^TH^{(i)}_2 = \left[E^{(i)}_2\right]^*.
\]
The proposition has thus been proved.  \hfill$\Box$

\begin{rem}
\label{rem:partialadjy}
Just as in Case~1, for restricted full norms $H^{(1)}_2$ and $H^{(2)}_2$:
\begin{equation}
\label{eq:e1e2adjy}
\nonumber
\left[E^{(1)}_2\right]^* =
\begin{pmatrix}
\tilde{I}^{(1)}_2 &0    \\
0   &\chi \\
0   &0
\end{pmatrix},
\quad
\left[E^{(2)}_2\right]^* =
\begin{pmatrix}
0      &0 \\
1-\chi &0 \\
0      &\tilde{I}^{(2)}_2
\end{pmatrix}.
\end{equation}
The adjoint of $E_y$ acts an averaging operator in the $y$-direction. \hfill$\Box$
\end{rem}
\begin{rem}
\label{rem:struct}
From Sections~\ref{sec:structurex} and \ref{sec:structurey} it is evident that $H$, $D_x$ and $D_y$
have identical structure in both two-block cases.  This result will be used in the next section. 
\hfill$\Box$
\end{rem}

\subsection{Four-block difference operators}
\label{sec:fourblockdiffop}
We conclude the discussion on multiblock difference operators by considering the four-block 
case, in which the unit square $\Omega=[0, 1]\times [0, 1]$ is broken up into four equally 
sized subdomains
\begin{equation}
\label{eq:omega4}
\nonumber
\Omega=\cup\Omega_{ij}, \quad i,j=1,2, 
\end{equation}
where
\begin{align*}
\Omega_{11} &= [0,1/2]\times [0,1/2]  \\
\Omega_{21} &= [1/2,1]\times [0,1/2]  \\
\Omega_{12} &= [0,1/2]\times [1/2,1]  \\
\Omega_{22} &= [1/2,1]\times [1/2,1].
\end{align*}
The individual meshes are defined as
\begin{align*}
\Omega_{11}: \quad (x_k, y_l) &= (kh^{(1)}_1,       lh^{(1)}_2),       \quad 0\leq k\leq N^{(1)}_1, 0\leq l \leq N^{(1)}_2  \\
\Omega_{21}: \quad (x_k, y_l) &= (0.5 + kh^{(2)}_1, lh^{(1)}_2),       \quad 0\leq k\leq N^{(2)}_1, 0\leq l \leq N^{(1)}_2  \\
\Omega_{12}: \quad (x_k, y_l) &= (kh^{(1)}_1,       0.5 + lh^{(2)}_2), \quad 0\leq k\leq N^{(1)}_1, 0\leq l \leq N^{(2)}_2  \\
\Omega_{22}: \quad (x_k, y_l) &= (0.5 + kh^{(2)}_1, 0.5 + lh^{(2)}_2), \quad 0\leq k\leq N^{(2)}_1, 0\leq l \leq N^{(2)}_2. 
\end{align*}
These meshes fulfill the grid matching restrictions of Sections~\ref{sec:twoblock1diffop} 
and \ref{sec:twoblock2diffop}.  Next, we introduce two intermediate partitions of $\Omega$:
\begin{align}
\Omega_j &= \cup_i\Omega_{ij} \label{eq:omegahor} \\
\Omega^i &= \cup_j\Omega_{ij} \label{eq:omegaver}.
\end{align}
The first partition $\Omega_j$ corresponds to joining $\Omega_{ij}$ horizontally, which is 
discussed in Section~\ref{sec:twoblock1diffop}; $\Omega^i$ represents joining the subdomains 
$\Omega_{ij}$ vertically, cf.~Section~\ref{sec:twoblock2diffop}.

\begin{figure}
\centering
\includegraphics[width=7cm]{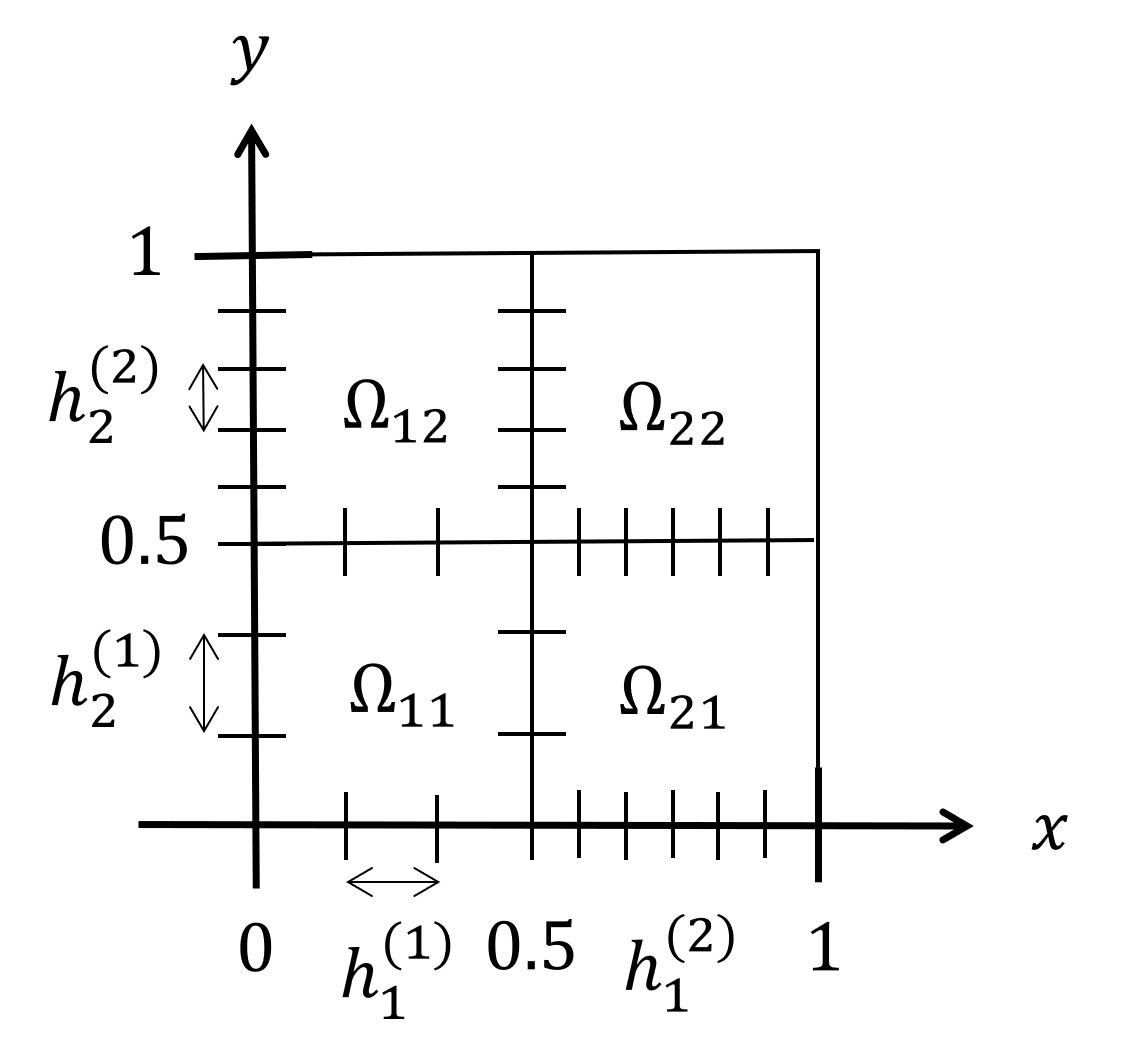}
\caption{Four blocks}
\label{fig:fourblockscase}
\end{figure}

Following the same cadence as in the two-block cases, for each domain we define grid 
vectors (\ref{eq:2dstate}), scalar products (\ref{eq:h2d}) (\ref{eq:2dinnerprod}) and 
difference operators (\ref{eq:dxdy}):
\begin{equation}
\label{eq:opdefij}
\nonumber
\Omega_{ij}: \quad u^{(ij)},v^{(ij)}, \quad H^{(ij)} = H^{(ij)}_xH^{(ij)}_y, \quad D^{(ij)}_x,D^{(ij)}_y,
\end{equation}
where we used Lemma~\ref{lemma:hxhydxdy}.  By (\ref{eq:hxhydef}):
\begin{align*}
H^{(ij)}_x &= I^{(ij)}_2 \otimes H^{(ij)}_1 \\
H^{(ij)}_y &= H^{(ij)}_2 \otimes I^{(ij)}_1.
\end{align*}
The requirement of having matching grid lines at the interfaces $\Omega_{1j}\cap\Omega_{2j}$
and $\Omega_{i1}\cap\Omega_{i2}$ implies the {\em necessary} conditions
\begin{align*}
I^{(ij)}_1 &= I^{(i)}_1 \\
I^{(ij)}_2 &= I^{(j)}_2.
\end{align*}
Furthermore, summation by parts in the two-block cases is possible iff 
\begin{align*}
H^{(ij)}_1 &= H^{(i)}_1 \\
H^{(ij)}_2 &= H^{(j)}_2.
\end{align*}
We thus end up with
\begin{align*}
H^{(ij)}_x &= I^{(j)}_2 \otimes H^{(i)}_1 \\
H^{(ij)}_y &= H^{(j)}_2 \otimes I^{(i)}_1.
\end{align*}
Finally, to preserve the tensor structure of $D_x, D_y$ in Propositions~\ref{prop:2ddxdystructx} 
and \ref{prop:2ddxdystructy}, we made the following {\em sufficient} assumptions on the 
one-dimensional difference operators:
\begin{align*}
D^{(ij)}_1 &= D^{(i)}_1 \\
D^{(ij)}_2 &= D^{(j)}_2.
\end{align*}
Hence,
\begin{align*}
D^{(ij)}_x &= I^{(j)}_2 \otimes D^{(i)}_1 \\
D^{(ij)}_y &= D^{(j)}_2 \otimes I^{(i)}_1.
\end{align*}

\subsubsection{The augmented state spaces $V_{+j}$}
\label{sec:2daugstatex}
We begin by dividing $\Omega$ into two vertically stacked multisets 
$\Omega_{+j}=\Omega_{1j} + \Omega_{2j}$:
\begin{equation}
\label{eq:2duvplusj}
\nonumber
u^{(+j)} \equiv
\begin{pmatrix}
u^{(1j)} \\
u^{(2j)}
\end{pmatrix}
\begin{matrix}
{\scriptstyle r_{1j}} \\
{\scriptstyle r_{2j}}
\end{matrix}, \quad
v^{(+j)} \equiv
\begin{pmatrix}
v^{(1j)} \\
v^{(2j)}
\end{pmatrix}
\begin{matrix}
{\scriptstyle r_{1j}} \\
{\scriptstyle r_{2j}}
\end{matrix}.
\end{equation}
The intermediate (augmented) state spaces $V_{+j}$ are defined as in 
Section~\ref{sec:twoblock1diffop}, and the scalar products are given by: 
\begin{equation}
\label{eq:2dhplusj}
\nonumber
(u^{(+j)},v^{(+j)})_{+j} = \left[u^{(+j)}\right]^TH^{(+j)}v^{(+j)},
 \quad
H^{(+j)} =
\overset{\raise3pt\hbox{$\scriptstyle c_{1j} \ \ \ \ \ \ c_{2j}$}}{
\begin{pmatrix}
H^{(1j)} \\
&H^{(2j)}
\end{pmatrix}}
\begin{matrix}
{\scriptstyle r_{1j}} \\
{\scriptstyle r_{2j}}
\end{matrix},
\end{equation}
where
\[
\begin{array}{l}
c_{ij} = r_{ij} = (N^{(i)}_1+1)(N^{(j)}_2+1).
\end{array}
\]

\subsubsection{The embedding operators $E^{(j)}_x$}
\label{sec:exembeddings}
Let $u^{(j)}, v^{(j)}$ be grid vectors on $\Omega_j$ (\ref{eq:omegahor}) as defined in 
(\ref{eq:2dstate}):
\[
u^{(j)}, v^{(j)} \in \mathbb{R}^{(N_1+1)(N^{(j)}_2+1)}, \quad N_1 \equiv N^{(1)}_1 + N^{(2)}_1.
\]
As usual, $\Omega_j$ is traversed horizontally and then vertically, that is
\[
u^{(j)}_l, v^{(j)}_l \in \mathbb{R}^{N_1+1}, \quad l=0, \ldots, N^{(j)}_2.
\]
Define mappings $E^{(j)}_x:\mathbb{R}^{(N_1+1)(N^{(j)}_2+1)}\rightarrow
\mathbb{R}^{(N_1+2)(N^{(j)}_2+1)}$:
\begin{equation}
\label{eq:expartj}
E^{(j)}_x = 
\begin{pmatrix}
E^{(1j)}_x \\
E^{(2j)}_x
\end{pmatrix}, \quad  E^{(ij)}_x \equiv I^{(j)}_2 \otimes E^{(i)}_1, \quad 
I^{(j)}_2\in\mathbb{R}^{(N^{(j)}_2+1)\times(N^{(j)}_2+1)};
\end{equation}
$E^{(i)}_1$, $i=1,2$, are defined by (\ref{eq:e1embed}) and (\ref{eq:e2embed}).

\begin{rem}
\label{rem:expartj}
The one-dimensional embedding (\ref{eq:epart}) is formally recovered by setting $N^{(j)}_2=0$.
\hfill$\Box$
\end{rem}
The inner products on $\mathbb{R}^{(N_1+1)(N^{(j)}_2+1)}\times\mathbb{R}^{(N_1+1)(N^{(j)}_2+1)}$ 
are then defined as:
\begin{equation}
\label{eq:2dinnerprodembedxj}
(u^{(j)},v^{(j)})_{H^{(j)}} \equiv (E^{(j)}_xu^{(j)})^TH^{(+j)}E^{(j)}_xv^{(j)},
\end{equation}
i.~e.,
\begin{equation}
\label{eq:2dinnerprodembedxja}
H^{(j)} = \left[E^{(j)}_x\right]^TH^{(+j)}E^{(j)}_x.
\end{equation}
We thus conclude that $E^{(j)}_x$ are mappings between the inner product spaces $V_j$ 
and $V_{+j}$, i.~e., $E^{(j)}_x:V_j\rightarrow V_{+j}$.  Hence,
\[
\left[E^{(j)}_x\right]^* = \left[H^{(j)}\right]^{-1}\left[E^{(j)}_x\right]^TH^{(+j)} =
\left[E^{(j)}_x\right]^+.
\]

\subsubsection{Multiblock difference operators $D_x$ and $D_y$}
\label{sec:2dfourblockxj}
Let $D^{(+j)}_x,D^{(+j)}_y:V_{+j} \rightarrow V_{+j}$:
\begin{equation}
\label{eq:2ddplusj}
\nonumber
D^{(+j)}_x = 
\overset{\raise3pt\hbox{$\scriptstyle c_{1j} \ \ \ \ \ \ c_{2j}$}}{
\begin{pmatrix}
D^{(1j)}_x \\
&D^{(2j)}_x
\end{pmatrix}}
\begin{matrix}
{\scriptstyle r_{1j}} \\
{\scriptstyle r_{2j}}
\end{matrix}, \quad
D^{(+j)}_y =
\overset{\raise3pt\hbox{$\scriptstyle c_{1j} \ \ \ \ \ \ c_{2j}$}}{
\begin{pmatrix}
D^{(1j)}_y \\
&D^{(2j)}_y
\end{pmatrix}}
\begin{matrix}
{\scriptstyle r_{1j}} \\
{\scriptstyle r_{2j}}
\end{matrix},
\end{equation}
which are obtained by applying (\ref{eq:2ddplus}) to $\Omega_j=\cup_i\Omega_{ij}$. 
\begin{define}
\label{def:2ddiffopembedxj}
Given the inner product space $V_j$ with inner product (\ref{eq:2dinnerprodembedxj}), the difference
operators $D^{(j)}_x,D^{(j)}_y:V_j\rightarrow V_j$ are defined as
\begin{align*}
D^{(j)}_x &\equiv \left[H^{(j)}\right]^{-1}\left[E^{(j)}_x\right]^TH^{(+j)}D^{(+j)}_xE^{(j)}_x \\
D^{(j)}_y &\equiv \left[H^{(j)}\right]^{-1}\left[E^{(j)}_x\right]^TH^{(+j)}D^{(+j)}_yE^{(j)}_x.
\end{align*}
\ \hfill$\Box$
\end{define}
It follows  immediately from Proposition~\ref{prop:2ddiffopembedx} that $D^{(j)}_x$ and
$D^{(j)}_y$ satisfy summation by parts in their respective domains $\Omega_j$.  Furthermore,
\begin{align*}
H^{(j)}            &= H^{(j)}_xH^{(j)}_y = H^{(j)}_yH^{(j)}_x  \\
D^{(j)}_xH^{(j)}_y &= H^{(j)}_yD^{(j)}_x  \\
D^{(j)}_yH^{(j)}_x &= H^{(j)}_xD^{(j)}_y  \\
H^{(j)}_x          &= I^{(j)}_2 \otimes H_1 \\
H^{(j)}_y          &= H^{(j)}_2 \otimes I_1 \\
D^{(j)}_x          &= I^{(j)}_2 \otimes D_1 \\
D^{(j)}_y          &= D^{(j)}_2 \otimes I_1, 
\end{align*}
where $H_1$, $D_1$ correspond to the one-dimensional operators in (\ref{eq:h}) 
and (\ref{eq:diffopembed}); $I_1\in\mathbb{R}^{(N_1+1)\times (N_1+1)}$.

We now apply the arguments of Section~\ref{sec:twoblock2diffop} to the domains 
$\Omega_j=\Omega_{1j}\cup\Omega_{2j}$.  Hence, $V_+$, $H^{(+)}$, $V$, $E_y:V\rightarrow V_+$ 
and $H$ are defined as in Section~\ref{sec:twoblock2diffop}.  For instance,
\[
H^{(+)} = 
\begin{pmatrix}
H^{(1)} \\
&H^{(2)}
\end{pmatrix},
\]
where $H^{(j)}, j=1,2$ are given by (\ref{eq:2dinnerprodembedxja}).  The scalar product in $V$
can thus be represented as
\begin{equation}
\label{eq:4blockinnerprodembed1}
H = E^T_yH^{(+)}E_y, 
\end{equation}
see (\ref{eq:eypart}) for the definition of $E_y$.  
The operators $D_x,D_y:V\rightarrow V$ are exactly as in Definition~\ref{def:2ddiffopembedy}.  
By Proposition~\ref{prop:2ddiffopembedy}, they satisfy summation by parts.  Finally, 
by Propositions~\ref{prop:2dhstructy} and \ref{prop:2ddxdystructy}:
\begin{align}
H      &= H_xH_y = H_yH_x \nonumber \\
D_xH_y &= H_yD_x \nonumber \\
D_yH_x &= H_xD_y \nonumber \\
H_x    &= I_2 \otimes H_1 \label{eq:2dembedxdystruct}\\
H_y    &= H_2 \otimes I_1 \nonumber \\
D_x    &= I_2 \otimes D_1 \nonumber \\
D_y    &= D_2 \otimes I_1. \nonumber
\end{align}
This time, $H_2$, $D_2$ represent the one-dimensional operators in (\ref{eq:h}) 
and (\ref{eq:diffopembed}); $I_2\in\mathbb{R}^{(N_2+1)\times (N_2+1)}$.

\subsubsection{Multiblock difference operators $D_x$ and $D_y$ revisited}
\label{sec:2dfourblockxi}
Instead of dividing $\Omega$ in the $y$-direction, we take a different route and begin with
two horizontal slabs $\Omega^i$ (\ref{eq:omegaver}), $\Omega=\Omega^1\cup\Omega^2$.  This will lead to augmented state spaces
$V_{i+}$ and $\tilde{V}_+$ with grid vectors $u^{(i+)}$, $\tilde{u}^{(+)}$ and corresponding 
scalar products $H^{(i+)}$ and $\tilde{H}^{(+)}$.  These spaces are {\em not} the same as $V_{+j}$ 
and $V_+$ encountered in the previous section. 
 
Given the state spaces $V_{ij}$, we define scalar products and state vectors in $V_{i+}$ 
represented by the matrices
\[
H^{(i+)} = 
\begin{pmatrix}
H^{(i1)} \\
&H^{(i2)}
\end{pmatrix},
\quad
u^{(i+)} = 
\begin{pmatrix}
u^{(i1)} \\
u^{(i2)}
\end{pmatrix},
\]
from which we derive the intermediate state spaces $\tilde{V}_i$, embedding operators 
$\tilde{E}^{(i)}_y:\tilde{V}_i \rightarrow V_{i+}$ 
\begin{equation}
\label{eq:eyijembed}
\tilde{E}^{(i)}_y = 
\begin{pmatrix}
\tilde{E}^{(i1)}_y \\
\tilde{E}^{(i2)}_y
\end{pmatrix}, \quad  \tilde{E}^{(ij)}_y \equiv E^{(j)}_2 \otimes I^{(i)}_1, \quad 
I^{(i)}_1\in\mathbb{R}^{(N^{(i)}_1+1)\times(N_1^{(i)}+1)},
\end{equation}
and scalar products
\begin{equation}
\label{eq:hitildeinnerprod}
\tilde{H}^{(i)} = \left[\tilde{E}^{(i)}_y\right]^TH^{(i+)}\tilde{E}^{(i)}_y.
\end{equation}
State vectors in $\tilde{V}_i$ are denoted by $\tilde{u}^{(i)}$.

Next, let $\tilde{E}_x:\tilde{V} \rightarrow \tilde{V}_+$ be the embedding of the 
final state space $\tilde{V}$ into $\tilde{V}_+$.  The inner products of $\tilde{V}$ and 
$\tilde{V}_+$ correspond to
\begin{equation}
\label{eq:4blockinnerprodembed2}
\tilde{H} = \tilde{E}_x^T\tilde{H}^{(+)}\tilde{E}_x,
\quad
\tilde{H}^{(+)} = 
\begin{pmatrix}
\tilde{H}^{(1)} \\
&\tilde{H}^{(2)}
\end{pmatrix}.
\end{equation}
State vectors in $\tilde{V}$, $\tilde{V}_+$ are denoted by $\tilde{u}$ and
\[
\tilde{u}^{(+)} = 
\begin{pmatrix}
\tilde{u}^{(1)} \\
\tilde{u}^{(2)}
\end{pmatrix}.
\]

\begin{lemma}
\label{lemma:gridvector}
Let $u\in V$ and $\tilde{u}\in\tilde{V}$ be two grid vectors describing the same 
state $u_{ij}$ defined on $\Omega=\cup\Omega_{ij}$.  Then $u=\tilde{u}$.
\end{lemma}

\noindent
{\bf Proof}:
The grid vector $u$ is obtained from two intermediate grid vectors $u^{(j)}$:
\[
u = 
\begin{pmatrix*}[l]
u^{(1)} \\
u^{(2)}[1\!\!:]
\end{pmatrix*},
\]
where 
\[
u^{(j)}_k = 
\begin{pmatrix*}[l]
u^{(1j)}_k \\
u^{(2j)}_k[1\!\!:]
\end{pmatrix*}, \quad 0\leq k \leq N^{(j)}_2.
\]
Similarly,
\[
\tilde{u}^{(i)} =
\begin{pmatrix*}[l]
u^{(i1)} \\
u^{(i2)}[1\!\!:]
\end{pmatrix*}.
\]
The final grid vector $\tilde{u}$ is obtained by interlacing the block rows of $\tilde{u}^{(1)}$
and $\tilde{u}^{(2)}$ omitting the first element of each block row in $\tilde{u}^{(2)}$:
\[
\tilde{u} =
\begin{pmatrix*}[l]
u^{(11)}_0 \\
u^{(21)}_0[1\!\!:] \\
\vdots \\
u^{(12)}_{N^{(2)}_2} \\
u^{(22)}_{N^{(2)}_2}[1\!\!:]
\end{pmatrix*} =
\begin{pmatrix*}[l]
u^{(1)} \\
u^{(2)}[1\!\!:]
\end{pmatrix*} = u,
\]
which proves the lemma. \hfill$\Box$

Applying the machinery of Sections~\ref{sec:twoblock1diffop} and \ref{sec:twoblock2diffop} to 
the alternate partitioning $\Omega=\Omega^1\cup\Omega^2$ will lead to the same result 
(\ref{eq:2dembedxdystruct}) replacing $H_{x,y}$, $D_{x,y}$ with $\tilde{H}_{x,y}$ and 
$\tilde{D}_{x,y}$ .  Thus, by the last four equations of (\ref{eq:2dembedxdystruct}):
\[
\tilde{H} = H, \quad \tilde{D}_x = D_x, \quad \tilde{D}_y = D_y,
\]
where we also used Lemma~\ref{lemma:gridvector}.  

We conclude this discussion by giving a direct proof that $\tilde{H} = H$.

\begin{lemma}
\label{lem:4blockinnerprodequiv} 
Let $H$ and $\tilde{H}$ be defined by (\ref{eq:4blockinnerprodembed1}) and (\ref{eq:4blockinnerprodembed2}).
Then $\tilde{H} = H$.
\end{lemma}

\noindent{\bf Proof}:  
Using (\ref{eq:hitildeinnerprod}) and (\ref{eq:4blockinnerprodembed2}):
\[
\tilde{H} = 
\sum^2_{i=1}\left[\tilde{E}^{(i)}_x\right]^T\left[\tilde{E}^{(i)}_y\right]^TH^{(i+)}\tilde{E}^{(i)}_y\tilde{E}^{(i)}_x.
\]
Thus, by (\ref{eq:eyijembed}):
\[
\tilde{H} = 
\sum^2_{i,j=1}\left[\tilde{E}^{(i)}_x\right]^T\left[\tilde{E}^{(ij)}_y\right]^TH^{(ij)}\tilde{E}^{(ij)}_y\tilde{E}^{(i)}_x,
\]
where
\[
\tilde{E}^{(ij)}_y = E^{(j)}_2 \otimes I^{(i)}_1, \quad
\tilde{E}^{(i)}_x = I_2 \otimes E^{(i)}_1,
\]
cf.~(\ref{eq:eyijembed}) and (\ref{eq:expart}).  Hence,
\[
\tilde{E}^{(ij)}_y\tilde{E}^{(i)}_x = E^{(j)}_2 \otimes E^{(i)}_1 =
\left(I^{(j)}_2 \otimes E^{(i)}_1\right)\left(E^{(j)}_2 \otimes I_1\right) = E^{(ij)}_xE^{(j)}_y,
\]
where we used (\ref{eq:expartj}) and (\ref{eq:eypart}).  Substituting this into the 
expression for $\tilde{H}$:
\[
\tilde{H} = 
\sum^2_{j=1}\left[E^{(j)}_y\right]^T\left[\sum^2_{i=1}\left[E^{(ij)}_x\right]^TH^{(ij)}E^{(ij)}_x\right]E^{(j)}_y =
\sum^2_{j=1}\left[E^{(j)}_y\right]^TH^{(j)}E^{(j)}_y = H,
\]
which proves the lemma. \ \hfill$\Box$

\section{Numerical results}
\label{sec:numerics}
In this section we verify the embedding method by considering Maxwell's equations in two space dimensions:
\begin{equation}
\label{eq:maxwell_pde}
\begin{alignedat}{4}
	C u_t &= A u_x + B u_y, \quad && x,y \in \Omega, \quad &&&t > 0 \\
	H(x,y,t) &= g(x,y,t), && x,y \in \Gamma, &&& t \geq 0 \\
	u(x,y,t) &= f(x,y), && x,y \in \Omega, &&& t = 0,
\end{alignedat}
\end{equation}
where the solution vector
\[
u = \begin{pmatrix}
E_x \\
H \\
E_y
\end{pmatrix}
\]
contains the $x$ and $y$ components of the electric field $E_{x,y}$ (not to be confused with the embedding 
operators $E_{x,y}$) and the magnetic field $H$.  The coefficient matrices are given by
\begin{equation}
\label{eq:abc_analytic}
\nonumber
	A = \begin{pmatrix}
		0 & 0 & 0 \\
		0 & 0 & -1 \\
		0 & -1 & 0
	\end{pmatrix}, \quad 
	B = \begin{pmatrix}
		0 & 1 & 0 \\
		1 & 0 & 0 \\
		0 & 0 & 0
	\end{pmatrix}, \quad \text{and} \quad 
	C = \begin{pmatrix}
		\epsilon & 0 & 0 \\
		0 & \mu & 0 \\
		0 & 0 & \epsilon
	\end{pmatrix}.
\end{equation}
The material parameters $\epsilon$ and $\mu$ will in general be space and time dependent functions 
but are considered constant in the following computations.  As boundary conditions we specify the 
magnetic field at all boundaries. The spatial domain is $\Omega \subset \mathbb{R}^2$ and its boundary 
is denoted $\Gamma$.

Well-posedness of \eqref{eq:maxwell_pde} can be analyzed using the energy method.  Multiplying \eqref{eq:maxwell_pde} with $u$ and integrating over $\Omega$ leads to
\[
	(u,C u_t) = (u, A u_x) + (u, B u_y).
\]
Integration by parts results in
\begin{equation}
    \label{eq:maxwell_anal_stab}
	\frac{d}{dt} \|u\|^2_C = \int _{\Gamma} 2H \left(E_x n_y - E_y n_x\right) \: ds,
\end{equation}
where $n = (n_x\ n_y)^T$ is the outward pointing normal and $ \|u\|^2_C  \equiv (u,Cu)$.  Inserting the homogeneous version of the boundary conditions, $g(x,y,t) = 0$, implies energy conservation:
\begin{equation}
	\label{eq:maxwell_phys_enest}
	\nonumber
	\frac{d}{dt} \|u\|^2_C = 0,
\end{equation}
which is enough to prove well-posedness.

\begin{figure}
\begin{subfigure}{.5\textwidth}
	\centering
	\includegraphics[width=0.9\textwidth]{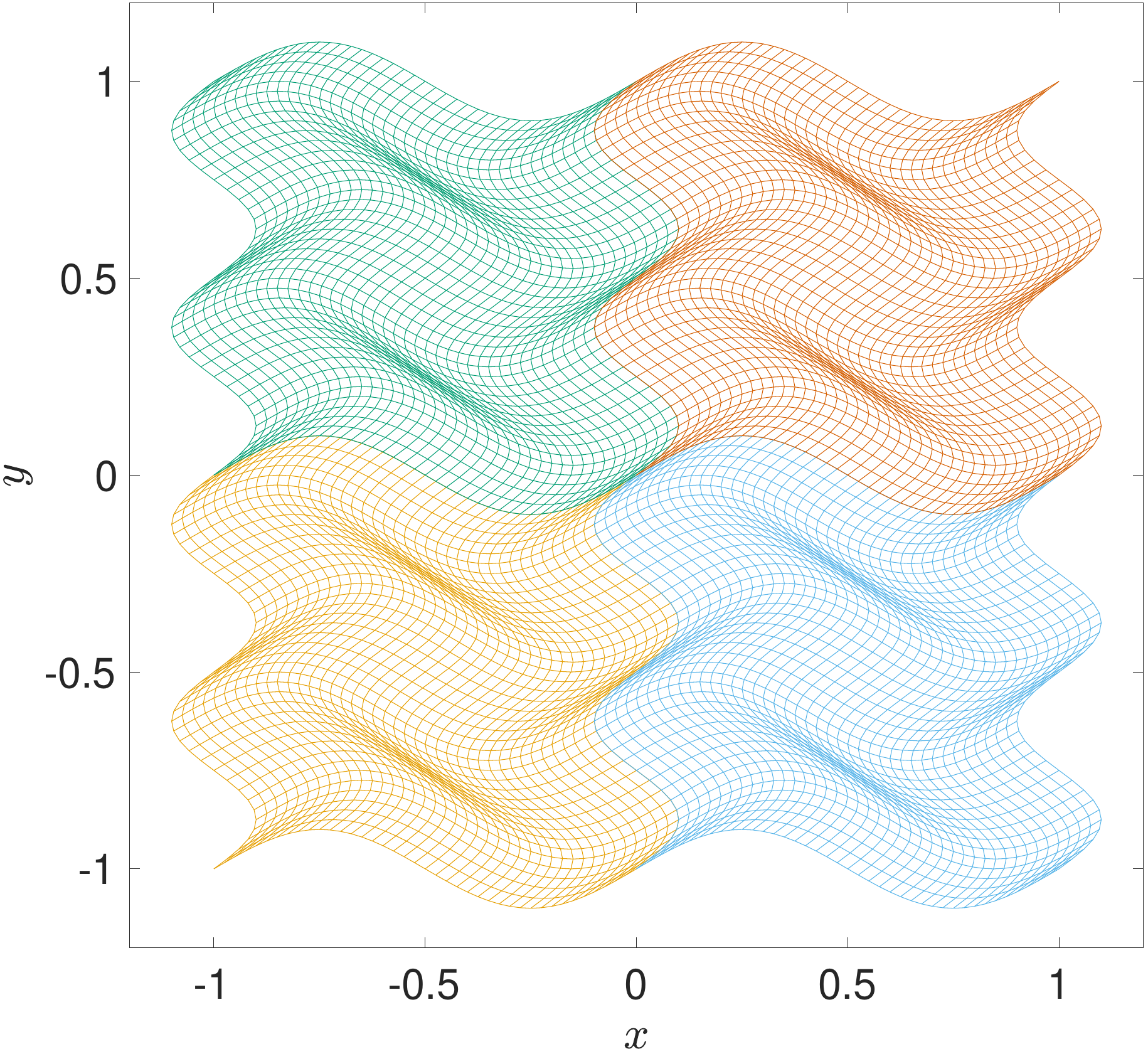}
	\caption{Physical domain $\Omega$.}
	\label{fig:phys_domain}
\end{subfigure}
\begin{subfigure}{.5\textwidth}
	\centering
	\includegraphics[width=0.9\textwidth]{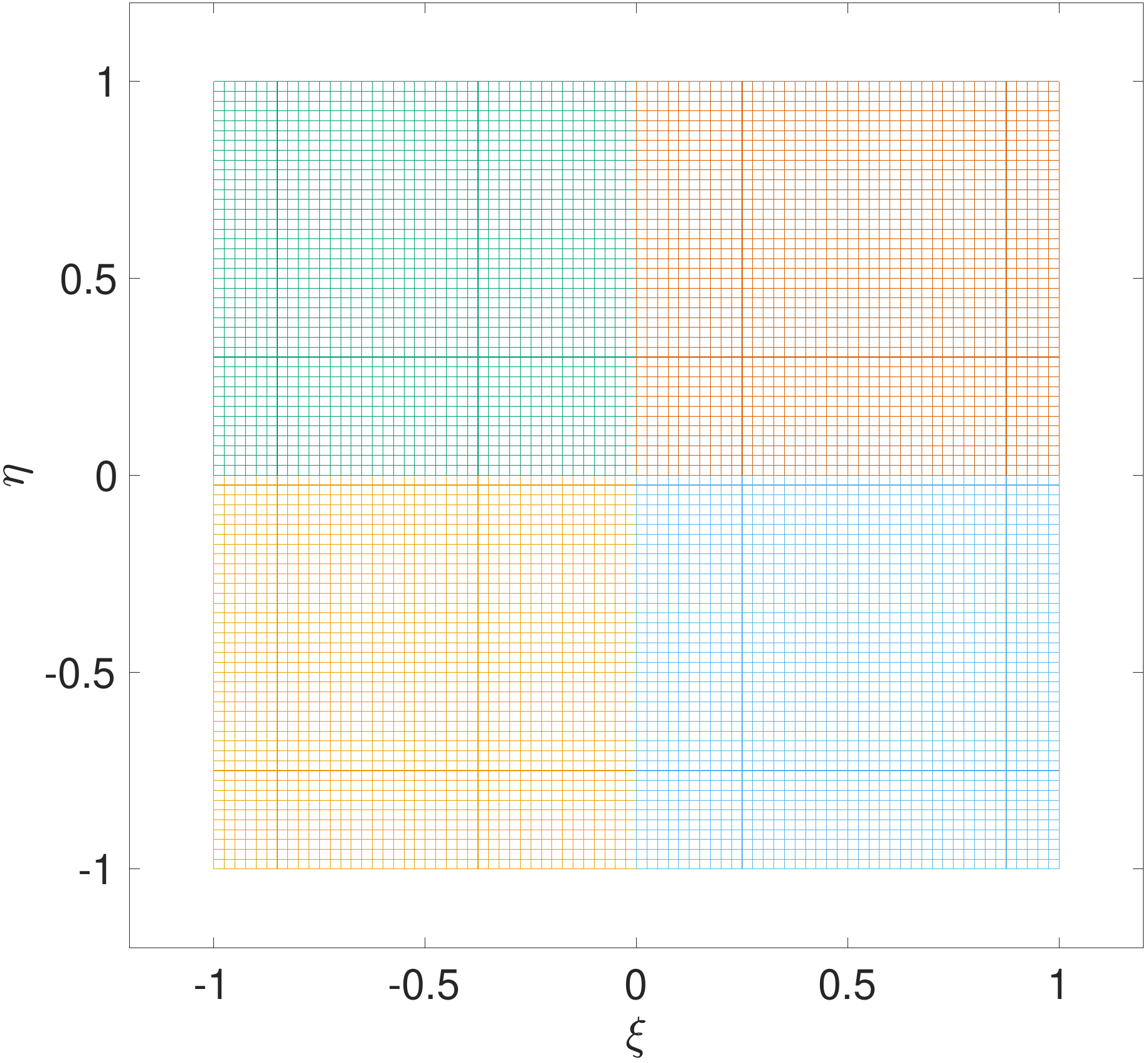}
	\caption{Computational domain $\Omega'$.}
	\label{fig:comp_domain}
\end{subfigure}
\caption{Physical and computational domains}
\label{fig:domain}
\end{figure}

To test the embedding method, we consider the multiblock curvilinear domain shown in Figure \ref{fig:phys_domain}.  
The physical domain is rectified using a reference domain $\Omega' = [-1,1]\times[-1,1]$, see Figure \ref{fig:comp_domain}, 
and a diffeomorphism $x = x(\xi,\eta)$ and $y = y(\xi,\eta)$, where $\xi,\eta \in \Omega'$.  The state 
$u(x,y,t)\in \mathbb{R}^3$ is represented by
\[
v(\xi,\eta, t) \equiv u(x(\xi,\eta), y(\xi,\eta), t)
\]
in the computational domain.  To reduce notational complexity, we will use $u$ to denote the state in $\Omega$ and
$\Omega'$.  The chain rule can then be expressed as
\begin{equation}
\label{eq:chainrule}
\nonumber
	\begin{alignedat}{2}
		u_x &= \frac{1}{2} J^{-1} \left[(y_\eta u)_\xi  +  y_\eta u_\xi - (y_\xi u)_\eta - y_\xi u_\eta \right] \\
		u_y &= \frac{1}{2} J^{-1} \left[-(x_\eta u)_\xi  -  x_\eta u_\xi + (x_\xi u)_\eta + x_\xi u_\eta \right],
	\end{alignedat}
\end{equation}
where $J \equiv x_\xi y_\eta - x_\eta y_\xi$.  

We will recast Maxwell's equations into a form amenable to summation by parts for curvilinear domains.  
Following the approach in \cite{po:spps2}:
\begin{equation}
\label{eq:maxwell_pde_mod}
\nonumber
\begin{alignedat}{4}
	C u_t &= A\frac{1}{2} J^{-1} \left[(y_\eta u)_\xi  +  y_\eta u_\xi - (y_\xi u)_\eta - y_\xi u_\eta \right] \\
		  &+ B\frac{1}{2} J^{-1} \left[-(x_\eta u)_\xi  -  x_\eta u_\xi + (x_\xi u)_\eta + x_\xi u_\eta \right]
		\quad && \xi,\eta \in \Omega', \quad &&&t > 0 \\
	H &= g(x(\xi,\eta),y(\xi,\eta),t), && \xi,\eta \in \Gamma', &&& t \geq 0 \\
	u &= f(x(\xi,\eta),y(\xi,\eta)), && \xi,\eta \in \Omega', &&& t = 0.
\end{alignedat}
\end{equation}
Note that $J=\xi_x\eta_y - \xi_y\eta_x$ in \cite{po:spps2}.  This problem will now be solved in the 
four-block domain $\Omega = \cup_{ij}\Omega_{ij}$, cf. Sec.~\ref{sec:fourblockdiffop}.

\subsection{Inner product in state space $V$}
\label{sec:discretization}
The computational domain $\Omega'$ is made up of four equally sized subdomains $\Omega'_{ij}$, see 
Figure~\ref{fig:comp_domain}:
\[
\Omega' = \cup_{ij}\Omega'_{ij}, \quad N^{(i)}_j = N \Longrightarrow 
h^{(i)}_j = h = \frac{1}{N}, \quad i,j = 1,2.
\]
Let
\begin{equation}
\label{eq:phys_grid_coord}
x \equiv \left(x(ih,jh)\right), \ y \equiv \left(y(ih,jh)\right) \in \mathbb{R}^M,
\end{equation}
represent the coordinates of the physical grid.  The previous conventions for row access \eqref{eq:2dstate} 
and column access \eqref{eq:2daltstate} apply.  Define the metric coefficients
\begin{equation}
\label{eq:matric_coeff}
\nonumber
\begin{alignedat}{2}
x_\xi &\equiv \left(x_\xi(ih,jh)\right), \quad x_\eta &\equiv \left(x_\eta(ih,jh)\right) \\ 
y_\xi &\equiv \left(y_\xi(ih,jh)\right), \quad y_\eta &\equiv \left(y_\eta(ih,jh)\right)
\end{alignedat} \quad \in \mathbb{R}^M
\end{equation}
with the corresponding matrix versions:
\begin{equation}
\label{eq:matric_coeff_matrix}
\begin{alignedat}{2}
X_\xi &\equiv \text{diag}\left(x_{\xi,ij}\right)\otimes I, 
\quad X_\eta &\equiv \text{diag}\left(x_{\eta,ij}\right)\otimes I \\ 
Y_\xi &\equiv \text{diag}\left(y_{\xi,ij}\right)\otimes I, 
\quad Y_\eta &\equiv \text{diag}\left(y_{\eta,ij}\right) \otimes I 
\end{alignedat},
\quad I \in \mathbb{R}^{3\times 3}.
\end{equation}

Since the analytic scalar product satisfies
\[
(u,v) = \int_\Omega u^Tv\:dS = \int_{\Omega'}u^TvJ\:dS',
\]
it is natural to define the corresponding discrete scalar product $(\cdot,\cdot):V\times V\rightarrow\mathbb{R}$ as
\begin{equation}
\label{eq:2dinnerprod_curv}
\nonumber
(u,v) \equiv u^TJH v,
\end{equation}
where $u,v\in\mathbb{R}^M$, $M=(N_1+1)^2$, $N_1 = 2N$; $H=H_1\otimes H_1$ is defined by (\ref{eq:4blockinnerprodembed1}), or equivalently, (\ref{eq:4blockinnerprodembed2}), cf.~Lemma~\ref{lem:4blockinnerprodequiv}, where $H_1$ are one-dimensional
{\em diagonal} norms (\ref{eq:innerprodembed}) of size $(N_1+1)\times(N_1+1)$, since the corresponding one-dimensional norms 
$H^{(ij)}_{1,2}$ are assumed to be diagonal and identical for all subdomains $\Omega'_{ij}$.  Furthermore, 
$H^{(ij)}_{1,2}$ have been constructed such that the structural requirement \eqref{eq:norm} holds.  To be clear,
the diagonal elements of $H$ are of the form $h_ih_jI$, $I\in \mathbb{R}^{3\times 3}$, $0\leq i,j \leq N_1$.  The Jacobian matrix $J$ is defined as
\begin{equation}
\label{eq:jacobian}
J \equiv X_\xi Y_\eta - X_\eta Y_\xi.
\end{equation}
We can now define the discrete scalar product $(\cdot,\cdot)_C:V_C\times V_C\rightarrow\mathbb{R}$ 
that corresponds to \eqref{eq:maxwell_anal_stab}:
\begin{equation}
\label{eq:2dinnerprod_curv_maxwell}
\nonumber
(u,v)_C \equiv u^TCJH v, \quad C = I_M \otimes
\begin{pmatrix}
		\epsilon & 0 & 0 \\
		0 & \mu & 0 \\
		0 & 0 & \epsilon
\end{pmatrix}.
\end{equation}
The matrices $C, J, H$ commute since they are diagonal.  This completes the construction of the inner product space
for the state space of the Maxwell equations.

\subsection{Inner product in boundary state space $V_\Gamma$}
\label{sec:bdiscretization}
The arc length matrices $S^{(k)}\in\mathbb{R}^{3(N_1+1)\times 3(N_1+1)}$ will be needed in the boundary scalar products
$\langle\cdot ,\cdot\rangle_k$:
\begin{equation}
\label{eq:arc_length_matrix}
\nonumber
\begin{array}{ll}
S^{(1)} \equiv \text{diag}\left(|r_{\xi,i0}|\right)\otimes I,
&S^{(2)} \equiv \text{diag}\left(|r_{\eta,N_1j}|\right)\otimes I \\
S^{(3)} \equiv \text{diag}\left(|r_{\xi,iN_1}|\right)\otimes I, 
&S^{(4)} \equiv \text{diag}\left(|r_{\eta,0j}|\right)\otimes I,
\end{array}
\end{equation}
where
\[
|r_{\xi,ij}| \equiv \left[x^2_{\xi,ij} + y^2_{\xi,ij}\right]^{1/2}, \quad
|r_{\eta,ij}| \equiv \left[x^2_{\eta,ij} + y^2_{\eta,ij}\right]^{1/2},
\]
represent the length of the tangent (arc length) evaluated at each grid point.

The line integral \eqref{eq:maxwell_anal_stab} suggests that the discrete boundary scalar product 
$\langle\cdot,\cdot\rangle_+:V_{\Gamma_+}\times V_{\Gamma_+} \rightarrow \mathbb{R}$ 
be defined as \cite{po:spps2}
\begin{equation}
\label{eq:2dbinnerprod_curv}
\nonumber
\langle u,v\rangle_+ \equiv
\sum^{4}_{k=1} \langle u,v\rangle_k,
\end{equation}
where
\begin{equation}
\label{eq:2dbinnnerprod_curv_parts}
\nonumber
\begin{array}{ll}
	\langle u,v\rangle_1 \equiv u^T_0H_1S^{(1)}v_0, \quad 
   &\langle u,v\rangle_2 \equiv \left[u^{N_1}\right]^TH_1S^{(2)}v^{N_1}  \\
    \langle u,v\rangle_3 \equiv u^T_{N_1}H_1S^{(3)}v_{N_1}, \quad 
   &\langle u,v\rangle_4 \equiv \left[u^0\right]^TH_1S^{(4)}v^0.
\end{array}
\end{equation}
Hence, cf.~\eqref{eq:bstateinnerprodplus},
\begin{equation}
\label{eq:bstateinnerprodplus_curv}
\langle u, v \rangle_+  \equiv u^T_{\Gamma_+}H^{(+)}_\Gamma S^{(+)}v^T_{\Gamma_+},
\end{equation}
where $S^{(+)} \in \mathbb{R}^{12(N_1+1)\times 12(N_1+1)}$ is given by
\begin{equation}
\label{eq:sgamma+_curv}
\nonumber
S^{(+)} \equiv
\begin{pmatrix}
S^{(1)} \\
&S^{(2)} \\
&&J_{N_1+1}S^{(3)}J_{N_1+1} \\
&&&J_{N_1+1}S^{(4)}J_{N_1+1}
\end{pmatrix} ;
\end{equation}
$J_{N_1+1}\in\mathbb{R}^{3(N_1+1)\times 3(N_1+1)}$ is the block anti-diagonal permutation matrix of \eqref{eq:hgamma+}.
Since $H$ satisfies \eqref{eq:norm}, it follows that $J_{N_1+1}H_1J_{N_1+1}=H_1$, whence
\begin{equation}
\label{eq:hgamma+_curv_simple}
\nonumber
H^{(+)}_\Gamma  =
\begin{pmatrix}
H_1 \\
&H_1 \\
&&H_1 \\
&&&H_1
\end{pmatrix} \in \mathbb{R}^{12(N_1+1)\times 12(N_1+1)}.
\end{equation}

\begin{rem}
\label{rem:hgamma+_curv}
From a summation-by-parts point of view, it does not matter if one uses $u_{\Gamma_{3,4}}$ \eqref{eq:2dbstate},
or if one employs the corresponding vector components $u_N, u^0$.  The arc lengths $S_{3,4}$ must still be 
properly ordered, however. \hfill$\Box$
\end{rem}
We are now in a position to define $\langle\cdot, \cdot\rangle: V_\Gamma \times V_\Gamma \rightarrow \mathbb{R}$:
\begin{equation}
\label{eq:bstateinnerprod_curv}
\langle u, v \rangle \equiv \langle Eu, Ev \rangle_+ = u^T_\Gamma E^TH^{(+)}_\Gamma S^{(+)}Ev_\Gamma;
\end{equation}
$u_\Gamma,v_\Gamma \in \mathbb{R}^{12N_1}$ are grid vectors on $\Gamma = \cup_i\Gamma_i$; the embedding
$E\in\mathbb{R}^{12(N_1+1)\times 12N_1}$ is defined by \eqref{eq:eembedboundary}.  Hence,  
$E:V_\Gamma \rightarrow V_{\Gamma_+}$ is an isometric isomorphism.   Furthermore, by  \eqref{eq:bstateinnerprodplus_curv}
and \eqref{eq:bstateinnerprod_curv}:
\[
\langle u, v \rangle = \langle u, v \rangle_+.
\]
This expression defines the arc length operator $S:V_\Gamma \rightarrow V_\Gamma$ through the relation
$H_\Gamma S \equiv E^TH^{(+)}_\Gamma S^{(+)}E$, where $H_\Gamma$ is defined via \eqref{eq:bstateinnerprod}.  Thus,
\begin{equation}
\label{eq:sgamma_curv}
S \equiv H^{-1}_\Gamma E^TH^{(+)}_\Gamma S^{(+)}E = E^+S^{(+)}E.
\end{equation}
If $x=\xi$ and $y=\eta$, we recover \eqref{eq:bstateinnerprod}.

Summation by parts on curvilinear grids requires explicit knowledge of the outward unit normals 
$n^{(k)} = (n^{(k)}_x\ n^{(k)}_y)^T$ for each boundary segment $\Gamma_k$:
\begin{equation}
\label{eq:out_normal}
\nonumber
\begin{alignedat}{1}
 n^{(1)}_i = \begin{pmatrix} n^{(1)}_{x,i}&n^{(1)}_{y,i}\end{pmatrix}^T
 &\equiv \begin{pmatrix} y_{\xi,i0} &-x_{\xi,i0}\end{pmatrix}^T/|r_{\xi,i0}| \\
 n^{(2)}_j = \begin{pmatrix} n^{(2)}_{x,j}&n^{(2)}_{y,j}\end{pmatrix}^T
 &\equiv \begin{pmatrix} y_{\eta,N_1j} &-x_{\eta,N_1j}\end{pmatrix}^T/|r_{\eta,N_1j}| \\
 n^{(3)}_i = \begin{pmatrix} n^{(3)}_{x,i}&n^{(3)}_{y,i}\end{pmatrix}^T 
 &\equiv \begin{pmatrix} -y_{\xi,iN_1} &x_{\xi,iN_1}\end{pmatrix}^T/|r_{\xi,iN_1}| \\
 n^{(4)}_j = \begin{pmatrix} n^{(4)}_{x,j}&n^{(4)}_{y,j}\end{pmatrix}^T
 &\equiv \begin{pmatrix} -y_{\eta,0j} &x_{\eta,0j}\end{pmatrix}^T/|r_{\eta,0j}|.
\end{alignedat}
\end{equation}
Next, we construct the corresponding outward normal matrices $N^{(i)}_x, N^{(i)}_y$:
\begin{equation}
\label{eq:out_normal_matrix}
\nonumber
N^{(k)}_x \equiv \text{diag}\left(n^{(k)}_{x,i}\right)\otimes I, \quad
N^{(k)}_y \equiv \text{diag}\left(n^{(k)}_{y,i}\right)\otimes I, \quad I\in \mathbb{R}^{3\times 3}.
\end{equation}
Completely analogous to the arc length operator $S^{(+)}: V_{\Gamma_+} \rightarrow V_{\Gamma_+}$, we define
$N^{(+)}_x, N^{(+)}_y: V_{\Gamma_+} \rightarrow V_{\Gamma_+}$ representing the outward normal for each boundary point:
\begin{equation}
\label{eq:n+xy_curv}
\nonumber
N^{(+)}_{x,y} =
\begin{pmatrix}
N^{(1)}_{x,y} \\
&N^{(2)}_{x,y} \\
&&J_{N_1+1}N^{(3)}_{x,y}J_{N_1+1} \\
&&&J_{N_1+1}N^{(4)}_{x,y}J_{N_1+1}
\end{pmatrix}.
\end{equation}
Define the outward normal operators $N_x, N_y: V_\Gamma \rightarrow V_\Gamma$ as
\begin{equation}
\label{eq:nxy_curv}
\begin{alignedat}{1}
N_x &\equiv S^{-1}E^+S^{(+)} N^{(+)}_xE \\
N_y &\equiv S^{-1}E^+S^{(+)} N^{(+)}_yE.
\end{alignedat}
\end{equation}
If 
\begin{equation}
\label{eq:sgamma_embed}
\begin{alignedat}{1}
s^{(1)}_0 &= s^{(4)}_N \\
s^{(2)}_0 &= s^{(1)}_N \\
s^{(3)}_0 &= s^{(2)}_N \\
s^{(4)}_0 &= s^{(3)}_N,
\end{alignedat}
\end{equation}
it follows that $E^+S^{(+)} = SE^+$, where  $S$ is defined by \eqref{eq:sgamma_curv}.
Conversely, if $E^+S^{(+)} = SE^+$, then \eqref{eq:sgamma_embed} is implied.  Thus,
if \eqref{eq:sgamma_embed} holds, then \eqref{eq:nxy_curv} simplifies to
\begin{equation}
\label{eq:nxy_curv_simplified}
\nonumber
\begin{alignedat}{1}
N_x &= E^+N^{(+)}_xE \\
N_y &= E^+N^{(+)}_yE.
\end{alignedat}
\end{equation}

\begin{rem}
\label{rem:sgamma_embed}
If one interprets $S$ and $S^{(+)}$ as vectors $s$ and $s^{(+)}$ instead of operators, then 
\eqref{eq:sgamma_embed} simply states that $s^{(+)}=s^{(e)}=Es$, i.~e., $s^{(+)}$ is the embedding 
of $s$, where the arc length vector $s$ is uniquely defined for each boundary point in 
$\Gamma=\cup_i\Gamma_i$. In our case, \eqref{eq:sgamma_embed} implies that $|r_\xi| = |r_\eta|$ 
at the corners of $\Gamma$.  Note that this condition is trivially satisfied for the 
standard Cartesian grid, where $x=\xi, y=\eta$.  \hfill$\Box$
\end{rem}

\subsection{Summation by parts in curvilinear domains}
\label{sec:sbp_curv}
In light of the chain rule, it makes sense to define $D_{x,y}:V \rightarrow V$:
\begin{equation}
\label{se:chainrule_disc}
\nonumber
	\begin{alignedat}{2}
		D_x &\equiv \frac{1}{2} J^{-1} (Y_\eta D_\xi  + D_\xi Y_\eta - Y_\xi D_\eta - D_\eta Y_\xi) \\
		D_y &\equiv \frac{1}{2} J^{-1} (X_\xi D_\eta  + D_\eta X_\xi - X_\eta D_\xi - D_\xi X_\eta),
	\end{alignedat}
\end{equation}
where $X_\xi$, $X_\eta$, $Y_\xi$, $Y_\eta$ and $J$ are defined according to \eqref{eq:matric_coeff_matrix} 
and \eqref{eq:jacobian}.
$D_{\xi,\eta}$ is shorthand for $D_{\xi,\eta}\otimes I$.  They correspond to $D_{x,y}$ of 
Section~\ref{sec:fourblockdiffop} and share all structural results with $D_{x,y}$.  Straightforward but 
somewhat tedious computations show that
\begin{equation}
	\label{eq:maxwell_sbp_props}
	\begin{alignedat}{2}
		(u,D_x v) + (D_x u,v) &= \sum^4_{k=1} \langle u,N^{(k)}_x v\rangle_k
		                       = \langle u,N_x v\rangle, \quad &&\forall u,v \in V \\
		(u,D_y v) + (D_y u,v) &= \sum^4_{k=1} \langle u,N^{(k)}_y v\rangle_k
		                       = \langle u,N_y v\rangle, \quad&&\forall u,v \in V.
	\end{alignedat}
\end{equation}

In practical situations it can happen that the grid vectors $x$ and $y$ are given, but no analytic expressions
for $x_\xi$, $y_\xi$, $x_\eta$ and $y_\eta$ are known.  In such cases one can compute 
the metric coefficients numerically from \eqref{eq:phys_grid_coord}, cf.~\cite{ALUND2019209}:
\begin{equation}
\label{eq:matric_coeff_disc}
\nonumber
\begin{alignedat}{2}
x_\xi &\equiv D_\xi x, \quad x_\eta &\equiv D_\eta x \\ 
y_\xi &\equiv D_\xi y, \quad y_\eta &\equiv D_\eta y
\end{alignedat} \quad \in \mathbb{R}^M.
\end{equation}
All other definitions remain unchanged.  The summation-by-parts rule \eqref{eq:maxwell_sbp_props} still holds.
We have used this approach to obtain the results in Section~\ref{sec:convergence}

\subsection{Stability of the semidiscrete Maxwell's equations}
\label{sec:maxwell_semi}
Let
\begin{equation}
\label{eq:abc}
\nonumber
\begin{alignedat}{1}
	A &\equiv I \otimes \begin{pmatrix}
		0 & 0 & 0 \\
		0 & 0 & -1 \\
		0 & -1 & 0
	\end{pmatrix} \\
	B &\equiv I \otimes \begin{pmatrix}
		0 & 1 & 0 \\
		1 & 0 & 0 \\
		0 & 0 & 0
	\end{pmatrix} \\
	C &\equiv I \otimes \begin{pmatrix}
		\epsilon & 0 & 0 \\
		0 & \mu & 0 \\
		0 & 0 & \epsilon
	\end{pmatrix}
	\end{alignedat} \quad I\in \mathbb{R}^{M\times M},
\end{equation}
where $M = (N_1+1)^2, N_1=2N$ is the total number of grid points.  The spatial discretization 
of \eqref{eq:maxwell_pde} is given by
\[
	\begin{alignedat}{2}
		Cv_t &= AD_x v + BD_y v, \quad &&t > 0 \\
		L v &= g(t), && t \geq 0 \\
		v &= v_0, && t = 0,
	\end{alignedat}
\]
where the discrete boundary operator $L:V\rightarrow V_{\Gamma_+}$ is defined as in 
\eqref{eq:2dcharbcdiscrete}, \eqref{eq:2dcharbcdiscrete1}:
\[
L \equiv 
\begin{pmatrix}
L_1 \\
L_2 \\
L_3 \\
L_4
\end{pmatrix},\qquad
\begin{alignedat}{1}
L_1 &= I_0     \otimes I       \otimes L_0 \\
L_2 &= I       \otimes I_{N_1} \otimes L_0 \\
L_3 &= I_{N_1} \otimes I       \otimes L_0 \\
L_4 &= I       \otimes I_0     \otimes L_0
\end{alignedat}, \qquad
L_0 = 
\begin{pmatrix}
0 &0 &0 \\
0 &1 &0  \\
0 &0 &0
\end{pmatrix};
\]
$I$ is the $(N_1+1) \times (N_1+1)$ identity matrix, $I_{0,N_1}$ are the first and last and rows of $I$. 
The boundary conditions $Lv = g(t)$ are imposed using the simplified projection method \eqref{eq:simplesemiibvp}. 
The resulting scheme is given by
\begin{equation}
	\label{eq:maxwell_ode_system}
	 \begin{alignedat}{2}
	 	w_t &= Q w + G(t), \quad && t > 0 \\
	 	w &= P v_0, && t = 0,
	 \end{alignedat}
\end{equation}
where 
\[
	\begin{alignedat}{2}
		Q &= PC^{-1}\left(AD_x + BD_y\right)P \\
		G(t) &= PC^{-1}\left(AD_x + BD_y\right)L^+g(t).
	\end{alignedat}
\]
The approximate solution $v$ is obtained using \eqref{eq:vdef}:
\[
	v = w + L^+g(t),
\]
where we also used $Pw=w$.  

We notice that $H_c\equiv CJH$ is symmetric positive definite (SPD) since all matrices are diagonal
with positive diagonal elements.  Hence, 
\[
LH_c = \bar{H}L
\]
for some $\bar{H}>0$, cf.~Sec.~\ref{sec:simpleprojrevisited}, and so
\[
L^+ = L^T\left(LL^T\right)^+.
\]
All rows of $L$ are orthogonal except for the rows corresponding to the four corners, in which case the same boundary
condition is enforced twice.  Removing the extraneous boundary conditions (superfluous rows of $L$) leads
to a very simple expression for $P$:
\[
P = I - L^TL.
\]
This projection operator $P$ is self-adjoint with respect to $(\cdot,\cdot)_C$

To prove stability of the ODE system \eqref{eq:maxwell_ode_system}, it is sufficient to 
consider the homogeneous problem, i.~e., $G(t)=0$.  Thus,
\begin{equation}
\label{eq:energy_est_maxwell}
\begin{alignedat}{2}
	(w,w_t)_C &= (w,PC^{-1}\left(AD_x + BD_y\right)Pw)_C \\
	&= (u,\left(AD_x + BD_y\right)u),
\end{alignedat}
\end{equation}
where have defined the temporary variable $u\equiv Pw$.  Substituting $v\rightarrow Au$ and 
$v\rightarrow Bu$ in \eqref{eq:maxwell_sbp_props}:
\[
\begin{alignedat}{1}
(u,AD_xu) + (AD_xu,u) &= \langle u,AN_x u\rangle \\
(u,BD_yu) + (BD_yu,u) &= \langle u,BN_y u\rangle.
\end{alignedat}
\]
Note that $A$ commutes with $D_x$, $H$ and $J$ 
(similar relations are true for $B$ as well).  Hence, by \eqref{eq:energy_est_maxwell}:
\[
\frac{d}{dt}\|w\|^2_C = \langle u,\left[AN_x + BN_y \right]u\rangle.
\]
The right member is made up of exactly the same kind of terms as the integrand of \eqref{eq:maxwell_anal_stab}.
Since $Lu=LPw=0$, it follows that the boundary terms vanish identically, whence
\[
\frac{d}{dt}\|w\|^2_C = 0,
\]
which proves stability of \eqref{eq:maxwell_ode_system}.

\begin{rem}
\label{rem:maxwell_stab}
It is not necessary to require that $A$ and $B$ be constant in space to prove stability, 
see \cite{po:spps2} for details. \hfill$\Box$
\end{rem}

\begin{figure}
\begin{subfigure}{.5\textwidth}
	\centering
	\includegraphics[width=0.9\textwidth]{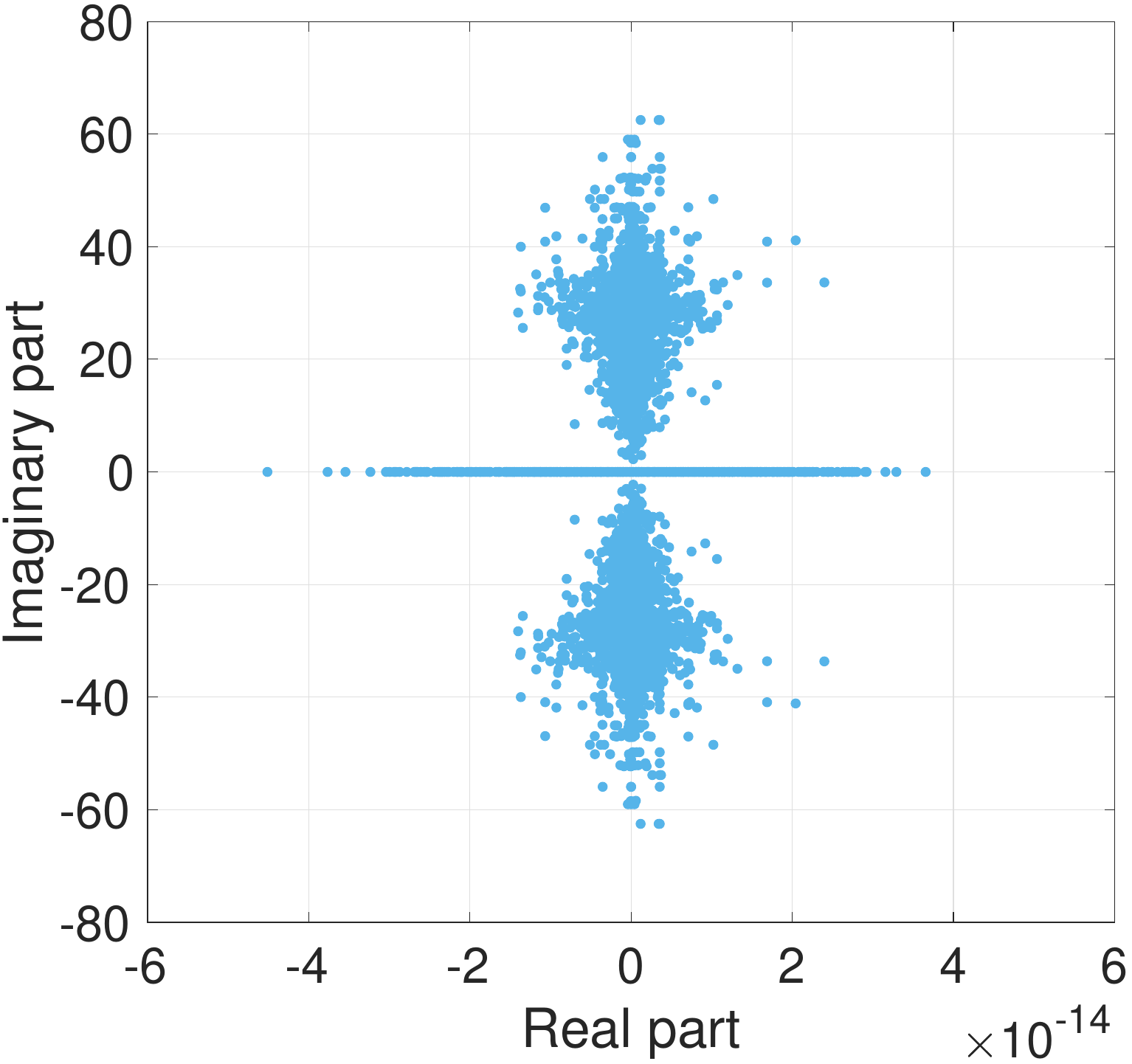}
	\caption{4th order operators.}
\end{subfigure}
\begin{subfigure}{.5\textwidth}
	\centering
	\includegraphics[width=0.9\textwidth]{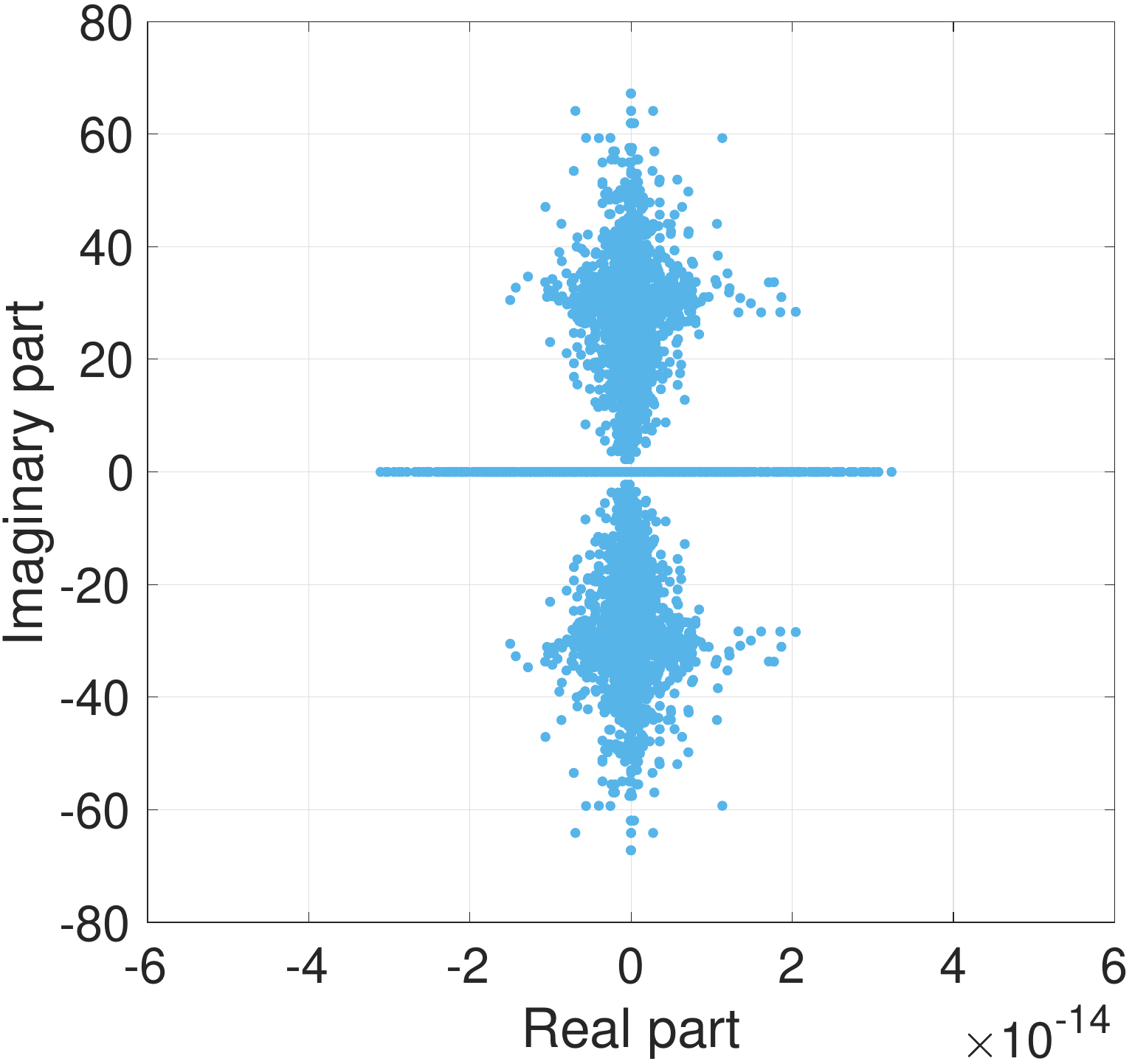}
	\caption{6th order operators.}
\end{subfigure}
\caption{Spectrum of $Q$ for 4th and 6th order SBP operators.}
\label{fig: maxwell_spectrum}
\end{figure}

The stability of \eqref{eq:maxwell_ode_system} can also be verified numerically by studying the 
eigenvalues of $Q$. In Figure \ref{fig: maxwell_spectrum} the real and imaginary eigenvalues are 
plotted with $41 \times 41$ points in each dimension, corresponding to a total of 5,043 degrees 
of freedom. Clearly, $Q$ has only imaginary eigenvalues, which again indicates stability and 
energy conservation.

\subsection{Convergence results}
\label{sec:convergence}
To evaluate the accuracy properties of the scheme, we use an analytical solution given by
\begin{equation}
	\label{eq:maxwell_analytical}
	\nonumber
	\begin{alignedat}{2}
	E_x(x,y,t) &= -\frac{4}{5}\cos(3x + 4y - 5t) \\
	H(x,y,t) &= \cos(3x + 4y - 5t) \\
	E_y(x,y,t) &= \frac{3}{5}\cos(3x + 4y - 5t),
	\end{alignedat}
\end{equation}
which corresponds to choosing $\epsilon = 1/5$ and $\mu=5$.  The boundary and initial data 
are obtained from the above expressions. The time stepping is done using the classical 
fourth order explicit Runge-Kutta method, with the time step given by 
$\Delta t = \frac{h_{\text{\tt min}}}{10}$ where $h_{\text{\tt min}}$ is the smallest 
spatial step size. This choice of time step ensures that the temporal error is negligible 
in relation to the spatial error.  The error $e_{p,M}$ is measured in the discrete $L^2$-norm:
\[
	e_{p,M} \equiv \|v - v_{\text{\tt exact}}\| = \sqrt{(v - v_{\text{\tt exact}})^T JH (v - v_{\text{\tt exact}})}, \quad p = 2,4,6,
\]
at $t = 1$, and the convergence as
\[
	q_p \equiv \frac{\log \frac{e_{p,M_1}}{e_{p,M_2}}}{\log\left(\frac{M_2}{M_1}\right)^{1/2}},
\]
where the subscript $p$ indicates the interior accuracy of the SBP operators.

In Table \ref{tabl:errconv}, the error and convergence for varying grid resolutions is 
presented for SBP operators of interior orders 2, 4, and 6. The boundary accuracies of the 
SBP operators are 1, 2, and 3. We see that for all operators the global convergence rate 
is one higher than the boundary accuracy, which is in accordance with the theoretical 
convergence analysis found in \cite{SVARD2021110020}. 

\begin{table}
\centering
\caption{Error (in base 10 logarithm) and convergence of Maxwell simulation with interface 
conditions imposed using the embedding method and 2nd, 4th, and 6th order SBP operators.}
\label{tabl:errconv}
\begin{tabular}{|c||c|c||c|c||c|c|} 
\hline
$N$ & $\log_{10}(e_{2,M})$ & $q_2$ & $\log_{10}(e_{4,M})$ & $q_4$ & $\log_{10}(_{6,M})$ & $q_6$   \\ 
\hline
40  & -1.45 & - 	 & -2.15 & - & -2.33 & - 	  \\ 
\hline
120  & -2.39 & 1.99 & -3.54 & 2.93 & -4.15 & 3.85  \\ 
\hline
200 & -2.83 & 1.98 & -4.20 & 2.97 & -5.08 & 4.18  \\ 
\hline
280 & -3.11 & 1.98 & -4.63 & 2.98 & -5.69 & 4.16   \\ 
\hline
360 & -3.33 & 1.98 & -4.96 & 2.98 & -6.14 & 4.13  \\ 
\hline
440 & -3.50 & 1.98 & -5.22 & 2.98 & -6.49 & 4.10  \\ 
\hline
520 & -3.65 & 1.98 & -5.43 & 2.98 & -6.79 & 4.08  \\ 
\hline
600 & -3.77 & 1.98 & -5.62 & 2.98 & -7.04 & 4.06  \\ 
\hline
\end{tabular}
\end{table}

\section{Discussion and conclusions}
\label{sec:disc}
In the present work, we have taken a vector space centric approach when discussing
summation-by-parts operators and the implementation of analytic boundary conditions.
The difference operators and boundary operators are regarded as mappings 
$D:V\rightarrow V$ and $L:V\rightarrow V_\Gamma$.  The inner product of the state
space $V$ is given by the summation-by-parts norms $H$; the scalar product of the
boundary state $V_\Gamma$ is implicitly determined by $H$ via summation by parts.
With these definitions in place, it is possible to give a formal definition of the
adjoint operators $D^*$ and $L^*$.

We have also shown how the pseudoinverse of the boundary operator can be used to 
generalize the implementation of boundary conditions as a projection:
\[
P = I - L^+L.
\]
The above expression is valid for {\em any} linear boundary operator $L$ regardless of 
rank.  This facilitates theoretical analysis in the presence of corners, which 
potentially may cause rank deficient, or near rank deficient, boundary operators.
The projection $P$ is not uniquely determined in general.  We used this fact to our 
advantage to simplify the expression for $L^+$ as much as possible, see 
Sec.~\ref{sec:psimplified}.  It was shown that one can always choose $H_\Gamma=I$ when
constructing the boundary projection.  Under certain circumstances, the 
boundary projection is completely independent of $H_\Gamma$ {\em and} $H$, 
cf.~\eqref{eq:psimplified2}, thus extending the conclusions of \cite{po:spps1}
to the general, possibly rank deficient, case.  The pseudoinverse provides a concise
way of representing the boundary data $Lv = g$ as a state vector defined on $\Omega$:
\[
v = w + L^+g,
\]
where $w$ solves the simplified semidiscrete equations \eqref{eq:simplesemiibvp}.

The embedding operator $E$ introduced in Sec.~\ref{sec:embedding} provides a convenient
mechanism for extending summation-by-parts operators defined in multiple domains
to a single operator defined in the union of the individual domains.  Given two
difference operators $D^{(i)}:V_i \rightarrow V_i$, the resulting multidomain
operator $D:V \rightarrow V$can be expressed as
\[
D = E^+D^{(+)}E,
\]
where the embedding $E:V \rightarrow V_+$ is defined in  \eqref{eq:epart} - 
\eqref{eq:e2embed}; $E^+$ is the pseudoinverse of $E$.  The new operator $D$
will inherit all properties of its constituent operators $D^{(1)}$ and $D^{(2)}$,
most notably summation by parts and accuracy.  We have also demonstrated how to
construct embedding operators in two space dimensions.  The results from the one-dimensional
theory carry over word for word.  In summary, given $D^{(i)}, H^{(i)}, V_i$, 
the construction of the multidomain difference operator $D$ follows the same pattern 
regardless of dimensionality:

\begin{enumerate}
\item Form the inner product $H^{(+)}$ in the extended state space $V_+$ using the 
inner products $H^{(i)}$ of the existing state spaces $V_i$.
\item Construct the embedding  $E:V \rightarrow V_+$.
\item Define the scalar product in $V$ as $H = E^TH^{(+)}E$.
\item Define the extended difference operator $D^{(+)}:V_+ \rightarrow V_+$.
\item Let $D:V \rightarrow V$ be defined as $D = E^+D^{(+)}E$.
\end{enumerate}

To illustrate the theory, we have implemented a test scenario involving the
two-dimensional Maxwell's equations on four curvilinear domains, which are joined
together using two-dimensional embedding operators $E_x,E_y$.  The boundary conditions
are implemented using a projection operator.  The rate of convergence
as measured for 2nd, 4th, and 6th-order accurate methods agrees very well with the 
theoretical convergence analysis.

The overarching goal of this work is to lay the theoretical foundation for
a "plug \& play" methodology in which stability of the semidiscrete problem follows more or less 
directly from well-posedness of the analytic problem.  As soon as the analytic boundary 
conditions are known, there is recipe for how to construct the corresponding projection,
such that it leads to a stable approximation.  Similarly, given the norms 
and difference operators $H^{(i)}, D^{(i)}$ that satisfy summation by parts on their 
respective domains $\Omega$, we can construct 
\[
H = E^TH^{(+)}E,\quad D = E^+D^{(+)}E,
\]
such that $D$ satisfies summation by parts on $\Omega = \cup_i\Omega_i$ with respect
to $H$.  In a coming study we will extend this to more general domains and look
deeper into algorithmic aspects, such as blocking and memory efficiency.  The
latter is important when solving large problems in multithreaded compute environments.

\bibliography{references} 
\bibliographystyle{plain}

\end{document}